\documentclass[letterpaper]{amsart}
\pdfoutput=1
\usepackage[T1]{fontenc}
\usepackage[margin=1in]{geometry}
\usepackage{amssymb, latexsym, amsmath, verbatim, amsthm,}
\usepackage{xfrac}
\usepackage{enumitem}
\usepackage{mathtools}
\usepackage{graphicx}
\usepackage{tikz}
\usetikzlibrary{calc,positioning}
\usepackage{tikz-cd}
\usepackage{microtype}
\usepackage{hyperref}
\allowdisplaybreaks

\title{Hierarchical hyperbolicity of Coxeter groups}
 
\author{Christopher H. Cashen} 
\address{
  TU Wien\\ Institute of Discrete Mathematics and Geometry\\ Wiedener
  Hauptstrasse~8-10\\ 1040 Vienna \\ Austria\\ \href{https://orcid.org/0000-0002-6340-469X}{\includegraphics[scale=.75]{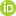}
  0000-0002-6340-469X}}
\email{christopher.cashen@tuwien.ac.at}
\author{Michelle Chu}
\address{University of Minnesota \\ School of Mathematics \\ 206 Church St. SE \\
Minneapolis, MN 55455}
\email{mchu@umn.edu}
\author{Jing Tao} 
\address{University of Oklahoma \\
    Department of Mathematics\\
    601 Elm Avenue Room 423\\Norman, OK 73019}
\email{jing@ou.edu}

\date{\today}
\thanks{Research supported by the Austrian Science
Fund (FWF) \href{https://doi.org/10.55776/PAT7799924}{10.55776/PAT7799924}, NSF DMS-2441034, and NSF DMS-2609746.}

\hypersetup{
    pdftitle={Hierarchical hyperbolicity of Coxeter groups},    
    pdfauthor={Christopher H. Cashen, Michelle Chu, Jing Tao},     
    pdfkeywords={Coxeter group, hierarchically hyperbolic group, relative
  hierarchical hyperbolicity, coarse median}, 
    colorlinks=true,       
    linkcolor=black,          
    citecolor=black,        
    filecolor=black,      
    urlcolor=black           
}

\theoremstyle{plain}
\newtheorem{theorem}{Theorem}[section]

\newtheorem{lemma}{Lemma}[section]
\newtheorem{proposition}{Proposition}[section]
\newtheorem{corollary}{Corollary}[section]

\theoremstyle{remark}

\newtheorem{remark}{Remark}[section]

\theoremstyle{definition}
\newtheorem{definition}{Definition}[section]
\newtheorem{example}{Example}[section]

\def\makeautorefname#1#2{\expandafter\def\csname#1autorefname\endcsname{#2}}
\let\fullref\autoref

\makeautorefname{theorem}{Theorem} 
\makeautorefname{lemma}{Lemma} 
\makeautorefname{proposition}{Proposition} 
\makeautorefname{corollary}{Corollary} 
\makeautorefname{definition}{Definition}
\makeautorefname{example}{Example}
\makeautorefname{section}{Section}
\makeautorefname{subsection}{Section}
\makeautorefname{subsubsection}{Section}
\makeautorefname{claim}{Claim}
\makeautorefname{question}{Question}
\makeautorefname{conjecture}{Conjecture}
\makeautorefname{remark}{Remark}

\makeatletter 
\let\c@lemma=\c@theorem 
\makeatother
\makeatletter 
\let\c@proposition=\c@theorem 
\makeatother
\makeatletter 
\let\c@corollary=\c@theorem 
\makeatother
\makeatletter 
\let\c@definition=\c@theorem 
\makeatother
\makeatletter 
\let\c@example=\c@theorem 
\makeatother
\makeatletter 
\let\c@conjecture=\c@theorem 
\makeatother
\makeatletter 
\let\c@remark=\c@theorem 
\makeatother
\makeatletter
\@addtoreset{claim}{proposition}
\makeatother
\makeatletter
\@addtoreset{claim}{lemma}
\makeatother
\mathtoolsset{centercolon} 

\newcommand{\interior}[1]{\mathring{#1}} 
\newcommand{\act}{\curvearrowright} 
\newcommand{\from}{\colon\thinspace} 
\newcommand{\onto}{\twoheadrightarrow} 
\newcommand{\into}{\hookrightarrow} 

\renewcommand{\setminus}{\smallsetminus}

\usepackage{bbm}

\renewcommand{\vec}[1]{{\bf #1}}
\newcommand{\Coxetermatrix}{M}
\newcommand{\Titsform}{B _{\! \Coxetermatrix}}
\newcommand{\canonicalrep}{\rho}

\DeclareFontFamily{U}{mathx}{}
\DeclareFontShape{U}{mathx}{m}{n}{<-> mathx10}{}
\DeclareSymbolFont{mathx}{U}{mathx}{m}{n}
\DeclareMathAccent{\widehat}{0}{mathx}{"70}
\DeclareMathAccent{\widecheck}{0}{mathx}{"71}

\newcommand{\Davis}{\Sigma}
\newcommand{\NR}{X_{\mathrm{N\! R}}}
\DeclareMathOperator{\stab}{Stab} 
\DeclareMathOperator{\pc}{Pc} 
\DeclareMathOperator{\GL}{GL} 
\DeclareMathOperator{\rank}{rank}
\newcommand{\transverse}{\pitchfork}
\DeclareMathOperator{\diam}{diam} 

\newcommand{\coneoff}{\widehat{\Davis}}

\newcommand{\dcomb}{d_1}

\newcommand{\ceq}{\asymp}

\newcommand{\centralizer}{Z}

\newcommand{\normalizer}{N}

\newcommand{\caprace}{\ell}

\newcommand{\gate}{\mathfrak{p}}
\newcommand{\walls}{\mathcal{M}}

\newcommand{\fibers}{\stab^{-1}}
\DeclareMathOperator{\sym}{Sym}
\DeclareMathOperator{\alt}{Alt}
\DeclareMathOperator{\orthog}{O}

\newcommand{\HHindex}{\mathfrak S}

\newcommand{\ww}{%
  \mathord{\mathchoice
    {\wwaux{.2}}%
    {\wwaux{.2}}%
    {\wwaux{.14}}%
    {\wwaux{.10}}}}
    
\newcommand{\wwaux}[1]{%
  \tikz[baseline=-0.1ex,scale=#1,line width=.45pt]
  \draw[thick] (0,1) -- (0,0) -- (.5,.65) -- (1,0) -- (1,1);}

\renewcommand{\top}{\ww}

\newcommand{\nested}{\sqsubseteq}
\makeatletter
\newcommand{\nestedneq}{\mathrel{\mathpalette\nestedneq@{}}}
\newcommand{\nestedneq@}[2]{%
  \ooalign{%
    $\m@th#1\sqsubseteq$\cr
    \hidewidth
    \raisebox{-0.38ex}{\scalebox{0.55}[0.55]{$\m@th#1/$}}%
    \hidewidth\cr
  }%
}
\makeatother
\newcommand{\orthogonal}{\perp}
\newcommand{\orthogonalto}[1]{\HHindex_{#1}^\orthogonal}

\newcommand{\coordinate}[1]{\mathcal{C}#1}

\newcommand{\domaindist}[3]{\delta_{#1}(#2,#3)}

\newcommand{\proj}[1]{\pi_{#1}}

\newcommand{\relprojection}[2]{\rho^{#1}_{#2}}

\newcommand{\nestedin}[1]{\HHindex_{#1}}

\newcounter{hhsaxiom}

\newenvironment{hhsaxioms}
{%
  \begin{list}{\textup{(\arabic{hhsaxiom})}}{%
    \usecounter{hhsaxiom}%
    \setlength{\leftmargin}{2.5em}%
    \setlength{\labelwidth}{2em}%
    \setlength{\labelsep}{0.5em}%
  }%
}
{%
  \end{list}%
}

\newcommand{\axiom}[1]{%
  \item\textbf{#1.}\ %
}

\begin{document}
\begin{abstract}
We show that for every finite rank Coxeter system, the Coxeter group
is hierarchically hyperbolic relative to its nonhyperoctahedral
irreducible higher rank affine parabolic subgroups.
This is sharp---the Coxeter group admits a hierarchically hyperbolic
group structure if and only if it has no such parabolic subgroups, and
we give a complete diagrammatic characterization of when this is true.
Moreover, this condition is equivalent to the existence of a coarsely
equivariant coarse median on the Coxeter group. 

We also show that the Coxeter group is virtually
a colorable hierarchically hyperbolic group. This is achieved by passing to a congruence subgroup of the Tits representation.

As consequences, the Davis complex is a hierarchically hyperbolic
space and the group admits a coarse median of rank
equal to its algebraic rank. Furthermore, the Coxeter group is quasimedian
quasiisometric to a finite dimensional CAT(0) cube complex and is
virtually coarsely injective. 

Two appendices prove results of independent interest that are not
Coxeter-specific. 
The first promotes a relative HHS structure with HHS peripheral coordinate
spaces to an HHS structure. The second shows that a virtually
Abelian group admits a coarsely equivariant coarse median precisely
when it is hyperoctahedral.
\end{abstract}

\maketitle



\section{Introduction}\label{sec:intro}
    Coxeter groups have a rich geometry arising from their wall structure
and the structure of their parabolic subgroups.
In a previous paper \cite{CasChuTao} we showed that hyperbolic aspects
of this geometry can be understood in close analogy to those of the
mapping class group of a hyperbolic surface.
By coning off the nonhyperbolic regions in the Davis complex
$\Davis=\Davis(W,S)$ of a finite rank Coxeter system $(W,S)$, we
produced a hyperbolic space $\coneoff$ playing the role of the
curve graph of the surface: it recognizes the Morse quasigeodesics and stable
subspaces of the Davis complex, and provides the largest
acylindrically hyperbolic structure on the Coxeter group.
However, coning off the nonhyperbolic regions collapses them to
bounded diameter, effectively discarding the geometry of unbounded
regions in $\Davis$.
The purpose of this paper is to recover that missing geometry, in a
way that further strengthens the mapping class group analogy.

For mapping class groups, a vertex in the curve graph corresponds to a
curve in the surface.
Cutting the surface along that curve gives
a subsurface whose components have strictly lower complexity
and are equipped with their own mapping class groups and curve graphs.
This leads to an inductive, hierarchical description of mapping class
group geometry that uses subsurface projection to recover the geometry
lost in passing to the curve graph \cite{MasMin99, MasMin00}. 
For a Coxeter group, a cone vertex of $\coneoff$ is associated
to a parabolic subcomplex of $\Davis$.
This is the Davis complex of a parabolic subgroup of $(W,S)$, whose
irreducible factors are Coxeter groups of strictly lower rank, suggesting an analogous
hierarchical induction to recover the
geometry lost in passing to $\coneoff$.
The notion of a hierarchically hyperbolic space formalizes these hierarchical constructions.

A \emph{hierarchically hyperbolic space} (HHS) is a quasigeodesic
metric space described by a system of hyperbolic coordinate spaces
that are indexed by hierarchically related ``domains''.
A set of axioms governs the interactions between domains,
their coordinate spaces, and projections to the coordinate spaces, in
such a way that consistent coordinate data can be coarsely realized by
a point in the original space. 
Moreover, a distance formula recovers the original metric, up to
quasiisometry, from the sum of sufficiently large coordinate
distances.
A \emph{hierarchically hyperbolic group} (HHG) is a group with an
HHS structure on its Cayley graph that
has additional equivariance properties.
HHSs were introduced to provide a common framework
for studying the geometry of mapping class groups and cocompactly cubulated
groups \cite{BehHagSis17}, thus providing a bridge to transfer results
from one setting to the other.

A \emph{relative HHS/HHG} structure is a relaxation of the 
requirements that allows some minimal domains to be designated
``peripheral''.
The coordinate space of a peripheral domain is not required to be hyperbolic.

Not all Coxeter groups admit an HHG structure
\cite{PetSpr23,AzuHag}.
The attempted analogy with mapping class groups breaks down when we
encounter irreducible higher rank affine parabolic subgroups, which
are not hyperbolic, but also do not decompose further by
induction on parabolic subgroups.
We build a
relative HHG structure by stopping when we reach such parabolic subgroups:
\begin{theorem}\label{thm:Cox_rHHG} 
  Let $(W,S)$ be a finite rank Coxeter system. Then $W$ is
  hierarchically hyperbolic relative to the irreducible higher rank
  affine parabolic subgroups of $(W,S)$.
\end{theorem}
The relative HHG structure we define in \fullref{sec:Cox_HHS}
is naturally related to the
Coxeter structure: the domains are certain of the parabolic subgroups,
nesting of domains is subgroup containment, orthogonality of domains
is commutation of subgroups, and coordinate projections are induced by
gate projections in $\Davis$.
This parallels the mapping class group structure, in which domains are
subsurfaces, nesting is containment of subsurfaces, and orthogonality is
disjointness of subsurfaces. 

\medskip

The natural Coxeter/parabolic induction stops at irreducible higher
rank affine parabolics.
To proceed further we equip these
parabolic subgroups with HHS structures from their crystallographic
geometries.
There is a general scheme, see \fullref{sec:promotion}, for
promoting a relative HHS structure to an absolute HHS structure by
replacing each peripheral domain with an HHS structure on its
coordinate space. 
Making the result an HHG structure requires compatibility
with the group action.
In \fullref{sec:virtual_HHG} we consider two ways to arrange such
compatibility: pass to a finite index subgroup or restrict the types of affine parabolic subgroups that occur.

A crystallographic group acts on a Euclidean space with  a
finite-index translation subgroup and a finite linear part, its
``point group''. The desired group compatibility is possible precisely when
the parabolic subgroup is \emph{hyperoctahedral}:
there exists a choice of basis that is permuted, up to
sign, by the point group.

We show that for any odd prime $p$,  the $p$--principal congruence subgroup of the Tits representation of
$W$ is a finite index subgroup $K_{(p)}\leq W$ with the property that
the action of the $K_{(p)}$--stabilizer of every irreducible higher
rank affine parabolic subgroup factors through its translation
subgroup.
We then get an HHG structure on $K_{(p)}$ (in fact, a colorable HHG structure, see
\fullref{def:colorable}):
\begin{theorem}\label{thm:Cox_vHHG}
  Let $(W,S)$ be a finite rank Coxeter system. Then $W$ is virtually
  a colorable hierarchically hyperbolic group; that is, $W$ has a finite
  index subgroup that is a colorable HHG.
\end{theorem}
Since being an HHS is preserved by quasiisometry, it follows that $W$
and $\Davis$ have HHS structures:
\begin{corollary}\label{cor:Davis_HHS}
  The Davis complex of a finite rank Coxeter system is an HHS.
\end{corollary}

\fullref{thm:Cox_vHHG} and \fullref{cor:Davis_HHS} answer
questions from the hierarchical hyperbolicity problem list
\cite{HHSProblemList}.

\medskip

The other option is to resolve the hyperoctahedral irreducible higher rank affine
parabolic subgroups and keep the nonhyperoctahedral ones as
peripherals of a refined relative HHG structure on all of $W$: 
\begin{theorem}\label{thm:refined_rHHG}
  Let $(W,S)$ be a finite rank Coxeter system. Then $W$ is
an HHG relative to the nonhyperoctahedral
  irreducible higher rank affine parabolic subgroups of $(W,S)$.
  In particular, $W$ is an HHG if it has no
  nonhyperoctahedral
  irreducible higher rank affine parabolic subgroups.
\end{theorem}
In fact, \fullref{thm:refined_rHHG} is sharp, and we give a
complete, diagrammatic characterization of being an HHG for Coxeter
groups.
The argument hinges on a Coxeter
nonhyperoctahedral escalation lemma,
\fullref{lem:highest_nonhyperoctahedral}, and
\fullref{thm:vAbelian_coarse_median}, which characterizes the existence of
coarsely equivariant coarse medians on virtually Abelian groups:
one exists if and only if the group is cocompactly cubulable, which is
equivalent to hyperoctahedrality.
We use \emph{coarse median} in the sense of Bowditch \cite{Bow13b},
see \fullref{def:NWZ} and \fullref{def:cecm}.
\begin{theorem}\label{thm:diagram_HHG}
  Let $(W,S)$ be a finite rank Coxeter system. The following are
  equivalent:
  \begin{enumerate}
  \item $W$ admits an
    HHG structure.
  \item $W$ admits a coarsely equivariant coarse median.
    \item $(W,S)$ has no nonhyperoctahedral irreducible higher rank
      affine parabolic subgroup.
    \item The Coxeter diagram of $(W,S)$ contains
  no induced subdiagram of irreducible affine type $\widetilde{A}_2$,
  $\widetilde{A}_n$ for $n\geq 4$, $\widetilde{E}_6$, $\widetilde{E}_7$,
  $\widetilde{E}_8$, $\widetilde{F}_4$, or $\widetilde{G}_2$.
  \end{enumerate}
\end{theorem}
A Coxeter system is \emph{affine-free} if $(W,S)$ has no irreducible
higher rank affine parabolic subgroups.
Niblo and Reeves \cite{NibRee03} built a CAT(0) cube complex $\NR$
whose wall structure is $W$--equivariantly isomorphic to that of the
Davis complex of $(W,S)$.
The action $W\act \NR$ is cocompact if and only if
$(W,S)$ is affine-free \cite{CapMuh05}.
In particular, if $(W,S)$ is affine-free then one can pull back the
cubical median from $\NR$ to get a coarsely equivariant coarse median
on $W$. 
Fioravanti, Levcovitz, and Sageev, \cite[Section~1.3]{FioLevSag24},
say that when a Coxeter group $W$ is not affine-free it is not known if
$W$ admits a coarsely equivariant coarse median.
\fullref{thm:diagram_HHG} determines
precisely when such a median exists. 

Azuelos and Hagen conjecture that $W$ is an HHG if and only if it is
cocompactly cubulated \cite[Conjecture~1.5]{AzuHag}; see also
\cite{HHSProblemList}.
Every cocompactly cubulated group admits a coarsely equivariant coarse
median, so \fullref{thm:diagram_HHG} establishes one direction of their
conjecture. 
For the converse, if $W$ is an HHG then it has no parabolics of the affine
types forbidden by \fullref{thm:diagram_HHG}, but we do not know if
that is a sufficient condition for the existence of a
cocompact cubulation.
The obstruction to cocompact cubulation of Munro and
Petyt \cite{MunPet}, using richly branching flats of maximal
dimension, cannot apply to Coxeter groups: it obstructs the existence
of any coarse median, not just those that are coarsely equivariant.

\medskip

New Coxeter group results emerge from the abstract HHG framework:

\begin{theorem}\label{thm:Cox_coarse_median}
Every finitely generated Coxeter group has a
  coarse median of rank equal to its algebraic rank.
\end{theorem}

The \emph{algebraic rank} is the maximal rank of a free Abelian
subgroup.
For Coxeter groups, the algebraic rank is easily determined directly from
the Coxeter diagram \cite[Theorem~6.8.3]{Kra09}.
\fullref{thm:Cox_coarse_median} is optimal in terms of rank: since $W$ is a CAT(0) group, Abelian subgroups
are undistorted, so by \cite[Corollary~2.4]{Bow13b}, the algebraic rank
is a lower bound for the rank of any coarse median on $W$.
For HHGs, the rank of the HHS coarse median is equal to the algebraic
rank, \cite[Theorem~1.15 and Corollary~2.15]{BehHagSis21} and \cite[Theorem~7.14]{AzuHag}.
Algebraic rank is invariant upon
passage to a finite index subgroup, so \fullref{thm:Cox_coarse_median}
follows.

\fullref{thm:Cox_coarse_median} says in particular that every finitely
generated Coxeter group admits a coarse median. This is a nontrivial
result in itself: 
Haettel \cite{Hae16} shows that many thick affine buildings have no
coarse medians at all, so the existence of a coarse median does not
follow  just from
the fact that $\Davis$ is a CAT(0)  polyhedral complex.

A theorem of Petyt, \cite[Theorem~B]{Pet26}, turns the ``colorable''
conclusion of \fullref{thm:Cox_vHHG} into a cubical model:
\begin{theorem}\label{thm:Cox_quasicubical}
  Every finitely generated Coxeter group is
   quasiisometric to a finite dimensional CAT(0) cube
   complex via a quasiisometry that coarsely matches the
   cubical median with the
   coarse median on the group coming from its virtual HHG structure.
 \end{theorem}

Haglund and Wise \cite{HagWis10} proved that there is a finite index
subgroup of $W$ whose action on $\NR$ is ``special''.
Because of the lack of cocompactness in the non-affine-free case,
neither Niblo-Reeves nor Haglund-Wise provides a quasiisometry to a
CAT(0) cube complex for arbitrary finitely generated Coxeter groups.

A further consequence of the virtual HHG structure concerns actions on
injective metric spaces (also known as hyperconvex spaces).
A geodesic metric space is \emph{injective} if closed balls have the
Helly property: every family of
pairwise intersecting closed balls has nonempty intersection. 
\fullref{thm:Cox_vHHG} combined with a result of Haettel, Hoda, and
Petyt \cite[Lemma~3.10]{HaeHodPet23} gives:
\begin{theorem}\label{thm:Cox_virtually_coarsely_injective}
  Every finitely generated Coxeter group is \emph{virtually coarsely
    injective}: it has a finite index subgroup
  that acts properly and coboundedly by isometries on an injective
  metric space.
\end{theorem}
Passing to a finite index subgroup is necessary, in general, because Hoda
\cite[Theorem~7.1]{Hod23} shows that a crystallographic group is
coarsely injective if and only if it is hyperoctahedral, and 
\fullref{thm:hyperoct_classification} provides examples of
nonhyperoctahedral crystallographic Coxeter groups. 

On a finite dimensional CAT(0) cube complex, the path metric induced
by taking the $\ell^\infty$ metric on each cube gives an injective
space \cite[Lemma~7.5 and Theorem~7.8]{Bow20}, so this metric on the
cube complex of \fullref{thm:Cox_quasicubical} gives a nonequivariant injective model for a Coxeter group, whereas
\fullref{thm:Cox_virtually_coarsely_injective} gives an equivariant injective model for a finite index subgroup.

\section{Preliminaries}\label{sec:prelim}
    
\subsection{Coxeter groups}
We provide a quick recap.
For further background and justifications,
see \cite{CasChuTao}.

A \emph{Coxeter system} $(W,S)$ is a Coxeter group $W$ and a
fundamental generating set $S$ of involutions such that if $m_{st}$ is the order of
the element $st$, for $s,t\in S$, then $W$ is the group with presentation:
\[W=\langle S\mid (st)^{m_{st}} \text{ for all $s,t\in S$ with $m_{st}<\infty$}\rangle\]

A \emph{special subgroup} of $(W,S)$ is a subgroup $W_T$ generated by
a subset $T\subset S$. 
A \emph{parabolic subgroup} of $(W,S)$ is a conjugate of a special
subgroup.
If $T\subset S$ then $(W_T,T)$ is a Coxeter system,
\cite[Theorem~4.1.6(i)]{davisbook}.
By \cite[Lemma~5.3.6]{davisbook}, the intersection of two parabolic
subgroups is a parabolic subgroup.
Since $S$ is finite, it follows that given a set of elements $A\subset W$, there is a
 well-defined smallest parabolic subgroup containing $A$, called its
 \emph{parabolic closure}, and denoted $\pc(A)$.
 Moreover, from the proof of  \cite[Lemma~5.3.6]{davisbook} and the
 subsequent remark, it follows that for every parabolic
subgroup $P$ of $(W,S)$, the intersection $P\cap W_T$ is a parabolic
subgroup of both $(W,S)$ and $(W_T,T)$.
In particular, the parabolic subgroups of $(W_T,T)$ are precisely the
parabolic subgroups of $(W,S)$ that are contained in $W_T$.

The Coxeter presentation can also be encoded by a \emph{Coxeter
  diagram}, which is a graph with vertex set $S$, and an edge labelled
$m_{st}$ between vertices $s,t\in S$ when $m_{st}>2$.
If the Coxeter diagram is connected then the Coxeter system is called
\emph{irreducible}. The \emph{canonical product decomposition} of
$(W,S)$ is $W=W_{S_1}\times \cdots\times W_{S_n}$, where the $S_i$ are the vertex
sets of the connected components of the Coxeter diagram. 

For $T\subset S$, let $T^\perp:=\{s\in S\setminus T\mid \forall t\in T,\,st=ts\}$.

Coxeter groups are called \emph{spherical} or \emph{affine} if they act
as geometric reflection groups on a sphere or Euclidean space, respectively.
A Coxeter group is spherical if and only if it is finite.
Since parabolic subgroups are themselves Coxeter groups, this terminology is
also used to describe parabolic subgroups. 
By extension, a subset $T\subset S$ is called spherical, affine,
irreducible, etc,
according to the corresponding property for $W_T$.

The \emph{rank} of a Coxeter system is $|S|$. An irreducible \emph{higher rank}
  affine parabolic subgroup is one of rank at least 3, which is
equivalent to the corresponding Euclidean space having dimension at
least 2.

A Coxeter system is \emph{wide} if its canonical product decomposition
into irreducible factors either contains at least two nonspherical
factors or contains an irreducible higher rank affine factor. 
A nonspherical irreducible parabolic subgroup is \emph{pre-wide} if
its normalizer is wide.
\begin{lemma}[{\cite[Lemma~5.12]{CasChuTao}}]\label{lem:pre-wide}
  If $P$ is a nonspherical irreducible parabolic subgroup that is
  contained in a wide parabolic subgroup, then $P$ is pre-wide.
\end{lemma}

A finite rank Coxeter system $(W,S)$ admits a geometric action on a
piecewise Euclidean, CAT(0) polyhedral cell complex $\Davis$, its
\emph{Davis complex}. By fixing a base vertex, the 1--skeleton
$\Davis^{(1)}$ is identified with the Cayley graph of $W$ with respect
to $S$. The path metric on $\Davis^{(1)}$ is denoted $\dcomb$. 
For $T\subset S$, the subcomplex $w\Davis_T$ of $\Davis$ spanned by
the coset $wW_T$ is called a \emph{parabolic subcomplex}.
Its stabilizer $\stab(w\Davis_T)=wW_Tw^{-1} $ is a parabolic subgroup.
When $P$ is a parabolic subgroup we write $\Upsilon\in\fibers(P)$ to mean that $\Upsilon$ is a
parabolic subcomplex with $\stab(\Upsilon)=P$, and we say $\Upsilon$ is a \emph{$P$--fiber}.
When $\Xi\subset\Upsilon$ are parabolic subcomplexes,
$\stab(\Xi)\leq\stab(\Upsilon)$.
Conversely if $P\subset Q$ are parabolic subgroups, then for every
$\Xi\in\fibers(Q)$ there exists $\Upsilon\in\fibers(P)$ with $\Upsilon\subset\Xi$.

If $\Upsilon$ is a parabolic subcomplex of $\Davis$ then
$\Upsilon^{(1)}$ is a convex subgraph of $\Davis^{(1)}$ and closest
point projection $\gate_\Upsilon\from\Davis^{(0)}\to\Upsilon^{(0)}$ is a \emph{gate map}, which means that for all
$x\in\Davis^{(0)}$ and $y\in \Upsilon^{(0)}$ we have:
\[\dcomb(x,y)=\dcomb(x,\gate_\Upsilon(x))+\dcomb(\gate_\Upsilon(x),y)\]
It follows easily that $\gate_\Upsilon$ is $1$--Lipschitz on
$\Davis^{(0)}$ and restricts to the identity on $\Upsilon^{(0)}$.
For $X\subset \Davis$, the notation $\gate_\Upsilon(X)$ means
$\bigcup_{x\in X\cap\Davis^{(0)}}\gate_\Upsilon(x)$.
We recall some further properties of gate maps to parabolic subcomplexes:
\begin{lemma}[{\cite[Proclamations~2.42--2.47]{CasChuTao}}]\label{lem:properties_of_gate_maps}
  For all $w\in W$, $x\in \Davis^{(0)}$, and parabolic subcomplexes
  $\Upsilon$ and $\Xi$:
  \begin{enumerate}
     \item $\Upsilon\subset\Xi\implies\gate_\Upsilon=\gate_\Upsilon\circ\gate_\Xi$ \label{item:gate_nested}
    \item A wall separates $x$ from $\gate_\Upsilon(x)$ if and only if
      it separates $x$ from $\Upsilon$. \label{item:wall_gate}
     \item $\gate_{w\Upsilon}(x)=w\gate_\Upsilon(w^{-1}x)$ \label{item:gate_action} 
  \item   If $T\subset S$ and $z\in W_{T^\perp}$ then
    $\gate_{w\Davis_T}(wzw^{-1}x)=\gate_{w\Davis_T}(x)$.\label{item:gate_perp}
  \item $\gate_\Upsilon(\Xi)$ spans a parabolic subcomplex of
    $\Upsilon$. \label{item:gate_residue}
    \item
      $\stab(\gate_\Xi(\Upsilon))=\stab(\gate_\Upsilon(\Xi))=\stab(\Upsilon)\cap\stab(\Xi)$\label{item:gate_bridge}
  \end{enumerate}
\end{lemma}

A conjugate of a generator is called a \emph{reflection}.
It acts on $\Davis$ as a reflection, a CAT(0)--isometric involution fixing a \emph{wall} pointwise and
exchanging the two complementary components of the wall, which are
called \emph{open halfspaces}.
Distinct walls of $\Davis$ are either disjoint or intersect
transversely, and their corresponding reflections generate an infinite
or finite dihedral group, respectively. 
For $x,y\in\Davis^{(0)}$, let $\walls(x,y)$ be the set of walls
separating $x$ from $y$; that is, the walls such that $x$ is in one
open halfspace and $y$ is in the other.
The combinatorial distance $\dcomb(x,y)$ from $x$ to $y$ in
$\Davis^{(1)}$ is exactly $|\walls(x,y)|$, and for every combinatorial
geodesic $\gamma$ from $x$ to $y$, the edges of $\gamma$ are in
bijection with walls of $\walls(x,y)$ by matching edge $e$ with the
unique wall $M$ separating the vertices of $e$. We say $M$ \emph{cuts}
$e$. Similarly, we say $M$ \emph{cuts} a parabolic subcomplex
$\Upsilon$ if it cuts some edge of $\Upsilon$.
A wall $M$ cuts $\Upsilon$ if and only if the reflection in $M$ belongs to the
stabilizer of $\Upsilon$.

A \emph{chain} is an ordered collection of distinct walls
$\mathcal{K}=\{K_i\}_{i\in I}$ for $I\subset \mathbb{Z}$ such that
there is a choice of open halfspaces $K_i^{-}$ and $K_i^+$ of wall
$K_i$ such that $i<j$ in $I$ implies $K_i^{-}\subset K_j^{-}$.
Its length is the number of walls it contains. 
All the chains in this paper have finite length.
In \cite{CasChuTao} we defined a \emph{Caprace constant}, $\caprace$,
depending on the Coxeter system.
The \emph{interior} $\interior{\mathcal{K}}$ of a chain is the
subchain obtained from $\mathcal{K}$ by discarding the
$(\caprace-1)$--outermost walls on each end, so if
$\mathcal{K}=\{K_i\}_{1\leq i\leq n}$, then
  $\interior{\mathcal{K}}=\{K_i\}_{\caprace\leq i\leq n+1-\caprace}$. 
A \emph{grid} is a pair $(\mathcal{V},\mathcal{H})$ such that
$\mathcal{V}$ and $\mathcal{H}$ are chains containing at least two
walls each, such that every wall of $\mathcal{V}$ is transverse to
every wall of $\mathcal{H}$.
 When $\mathcal{M}$ is a set of walls, let $\pc(\mathcal{M})$ be the
 parabolic closure of the corresponding set of reflections.
 When $\mathcal{V}$ is a chain of length at least 2,
 $\pc(\mathcal{V})$ is a nonspherical irreducible parabolic subgroup. 

 A \emph{reflection subgroup} is a
subgroup generated by reflections. The following is well-known:
\begin{lemma}\label{lem:reflection_splitting}
 A reflection subgroup of a direct product of special subgroups splits
 as a direct product of reflection subgroups of the factors.
 Consequently, if $P$ is an irreducible parabolic subgroup and
 $Q=Q_1\times Q_2$ is parabolic subgroup that splits as a product of
 parabolic subgroups, then:
 \[P\leq Q\iff (P\leq Q_1)\lor(P\leq Q_2)\]
\end{lemma}

\subsection{Essential parabolic subgroups}
In \cite{CasChuTao} we used known results about parabolic centralizers
and normalizers to say that a nonspherical irreducible
parabolic subgroup $P$ has an ``orthogonal complement'' $P^\perp$ that
is parabolic, equal to the centralizer $\centralizer_W(P)$ of $P$, and
such that $P\times P^\perp$ is parabolic and equal to the normalizer
$\normalizer_W(P)$ of $P$.
Then we used these facts to understand the structure of parabolic
subcomplexes stabilized by $P$.
Now, following Caprace and Marquis \cite{CapMar13} and Caprace,
Marquis, and Reid \cite{CapMarRei22}, we extend these results to a
wider class of parabolic subgroups.

\begin{definition}\label{def:essential}
  A parabolic subgroup $P$ of $(W,S)$ has a canonical
decomposition as a product of irreducible factors,
$P=Q_1\times\cdots\times Q_m\times F_1\times\cdots \times F_n$, where the
$Q_i$ are nonspherical and the $F_j$ are spherical.
Call $P^{ess}:=Q_1\times\cdots \times Q_m$ the \emph{essential part} of $P$.
Call $P$ \emph{essential} if $P=P^{ess}$.
\end{definition}
\begin{lemma}\label{lem:ess_monotone}
  If $P\leq Q$ are parabolic subgroups of $(W,S)$ then $P^{ess}\leq Q^{ess}$.
\end{lemma}
\begin{proof}
 Apply \fullref{lem:reflection_splitting} to the canonical product
 decompositions of $P$ and $Q$.
\end{proof}

\begin{lemma}[{\cite[Lemma~2.1]{CapMar13}, expanded proof as \cite[Lemma~3.3]{CapMarRei22}}]\label{lem:essential_normalizer}
  If $P$ is an essential parabolic subgroup of $(W,S)$ then there exist
  $w\in W$ and a unique $T\subset S$ such that $P=wW_Tw^{-1}$, 
  $\centralizer_W(P)=wW_{T^\perp}w^{-1}$, and
  $\normalizer_W(P)=wW_{T\sqcup T^\perp}w^{-1}$, so
  $P\times\centralizer_W(P)=\normalizer_W(P)$ and both
  $\centralizer_W(P)$ and $\normalizer_W(P)$ are parabolic subgroups.
\end{lemma}

\begin{definition}
  If $P$ is an essential parabolic subgroup of $(W,S)$, define $P^\perp:=\centralizer_W(P)$.
\end{definition}
This extends the $P^\perp$ notation for $P$ nonspherical and
irreducible from \cite{CasChuTao}.
In addition to brevity, the notation emphasizes that
$P^\perp$ is parabolic, whereas centralizers of arbitrary
parabolic subgroups need not be.

We now state two results generalizing results of \cite{CasChuTao} from
nonspherical irreducible parabolic subgroups to essential parabolic
subgroups. \fullref{lem:essential_normalizer} replaces the conclusions
that were previously obtained from the irreducibility hypothesis, and
the rest of the proofs go through unchanged.
\begin{lemma}[{cf \cite[Lemma~2.48]{CasChuTao}}]\label{lem:fiber_transitions}
  Let $P$ be an essential parabolic subgroup of $(W,S)$. Suppose
  $\Upsilon_0=w_0\Davis_{T_0}$ and $\Upsilon_1=w_1\Davis_{T_1}$ are
  elements of $\fibers(P)$.
  Then $T:=T_0=T_1$ and the unique shortest element $u$
  of $W_Tw_0^{-1}w_1W_T$ belongs to $W_{T^\perp}$ and satisfies:
  \begin{itemize}
  \item
    $\gate_{\Upsilon_1}(\Upsilon_0)=\Upsilon_0^{(0)}u=\Upsilon_1^{(0)}$
    \item Left-multiplication by $w_0uw_0^{-1}\in P^\perp$ and
      right-multiplication by $u\in W_T^\perp$ induce the same
      isomorphism $\Upsilon_0\to\Upsilon_1$.
  \end{itemize}
\end{lemma}

\begin{lemma}[{cf \cite[Corollary~2.17, Lemma~5.3]{CasChuTao}}]\label{lem:essential_normalizer_fibered}
  Let $P$ be an essential parabolic subgroup of $(W,S)$.
  \begin{itemize}
  \item  There is a
  unique parabolic subcomplex $\Upsilon_{N(P)}$ with stabilizer equal
  to $\normalizer_W(P)$.
  \item $P^\perp$ acts freely and transitively by left multiplication
    on $\fibers(P)$, and the action induces isomorphisms between the
    corresponding subcomplexes.
    \item Distinct elements of $\fibers(P)$ are disjoint, and
      $\Upsilon_{N(P)}^{(0)}=\coprod_{\Xi\in\fibers(P)}\Xi^{(0)}$.
      \item For every $\Xi\in\fibers(P)$, the composition of the
        inclusion $\Xi^{(0)}\into \Upsilon^{(0)}_{N(P)}$ with the
        quotient $\Upsilon^{(0)}_{N(P)}\onto
        P^\perp\backslash\Upsilon^{(0)}_{N(P)}$ is a bijection.
  \end{itemize}  
\end{lemma}

\subsection{Hyperoctahedrality}
The setup here follows Azuelos and Hagen \cite[Section~2.13]{AzuHag}.
A \emph{crystallographic group} is a group that acts
geometrically and faithfully on a Euclidean space $\mathbb{E}^n$.
A crystallographic group $G$ has a translation subgroup isomorphic to
$\mathbb{Z}^n$, and $G$ modulo its translation group is the finite
orthogonal ``point group''.
More generally, for $n\geq 1$, let $G$ be any group with a finite index $\mathbb{Z}^n$ subgroup.
By passing to a normal core, $G$ has a finite index normal subgroup
$A\cong \mathbb{Z}^n$.
Set $V=A\otimes_{\mathbb Z}\mathbb R$ and let $\rho\from G\to\GL(V)$
be the representation of $G$ induced by conjugation on $A$.
Its image $F=\rho(G)$ is finite, since $A\leq\ker\rho$ has finite
index in $G$.
We again call $F$ the ``point group''; in the general case it
may be a proper quotient of $G/A$.
When $G$ is crystallographic with translation subgroup $T$, the
representation $\rho$ induces the classical real holonomy representation,
which gives an isomorphism $G/T\cong F$.
\begin{definition}[{\cite[Definition~2.53]{AzuHag}, cf \cite[Definition~2.2]{Hag14crystal}}]\label{def:hyperoctahedral}
  $G$ is \emph{hyperoctahedral} if $V$ has a basis
$\{v_1,\ldots,v_n\}$ such that $F$ permutes
$\{\pm v_1,\ldots,\pm v_n\}$.
Equivalently, $F$ is conjugate in $\GL_n(\mathbb R)$ into the
``hyperoctahedral group'' of signed
permutation matrices $\orthog_n(\mathbb{Z})\cong
(\mathbb{Z}/2\mathbb{Z})^n\rtimes \sym_n$.
\end{definition}
A different choice of $A$ in $G$ yields a conjugate representation, so
$\rho$ is well-defined up to conjugacy, and hyperoctahedrality of $G$ is independent of the choice of $A$.

The following uses only classical material, but for
completeness of \fullref{thm:diagram_HHG} we work out the details.
\begin{theorem}[Hyperoctahedrality classification of irreducible
  affine Coxeter groups]\label{thm:hyperoct_classification}
  The irreducible affine Coxeter groups of types $\widetilde{A}_1$,
  $\widetilde{A}_3$, $\widetilde{B}_n$ for $n\geq2$, $\widetilde{C}_n$
  for $n\geq 2$, and
  $\widetilde{D}_n$ for $n\geq 4$ are hyperoctahedral.
  The irreducible affine Coxeter groups of types $\widetilde{A}_2$,
  $\widetilde{A}_n$ for $n\geq 4$, $\widetilde{E}_6$, $\widetilde{E}_7$,
  $\widetilde{E}_8$, $\widetilde{F}_4$, and $\widetilde{G}_2$ are
  nonhyperoctahedral.  
\end{theorem}
\begin{proof}
  The irreducible affine Coxeter groups were classified by Coxeter,
  see \cite[Figure~2 in Section~2.5]{Hum90}.
Humphreys \cite[Chapter~4]{Hum90} constructs a faithful geometric
action of the affine Coxeter group of type $\widetilde{X}_n$ on
$\mathbb{R}^n$ with point group the finite Weyl group $W(X_n)$ with its standard
reflection representation.
As a Coxeter group, $W(X_n)$ is the
irreducible spherical Coxeter group of type $X_n$.

Humphreys \cite[Section~2.10]{Hum90} works out concrete descriptions
of crystallographic root systems. 
The standard reflection representations for types $B_n$, $C_n$, and
$D_n$ permute the standard basis of the corresponding $\mathbb{R}^n$,
up to sign.
Thus, $\widetilde{B}_n$, $\widetilde{C}_n$, and $\widetilde{D}_n$ are hyperoctahedral.

For $\widetilde{A}_n$ the group $W(A_n)\cong \sym_{n+1}$ acts on
$\mathbb{R}^{n+1}$ by permuting its standard basis, but
the reflection representation acts on the $n$--dimensional subspace
$P$ of $\mathbb{R}^{n+1}$ consisting of vectors whose coordinates sum
to 0.
When $n=1$ the subspace $P$ is spanned by $p_1=e_1-e_2$ and the permutation
action on $\{e_1,e_2\}$ exchanges $p_1$ and $-p_1$, so $\widetilde{A}_1$ is hyperoctahedral.
For $n=3$, $P$ has an orthogonal basis $p_1=(1,1,-1,-1)$,
$p_2=(1,-1,1,-1)$, and $p_3=(1,-1,-1,1)$, and the set $\{\pm p_1,\pm p_2,\pm
p_3\}$ accounts for all of the $\binom{4}{2}=6$ vectors with two
coordinates equal to 1 and two equal to $-1$, so $W(A_3)$ permutes
them. Thus, $\widetilde{A}_3$ is hyperoctahedral. 

 $W(A_2)\cong \sym_3$ contains an element of order 3, so it
 cannot inject into the group $\orthog_2(\mathbb{Z})$ of order $8$.

 For $n\geq 4$, suppose $\phi\from W(A_n)\into
\orthog_n(\mathbb{Z})$.
Let $K$ be the kernel of the composition of $\phi$ with the quotient
by the normal $(\mathbb{Z}/2\mathbb{Z})^n $ factor of $\orthog_n(\mathbb{Z})$.
Then $K$ injects via $\phi$ into $(\mathbb{Z}/2\mathbb{Z})^n $, so it is a
normal 2--subgroup of $\sym_{n+1}$.
Since $n+1\geq 5$, the only normal subgroups of
$\sym_{n+1}$ are 1, $\alt_{n+1}$, and $\sym_{n+1}$.
The latter two are not 2--groups, so $K$ is trivial. This is
impossible, since $|W(A_n)|=|\sym_{n+1}|>|\sym_n|$.

For the exceptional types the point group is
bigger than the corresponding hyperoctahedral group, see \cite[Table~2 in
Section~2.11]{Hum90}.
\end{proof}

\begin{example}[{\cite{PetSpr23} and
  \cite[Figure~1]{AzuHag}}]\label{ex:hyperoctahedral_times_Z}
  The group $\widetilde{A}_2\times\mathbb{Z}$ is hyperoctahedral: its
  point group is $W(A_2)\cong\sym_3$ acting by the direct sum of its
  standard reflection representation and a trivial 1--dimensional
  representation.
  As described in the proof of \fullref{thm:hyperoct_classification},
  the reflection representation is the coordinate permutation action
  on the plane $P=\{(x_1,x_2,x_3)\in\mathbb{R}^3\mid
  x_1+x_2+x_3=0\}$.
  The linear isomorphism $P\oplus\mathbb{R}\cong
  \mathbb{R}^3:(p,t)\mapsto p+t(1,1,1)$ conjugates the point group
  action to the coordinate permutation action of $\sym_3$ on $\mathbb{R}^3$.
\end{example}

\begin{example}\label{ex:hyperoctahedral_times_Dinfty} 
  The group $\widetilde A_2\times\widetilde
A_1$ is not hyperoctahedral.
It is virtually $\mathbb{Z}^3$ with point group
$F\cong \sym_3\times \mathbb{Z}/2\mathbb{Z}$.
The conjugation of \fullref{ex:hyperoctahedral_times_Z} takes the reflection of $F$ coming
from the $\widetilde
A_1$ factor to the reflection through $P$, which is not a signed permutation
matrix. This shows why taking an $\widetilde
A_1$ factor instead of a $\mathbb{Z}$ factor breaks that specific
conjugation to signed permutation matrices. 
To see that no other such conjugation is possible, consider that 
the 3-dimensional hyperoctahedral
group is $\orthog_3(\mathbb{Z})\cong (\mathbb{Z}/2\mathbb{Z})^3\rtimes\sym_3\cong \sym_4\times
\mathbb{Z}/2\mathbb{Z}$.
The group $F$ is not conjugate into
$\orthog_3(\mathbb{Z})$ because in $F$ the order 3 rotation
coming from the $\widetilde A_2$ factor and the reflection coming from the
$\widetilde A_1$ factor commute.
Linear conjugation preserves being a reflection, but there is no pair
of commuting order 3 element and reflection in
$\orthog_3(\mathbb{Z})$---one can check every element $g$ of order 3
in $\orthog_3(\mathbb{Z})$ has centralizer $\langle
g\rangle\times \{\pm\mathrm{I}_3\}$, which contains no reflections.
\end{example}

\begin{remark}\label{rem:conjugation}
\fullref{ex:hyperoctahedral_times_Dinfty} illustrates the difference between $F$ being abstractly isomorphic to a subgroup of
$\orthog_n(\mathbb{Z})$ and being conjugate into
$\orthog_n(\mathbb{Z})$; 
$\orthog_3(\mathbb{Z})$ does contain a subgroup isomorphic
to $\sym_3\times \mathbb{Z}/2\mathbb{Z}$. 
\end{remark}

  \subsection{Hierarchical Hyperbolicity}\label{sec:HHS}
  Behrstock, Hagen, and Sisto defined hierarchically hyperbolic spaces
in \cite{BehHagSis17}.
In a sequel \cite{BehHagSis19} the defining axioms were refined.
  Even for the refined axioms there are two versions: 
\cite[Definition~1.1]{BehHagSis19} and a streamlined version
\cite[Section~1C]{BehHagSis19}, which are equivalent when
the Projection Axiom maps are coarsely surjective
\cite[Proposition~1.11]{BehHagSis19}.
The (relative) HHS structure that we will define on a Coxeter group has coarsely
surjective projections, so we state the streamlined definition, with
the coarse surjectivity hypothesis built in to the Projection Axiom.
\begin{definition}[{\cite[Section~1C variant]{BehHagSis19}, see also \cite{BerRus22,Per26}}]\label{def:HHS}
  A \emph{hierarchically hyperbolic space} (HHS) consists, for some
  parameter $E\geq 1$, of an  $(E,E)$--quasigeodesic
  metric space $X$, an \emph{index set} $\HHindex$ of \emph{domains}, and
  a collection of geodesic, $E$--hyperbolic \emph{coordinate
    spaces} $\{(\coordinate{P},\delta_P)\mid P\in \HHindex\}$, satisfying the
  following axioms:
  \begin{hhsaxioms}
  \axiom{Projection}\label{ax:projections}
    For every domain $P\in\HHindex$ there is a projection
    $\proj{P}\from X\to 2^{\coordinate{P}}\setminus\{\emptyset\}$ such
    that every point of $\coordinate{P}$ is within distance $E$ of
    some point in $\proj{P}(X)$ and such that
         for all $x,y\in X$, $\diam(\proj{P}(x)\cup\proj{P}(y))\leq Ed_X(x,y)+E$.

    \axiom{Nesting}\label{ax:nesting} If $\HHindex\neq\emptyset$ then it is partially ordered by the \emph{nesting}
   relation $\nested$, and has a unique maximal element.
   For $P\in\HHindex$, define $\nestedin{P}:=\{Q\in\HHindex\mid
   Q\nested P\}$.

   \axiom{Finite complexity}\label{ax:complexity}
   Chains in the partially ordered
   set $(\HHindex,\nested)$ have length at most $E$.
 \axiom{Orthogonality}\label{ax:orthogonality} $\HHindex$ has a symmetric, antireflexive
   \emph{orthogonality} relation $P\orthogonal Q$.
   Orthogonality is downward closed with respect to nesting and
   implies incomparability in the nesting ordering:
   \begin{enumerate}[label=(\alph*)]
   \item $(P\nested Q)\land(Q\orthogonal R)\implies P\orthogonal R$
     \item $P\orthogonal Q\implies (P\not\nested Q)\land(Q\not\nested P)$
     \end{enumerate}
     For $P\in\HHindex$, define
     $\orthogonalto{P}:=\{Q\in\HHindex\mid P\orthogonal Q\}$.

     \axiom{Containers}\label{ax:containers} If $P\nestedneq Q$ and $\nestedin{Q}\cap\orthogonalto{P}
  \neq\emptyset$ then there exists $R\nestedneq Q$, a \emph{container}
  of $P$ in $Q$, such that
  $\nestedin{Q}\cap\orthogonalto{P}\subset\nestedin{R}$.
  
  \axiom{Consistency}\label{ax:consistency} If $P$ and $Q$ are not orthogonal and not
   nested then they are \emph{transverse}, $P\transverse Q$.
  For all $x\in X$ and 
   $P,Q,R\in\HHindex$, if $P\nestedneq Q$ or $P\transverse Q$ there is
   a specified subset $\emptyset\neq\relprojection{P}{Q}\subset \coordinate{Q}$ of diameter at most
   $E$ and these subsets satisfy the following consistency conditions:
   \begin{enumerate}[label=(\alph*)]
   \item $\left(P\nested Q\nestedneq R\right)\lor\left(R\not\orthogonal
       P\nested Q \transverse R\right) \implies
     \domaindist{R}{\relprojection{P}{R}}{\relprojection{Q}{R}}\leq E$\label{ax:consistency:small} 
     \item $P\transverse Q\implies \min\{\domaindist{P}{x}{\relprojection{Q}{P}},\,\domaindist{Q}{x}{\relprojection{P}{Q}}\}\leq E$\label{ax:consistency:Behrstock}
   \end{enumerate}
   
\axiom{Bounded geodesic image}\label{ax:BGI}
If $P\nestedneq Q$ and there are points $x,y\in X$ with
$\domaindist{P}{x}{y}>E$ then every geodesic from $\proj{Q}(x)$ to
$\proj{Q}(y)$ in $\coordinate{Q}$ passes within distance $E$ of $\relprojection{P}{Q}$.

\axiom{Large links}\label{ax:largelink}
   For all $x,y\in X$ and
   $R\in\HHindex$, there exist
     domains $P_1,\dots,P_n$, properly nested in $R$, with $n\leq N:=E
     \domaindist{R}{x}{y}+E$, 
and such that:
   \begin{enumerate}[label=(\alph*)]
   \item $\forall
     i,\,\domaindist{R}{x}{\relprojection{P_i}{R}}\leq N$
   \item $\forall
     Q\in \HHindex,\, (Q\nestedneq R)\land(\domaindist{Q}{x}{y}\geq
     E)\implies \exists i,\,Q\nested P_i$
   \end{enumerate}
   
\axiom{Partial realization}\label{ax:realization}
For every collection $\{P_i\}$ of pairwise orthogonal domains and
points $p_i\in\proj{P_i}(X)\subset\coordinate{P_i}$ there exists $x\in
X$ such that for all $i$:
   \begin{enumerate}[label=(\alph*)]
\item $\domaindist{P_i}{p_i}{x}\leq E$\label{ax:realization:location}
  \item $\forall Q\in\HHindex,\,\left(P_i\nestedneq Q\right)\lor\left( P_i\transverse Q\right)\implies \domaindist{Q}{x}{\relprojection{P_i}{Q}}\leq E$\label{ax:realization:consistency}
\end{enumerate}

\axiom{Uniqueness}\label{ax:uniqueness}
There is a \emph{uniqueness function} $\theta\from
[0,\infty)\to[0,\infty)$ such that for all $x,y\in X$ and $r\geq 0$, if $d_X(x,y)\geq
\theta(r)$ then there exists $P\in\HHindex$ such that
$\domaindist{P}{x}{y}\geq r$.
\end{hhsaxioms}
\end{definition}

\begin{definition}[{cf \cite[Definition~7.1]{AbbBehDur21}, \cite[Definition~3.2]{Per26}}]
 An HHS has \emph{clean containers} if in the Container Axiom it is
 always possible to choose a
 container of $P$ in $Q$ that is orthogonal to $P$.
\end{definition}

\begin{definition}[{\cite[Definition~6.8]{BehHagSis19},\cite[Definition~2.17]{Per26}}]
  A \emph{relatively hierarchically hyperbolic space} is a space $X$
  with domains and coordinate spaces satisfying all of the
  requirements of a hierarchically hyperbolic space of \cite[Definition~1.1]{BehHagSis19}, except that
some of the $\nested$--minimal domains are designated
\emph{peripheral}, and the coordinate spaces for peripheral domains are not required to be hyperbolic.
\end{definition}

\begin{lemma}\label{relHHS_is_relHHS}
  If a space $X$ has a structure $\HHindex$ satisfying all of
  the conditions of \fullref{def:HHS}, except that some of the
  $\nested$--minimal domains are designated peripheral, and the
  coordinate spaces of peripheral domains are not required to be
  hyperbolic, then $\HHindex$ defines a relative hierarchically
  hyperbolic structure on $X$.
\end{lemma}
\begin{proof}
  \cite[Proposition~1.11]{BehHagSis19} promotes the streamlined
  HHS axioms (with uniformly coarsely surjective projections), as
  enumerated in \fullref{def:HHS}, to the full HHS
  axioms of \cite[Definition~1.1]{BehHagSis19}.
  The step that requires work is to define downward projections
$\relprojection{Q}{P}\from \coordinate{Q}\to 2^{\coordinate{P}}$ for
$P\nestedneq Q\in\HHindex$.
Their construction uses that $\coordinate{Q}$ is hyperbolic.
That is true here, since $P\nestedneq Q$ implies $Q$ is not
peripheral. 
\end{proof}

\begin{definition}[{\cite[Definition~2.27]{BerRus22}}]\label{def:relHHG}
  A finitely generated group $G$ is a \emph{relatively hierarchically
  hyperbolic group} if it has a Cayley graph $X$, with respect to some
  finite generating set, that admits a relatively hierarchically
  hyperbolic structure $(X,\HHindex)$ such that there is a cofinite action of $G$ on $\HHindex$ that preserves
  nesting and orthogonality, preserves the set of peripheral domains, and for which there exists $E\geq 0$  such that for all $g,h\in G$, 
  $P,Q\in\HHindex$, and $x\in X$, the following conditions are satisfied:
  \begin{itemize}
    \item There is an isometry $g_P\from
  \coordinate{P}\to\coordinate{g.P}$.
  \item $(gh)_P=g_{h.P}\circ h_P$
    \item $\domaindist{g.P}{g_P(\proj{P}(x))}{\proj{g.P}(gx)}\leq E$
      \item $\left(P\nestedneq
          Q\right)\lor\left(P\transverse Q\right)\implies
        \domaindist{g.Q}{g_Q(\relprojection{P}{Q})}{\relprojection{g.P}{g.Q}}\leq E$ 
  \end{itemize}
Saying $G$ is \emph{hierarchically hyperbolic relative to} $\mathfrak{P}$ identifies $\mathfrak{P}$
as the set of peripheral domains.
\end{definition}

\section{Coned off parabolic spaces}\label{sec:coneoff}
    The main character of \cite{CasChuTao} is the space
$\coneoff$ obtained from the 1--skeleton of $\Davis$ by coning off the
vertex set of every wide parabolic subcomplex. It is a hyperbolic,
Morse-recognizing space for $\Davis$, by \cite[Corollary~4.20]{CasChuTao}.

\begin{lemma}\label{lem:essential_with_wide_normalizer_equals_unbounded_in_coneoff}
  Let $(W,S)$ be a nonspherical, finite rank Coxeter system.
  Suppose $W_T$ is a nontrivial, essential special subgroup of
  $(W,S)$.
  \begin{itemize}
    \item If $\normalizer_W(W_T)$ is wide then $\Davis_T$ has diameter 1 in
      $\coneoff$. In particular, this holds when $W$ is wide. 
      \item If $\normalizer_W(W_T)$ is not wide then $\Davis_T$ has
        infinite diameter in $\coneoff$.
        \item If $W$ is not wide then $\coneoff$ is unbounded hyperbolic.
  \end{itemize}
\end{lemma}
\begin{proof}
  Suppose $W_T$ is a nontrivial, essential special subgroup. 
\fullref{lem:essential_normalizer} says
$\normalizer_W(W_T)=W_{T\sqcup T^\perp}=W_T\times W_{T^\perp}$.
If $\normalizer_W(W_T)$ is wide then $\Davis_{T\sqcup T^\perp}$ is
coned off in $\coneoff$, so $\Davis_T\subset\Davis_{T\sqcup T^\perp}$
has diameter 1 in $\coneoff$.
If $W_T$ is reducible then it is already wide, so $\normalizer_W(W_T)$
is too.
If $W_T$ is irreducible then \fullref{lem:pre-wide} says that if $W_T$
is contained in a
wide parabolic subgroup then $\normalizer_W(W_T)$ is wide.
This is the case if $W$ itself is wide.

If $\normalizer_W(W_T)$ is not wide then the previous paragraph says
that $W_T$ is irreducible and not contained in any wide parabolic subgroup.
Take a Coxeter element $g\in W_T$, that is, $g:=t_1\cdots t_n$,
where $T=\{t_1,\dots,t_n\}$.
Then $g$ is an infinite order element of $W$ that is not contained in
any wide parabolic subgroup, since such a subgroup would also contain
$\pc(g)=W_T$, \cite[Theorem~3.1]{Par07}. 
\cite[Proposition~4.5]{CapFuj10} says that an
infinite order element is rank-one if and only if it is not contained
in a wide parabolic subgroup. The Morse-recognizing property of
$\coneoff$ then implies that a rank-one element acts loxodromically on
$\coneoff$, so $n\mapsto g^n$ is a quasigeodesic. 
Thus, $\langle g\rangle\subset\Davis_T$ has
infinite diameter in $\coneoff$.

If $W$ is not wide then this argument applies to $W_T=W^{ess}$, so $\coneoff$
is unbounded.
\end{proof}

In \cite[Section~5.1]{CasChuTao} we also defined an auxiliary coned off graph for
parabolic subgroups/complexes in which all proper subparabolics are
coned off. 
Given a parabolic subcomplex $\Upsilon$, define its projection
  pseudo-metric on $\Davis^{(0)}$ by
$d_\Upsilon(x,y):=\dcomb(\gate_\Upsilon(x),\gate_\Upsilon(y))$.
Define its \emph{coned off graph} to be the metric graph $C(\Upsilon)$
obtained from $\Upsilon^{(1)}$ by coning off the 0-skeleton of every
proper parabolic subcomplex $\Xi$ of $\Upsilon$ by adding a new vertex
$v_\Xi$ and connecting $v_\Xi$ to every vertex in $\Xi^{(0)}$ by an
edge of length $\sfrac{1}{2}$.
We get another projection pseudo-metric on $\Davis^{(0)}$ by
$\delta_\Upsilon(x,y):=d_{C(\Upsilon)}(\gate_\Upsilon(x),\gate_\Upsilon(y))$.

A single nontrivial essential parabolic subgroup $P$ can have many
parabolic subcomplexes for which it is the stabilizer.
However, \fullref{lem:fiber_transitions} and
\fullref{lem:essential_normalizer_fibered} say these $P$--fibers
combine to span the unique parabolic subcomplex stabilized by
$\normalizer_W(P)=P\times P^\perp$, and the $P^\perp$--action permutes
the $P$--fibers freely and transitively, inducing isomorphisms between
fibers. Furthermore, it follows from 
\fullref{lem:properties_of_gate_maps}\eqref{item:gate_action} and
\eqref{item:gate_perp} that the $P^\perp$--action induces
canonical isomorphisms between $C(\Upsilon_0)$ and $C(\Upsilon_1)$ for
any $\Upsilon_0,\Upsilon_1\in\fibers(P)$.

This allows us to extend the definitions $\gate_\Upsilon$,
$d_\Upsilon$, $C(\Upsilon)$, and $\delta_\Upsilon$ for
$\Upsilon\in\fibers(P)$ to $\gate_P$, $d_P$, $C(P)$, and $\delta_P$.
Formally, set $\gate_P(x)$ to be the $P^\perp$--orbit of
$\gate_\Upsilon(x)$, set $d_P(x,y):=\dcomb(\gate_P(x),\gate_P(y))$,
define $C(P)$ to be the quotient of
$\coprod_{\Upsilon'\in\fibers(P)}C(\Upsilon')$ by the
$P^\perp$--action, and let
$\delta_P(x,y):=d_{C(P)}(\gate_P(x),\gate_P(y))$. 
Then $C(P)$ is canonically isomorphic to $C(\Upsilon)$, and none of
these definitions depend on the choice of $\Upsilon\in\fibers(P)$.

\begin{lemma}\label{lem:irreducible_coordinates_hyperbolic}
  There is an $\epsilon$ depending on $(W,S)$ such that if 
  $P$ is a nonspherical irreducible parabolic subgroup of $(W,S)$ that is not
  higher rank affine, then
  $C(P)$ is $\epsilon$--hyperbolic.

  Moreover, there exists $D$ depending on $(W,S)$ such that for all
  nonspherical irreducible parabolic
  subgroups $P$ of $(W,S)$ that are not higher rank affine, all
  $\Upsilon\in\fibers(P)$, all
  $x,y\in\Upsilon^{(0)}$, all combinatorial geodesics $\gamma$ from $x$
  to $y$ in $\Upsilon^{(1)}$, and all geodesics
  $\alpha$ from $x$ to $y$  in $C(\Upsilon)$, the Hausdorff distance
  between $\alpha$ and $\gamma$ in $C(\Upsilon)$ is at most $D$. 
\end{lemma}
\begin{proof}
  Let $\Upsilon\in\fibers(P)$.
  The parabolic subgroup $P$ with its inherited Coxeter system has Davis
  complex isomorphic to $\Upsilon$, so, 
  by \cite[Corollary~4.20]{CasChuTao}, $\widehat\Upsilon$ is an $\epsilon_{\widehat \Upsilon}$--hyperbolic graph.

The graph $C(P)$ is canonically isomorphic to the
graph $C(\Upsilon)$ obtained from $\Upsilon^{(1)}$
by coning off all proper parabolic subcomplexes.
We can perform these cone-offs in two stages, first coning off
proper wide parabolic subcomplexes to get $\widehat\Upsilon$ , and
then coning off the rest to get $C(\Upsilon)$.
If $\Xi$ is a proper parabolic subcomplex of $\Upsilon$ then $\Xi$ is
convex, so every pair of vertices in $\Xi$ is connected by a
combinatorial geodesic $\gamma$ of $\Davis$ that stays inside of $\Xi$.
In $\widehat{\Upsilon}$, the image of $\gamma$ is an unparameterized
uniform quasigeodesic
\cite[Corollary~4.5 and Corollary~4.20]{CasChuTao}, which is therefore uniformly Morse, since
$\widehat{\Upsilon}$ is hyperbolic.
It follows that the image of $\Xi$ in $\widehat{\Upsilon}$ is
quasiconvex.
Thus, $C(\Upsilon)$ is obtained from the hyperbolic graph
$\widehat{\Upsilon}$ by coning off a collection of uniformly
quasiconvex subspaces. 
A result of  Kapovich and Rafi \cite[Proposition~2.6]{KapRaf14} then
implies that $C(\Upsilon)$ is $\epsilon_\Upsilon$--hyperbolic,
for some $\epsilon_\Upsilon$ depending on $\epsilon_{\widehat
  \Upsilon}$ and the quasiconvexity constant.
Set $\epsilon_P:=\epsilon_\Upsilon$, which is independent of the
choice of $\Upsilon\in\fibers(P)$.
Moreover, conjugates of $P$ have fibers that differ by isometries of
$\Davis$, so $\epsilon_P$ only depends on the conjugacy class of $P$. 
Since there are only finitely many conjugacy classes of parabolic
subgroup of $(W,S)$, we get uniform $\epsilon$ by maximizing over the
finitely many values $\epsilon_P$. 

For the Hausdorff claim, $\gamma$ becomes a uniform
unparameterized quasigeodesic in the hyperbolic space
$\widehat{\Upsilon}$, so it is uniformly Hausdorff close to any
$\widehat{\Upsilon}$--geodesic $\beta$ with the same endpoints.
The map from $\widehat{\Upsilon}$ to $C(\Upsilon)$ is $1$--Lipschitz,
so $\gamma$ and $\beta$ remain uniformly Hausdorff close in
$C(\Upsilon)$.
The second part of \cite[Proposition~2.6]{KapRaf14} says that $\beta$
is uniformly Hausdorff close to any $C(\Upsilon)$--geodesic $\alpha$
with the same endpoints. This gives an $\Upsilon$--specific Hausdorff
distance bound, but again it actually is conjugacy invariant, so we
can maximize over the finitely many conjugacy classes of parabolics to
get a uniform bound.
\end{proof}

\begin{proposition}\label{lem:coordinate_diameters}
  Let $P$ be a nontrivial essential parabolic subgroup.
  \begin{itemize}
  \item If $P$ is reducible then $C(P)$ has diameter at most 3.
    \item If $P$ is irreducible affine then $C(P)$ is quasiisometric
      to $\mathbb{Z}^{\rank(P)-1}$.
      \item If $P$ is irreducible and not higher rank affine then
        $C(P)$ is hyperbolic and unbounded.
  \end{itemize}
\end{proposition}
The second and third cases overlap for infinite dihedral $P$.
\begin{proof}
  It suffices to assume $P=W_T$ is a special subgroup.
  If it is essential and reducible then let $W_{T_0}$ be one
  nonspherical irreducible factor of $W_T$ and let $T_1:=T\setminus T_0$.
For all $x,y\in W_T$, write $x^{-1}y=ab$ for $a\in W_{T_0}$ and $b\in
W_{T_1}$. Then $x$ and $xa$ are vertices contained in the proper
parabolic subcomplex $x\Davis_{T_0}$ of $\Davis_T$, while $xa$ and
$y=xab$ are contained in the proper parabolic subcomplex
$xa\Davis_{T_1}$ of $\Davis_T$. Both of these proper parabolic
subcomplexes are coned off in $C(\Davis_T)$, so
$d_{C(\Davis_T)}(x,y)\leq 2$.
Every cone vertex of $C(\Davis_T)$ is at distance $\sfrac{1}{2}$ from some
$\Davis_T^{(0)}$ vertex, so $3$ is an upper diameter bound for all
vertices of $C(\Davis_T)$.

Suppose $P=W_T$ is irreducible and affine. 
Every proper parabolic subgroup of an irreducible affine Coxeter group
is spherical, so $C(\Davis_T)$ is constructed from the 1--skeleton of
$\Davis_T$ by coning off vertices of spherical parabolic
subcomplexes.
Spherical parabolic subcomplexes of $\Davis$, hence of $\Davis_T$,
have uniformly bounded diameter, so $C(\Davis_T)$ is quasiisometric to
$\Davis_T$, hence is unbounded, since $W_T$ is nonspherical.
Since $W_T$ is affine, it has a finite index
translation subgroup isomorphic to $\mathbb{Z}^{|T|-1}$.
Thus, $C(\Davis_T)\cong C(W_T)$, $\Davis_T$, $W_T$, and
$\mathbb{Z}^{|T|-1}$ are all quasiisometric to one another.

Now suppose $P=W_T$ is irreducible and not higher rank affine. 
Let $T=\{t_0,\dots,t_{r-1}\}$.
Take a Coxeter element $c:=t_0\cdots t_{r-1}$ of $W_T$, which has
$\pc(c)=W_T$, \cite[Theorem~3.1]{Par07}.

By \cite[Corollary~4.20]{CasChuTao}, $\widehat{\Davis}_T$ is
hyperbolic, and the proof of
\fullref{lem:essential_with_wide_normalizer_equals_unbounded_in_coneoff}
shows that $c$ acts as a rank-one element on $\Davis_T$ and acts loxodromically on $\widehat{\Davis}_T$.
The proof of \fullref{lem:irreducible_coordinates_hyperbolic} shows
that $C(\Davis_T)$ is a hyperbolic graph obtained from
$\widehat{\Davis}_T$ by coning off further proper parabolic
subcomplexes of $\Davis_T$.
It remains to show that $C(\Davis_T)$ has unbounded
diameter, which we establish by showing that $c$ remains loxodromic in $C(\Davis_T)$.

Let $\gamma$ be the biinfinite path through 1 in
  $\Davis_T$ that starting from 1 in one direction has edge labels
  $t_0$, $t_1$, \dots,$t_{r-1}$, $t_0$,\dots, and in the other direction
  has edge labels $t_{r-1}$, $t_{r-2}$, \dots, $t_0$, $t_{r-1}$,\dots, so that
  $\gamma$ is an axis for the $c$--action.
  By \cite[Theorem~1]{Spe09}, $\gamma$ is a combinatorial geodesic.

Since $c$ acts as a rank-one element on $\Davis_T$, its axis $\gamma$ is a Morse quasigeodesic in the
  CAT(0) metric, hence is a Morse geodesic in the combinatorial metric.
  Thus, there exists $R$ such
  that every combinatorial geodesic segment with endpoints on $\gamma$ is contained
  in the $R$--neighborhood of $\gamma$, in the combinatorial metric $\dcomb$.

  For $n>0$ suppose $m:=d_{C(\Davis_T)}(1,c^n)$.
  By \cite[Lemma~5.21]{CasChuTao}, there is a factorization
  $c^n=h_1\cdots h_m$ such that for each $i$ there is $J_i\subsetneq
  T$ with $h_i\in W_{J_i}$ and a geodesic $\beta$ from 1 to $c^n$ that is a
  concatenation of geodesic subsegments $\beta_i$ labelled by $h_i$.

  Let $a_i$ and $b_i$ be the end vertices of $\beta_i$. There exist
  $a_i'$ and $b_i'$ on $\gamma$ at distance at most $R$ from $a_i$ and
  $b_i$, respectively.
  We can choose $a_1'=a_1$, $a_{i+1}'=b'_i$, and $b'_m=b_m$.
  
  Suppose, for some $i$, that $\dcomb(a_i',b_i')\geq r(2R+1)$.
  Subdivide the subsegment $\gamma_i$ of $\gamma$ between $a_i'$ and $b_i'$ into
  (possibly trivial) subsegments such that each edge is dual to a wall
  that cuts the proper parabolic subcomplex $\Upsilon_i$ of type $J_i$
  containing $\beta_i$, separated by single edges dual to a wall that
  does not cut $\Upsilon_i$.
  Now, a wall $M$ that separates $a_i'$ and $b_i'$ but does not cut
  $\Upsilon_i$ does not separate $a_i$ and $b_i$, so either $M$
  separates $a_i$ and $a_i'$ or $M$ separates $b_i$ and $b_i'$.
  There are at most $R$ walls in each of those cases, so $\gamma_i$ is
  subdivided into at most $2R+1$ subsegments with all edges cut by a
  wall that cuts $\Upsilon_i$, separated by $2R$ single edges that do
  not cut $\Upsilon_i$. By the distance hypotheses, at least one of
  the former has length at least $r$, so $\gamma_i$ contains a
  subsegment of length $r$ whose edges are dual to walls
  $M_1$, \dots, $M_r$ with each $M_i$ cutting $\Upsilon_i$.
  On one hand, that means that the reflections in the walls $M_j$ are
  all in $\stab(\Upsilon_i)\lneq W_T$.
On the other hand, if $g$ is the initial vertex of the edge cut by
  $M_1$, then there is $0\leq k<r$ such that, taking all subscripts modulo
  $r$, the reflection in $M_1$ is
  $gt_kg^{-1}$, the reflection in $M_2$ is $gt_kt_{k+1}t_kg^{-1}$,\dots, the reflection in $M_r$ is $gt_k\cdots
  t_{k+r-2}t_{k+r-1}t_{k+r-2}\cdots t_kg^{-1}$.
  Thus the subgroup generated by the reflections in the walls
  $M_1$,\dots,$M_r$ is $g\langle t_0,\dots,t_{r-1}\rangle g^{-1}=W_T$. 
  This is a contradiction, so $\dcomb(a_i',b_i')< r(2R+1)$.
Since this bound holds for each $i$:
  \[nr=|c^n|_T<m\cdot r(2R+1)\implies d_{C(\Davis_T)}(1,c^n)=m>\frac{n}{2R+1}\]
  Thus, 
  $c$ acts loxodromically on $C(\Davis_T)$, 
  which implies that $C(\Davis_T)\cong C(W_T)$ is unbounded.
\end{proof}

\section{Relative HHG structure on \texorpdfstring{$(W,S)$}{(W,S)}}\label{sec:Cox_HHS}
To prove \fullref{thm:Cox_rHHG}, we define a relative HHS structure on
$\Davis^{(1)}=\operatorname{Cay}(W,S)$ in
\fullref{sec:define_HHS_structure}, verify the HHS axioms in
\fullref{sec:verify_HHS}, and then check the relative HHG equivariance conditions  in \fullref{sec:verify_equivariance}.
    \subsection{Defining the relative HHS structure}\label{sec:define_HHS_structure}
Given a finite rank Coxeter system $(W,S)$, define the index set $\HHindex$ of domains to consist of:
\begin{enumerate}
\item A single, formal, top domain $\top$.
\item Every nontrivial, essential parabolic subgroup $P$ of
  $(W,S)$ such that $P\neq W^{ess}$ and $\normalizer_W(P)$ is wide.
\end{enumerate}

Define $\top\nested\top$ and $\top\not\orthogonal \top$, and for all
$P\in\HHindex\setminus\{\top\}$, define $\top\not\nested P\nested\top$ and $P\not\orthogonal \top$.
For $P,Q\in\HHindex\setminus\{\top\}$, define $P\nested Q$ if and only if
$P\leq Q$, and define $P\orthogonal Q$ if and only if $P$ and $Q$
commute (element-wise).

By Moussong's Theorem \cite{Mou88}, $W$ is hyperbolic if and only if
it has no wide parabolic subgroups. 
In this case, $\HHindex=\{\top\}$ and $W$ is already an HHG by taking $\coordinate{\top}=\Davis^{(1)}$. 
If $W^{ess}$ is irreducible higher rank affine then it has no
nonspherical proper parabolic subgroups, so again $\HHindex=\{\top\}$. Set
$\coordinate{\top}:=\Davis^{(1)}$. This is a trivial HHG structure on
$W$ relative to $\HHindex$.
In both of these cases, \fullref{thm:Cox_rHHG} is confirmed, so 
assume we are not in either. It follows that $\HHindex\neq\{\top\}$.

If $W$ is wide, equivalently, since we have assumed that $W^{ess}$ is
not higher rank affine, if $W^{ess}$ is reducible, define
$\coordinate{\top}$ to be a point, and define $\proj{\top}$ to be the
unique map $\Davis^{(0)}\to\coordinate{\top}$.

If $W$ is not wide then define $\coordinate{\top}:=\coneoff$, and define $\proj{\top}$ to be
the canonical map identifying $\Davis^{(0)}$ with the non-cone
vertices of $\coneoff$. 

For reducible $P\in\HHindex\setminus\{\top\}$, define
$\coordinate{P}$ to be a single point, and define $\proj{P}$ to be the
unique map $\Davis^{(0)}\to\coordinate{P}$.

For irreducible $P\in\HHindex\setminus\{\top\}$, define 
$\coordinate{P}:=C(P)$ and $\proj{P}:=\gate_P$, as in \fullref{sec:coneoff}.

The formal top domain $\top$ is in many ways playing the role of $W^{ess}$, and so
we will say $\top$ is ir/reducible if $W^{ess}$ is ir/reducible. 
We formally distinguish $\top$ from $W^{ess}$ because:
\begin{enumerate}
\item We do not require that $\normalizer_W(W^{ess})=W$ is wide.
  \item In the irreducible case,  $\coordinate{\top}$ is $\coneoff$, not 
 $C(W^{ess})$.
\end{enumerate}

 We chose to make the coordinate spaces for reducible domains single
 points, rather than the bounded spaces $\coneoff$ or $C(P)$ that
 would have resulted from a uniform definition.
 This emphasizes that they are present in the hierarchy for formal
 reasons (to provide a maximal element and to provide
containers).

\fullref{thm:Cox_rHHG} says the peripherals are the
irreducible higher rank affine parabolic subgroups. By construction,
these are all domains. We now check that they are precisely the
domains with non-hyperbolic coordinate spaces, and that they are all $\nested$--minimal:
\begin{lemma}\label{lem:nonaffine_domains_have_hyperbolic_coordinates}
  A domain $P\in\HHindex$ has non-hyperbolic coordinate space
  $\coordinate{P}$ if and only if $P$ is an irreducible higher
  rank affine parabolic subgroup, and every such domain is
  $\nested$--minimal.
  The other domains have uniformly hyperbolic coordinate spaces. 
\end{lemma}
\begin{proof}
   For reducible $P\in\HHindex$ the coordinate space $\coordinate{P}$
   is a single point, hence 0--hyperbolic.

   If $\top$ is irreducible then $\coordinate{\top}=\coneoff$ is unbounded hyperbolic, by
  \fullref{lem:essential_with_wide_normalizer_equals_unbounded_in_coneoff}.
   
  If $P$ is an irreducible affine parabolic subgroup then,
  as in the proof of \fullref{lem:coordinate_diameters}, $P$ has no
  proper nonspherical parabolic subgroups and is
  quasiisometric to $\mathbb{Z}^{\rank(P)-1}$.
  The former implies $P$ is $\nested$--minimal in $\HHindex$, and the
  latter implies $\coordinate{P}$ is non-hyperbolic when $\rank(P)\geq 3$.

  For non higher rank
  affine, irreducible $P\in\HHindex\setminus\{\top\}$, \fullref{lem:coordinate_diameters} says
  $\coordinate{P}$ is unbounded hyperbolic.
\end{proof}

\subsection{Verifying the HHS axioms}\label{sec:verify_HHS}
We verify the HHS axioms one at a time.
\subsubsection{Projections}
All the projections are $1$--Lipschitz and 1--coarsely surjective.

\subsubsection{Nesting}
The parabolic subgroups of $(W,S)$ are partially ordered by
inclusion. Our $(\HHindex,\nested)$ is a partially ordered subset of
that partially ordered set, to which we added a formal maximal
element $\top$.

\subsubsection{Finite complexity}
If $P\lneq Q$ are nontrivial essential parabolic subgroups then
$2\leq\rank(P)<\rank(Q)$.
Since non-top domains are proper nontrivial essential parabolic
subgroups properly contained in $W^{ess}$, their ranks are
between $2$ and $\rank(W^{ess})-1$, so a chain of them has
length at most $\rank(W^{ess})-2$. Adding $\top$ gives that
chains in $(\HHindex,\nested)$ have length at most
$\rank(W^{ess})-1\leq |S|-1$.

\subsubsection{Orthogonality}
By definition, for all $P\in \HHindex$, we have $P\nested \top$ and
$P\not\orthogonal\top$.
For $P,Q\in\HHindex\setminus\{\top\}$, if $P\nested Q$ then $P\leq Q$.
Let $P=P_1\times\cdots \times P_n$ be the canonical decomposition of $P$ into
nonspherical irreducible factors. If $P$ and $Q$ commute then, in
particular, $P_1$ commutes with itself, but that is impossible, since
nonspherical irreducible parabolic subgroups have trivial center.
Thus, $P\nested Q$ implies $P\not\orthogonal Q$.

We have shown $P\nested Q$ implies $P\not\orthogonal Q$, so,
contrapositively, if $P\orthogonal Q$ then $P\not\nested Q$ and
$Q\not\nested P$.

Now suppose  $P,Q,R\in \HHindex$ with $P\nested Q$ and
$Q\orthogonal R$. No domain is orthogonal to $\top$, so
$P,Q,R\in\HHindex\setminus\{\top\}$ are parabolic subgroups, and we
have that $P\leq Q$ and $Q$ and $R$ commute, so $P$ and $R$ commute, so $P\orthogonal R$.

\subsubsection{Clean containers}
Consider $P,Q\in\HHindex\setminus\{\top\}$, with $P\nestedneq Q$.
If $\nestedin{Q}\cap\orthogonalto{P}=\emptyset$ there is nothing to
prove, so suppose $R\in \nestedin{Q}\cap\orthogonalto{P}$.
We claim that $(Q\cap P^\perp)^{ess}$ is a clean container for $P$ in $Q$.
It is a proper essential parabolic subgroup that commutes with
$P$ and is contained in $Q$.
It is not equal to $Q$, since $P\nested Q$ implies $P$ and $Q$ do not commute.
Furthermore, $R$ is contained in $Q$ and commutes with $P$, so $R\leq
Q\cap P^\perp$, so, by \fullref{lem:ess_monotone}, $R=R^{ess}\leq
(Q\cap P^\perp)^{ess}$.
It remains only to show that $(Q\cap P^\perp)^{ess}$ is actually a
domain. It is nontrivial, since it contains $R$, and is a proper essential
parabolic subgroup by construction.
To conclude it is a domain, observe, by \fullref{lem:essential_normalizer}, that $\normalizer_W((Q\cap
P^\perp)^{ess})=(Q\cap
P^\perp)^{ess}\times\centralizer_W((Q\cap
P^\perp)^{ess})$, which is wide, because the first factor contains $R$ and the
second factor contains $P$, so both factors are nonspherical.

A similar argument shows that  $(P^\perp)^{ess}=(W^{ess}\cap P^\perp)^{ess}$ is a clean container for $P$ in $\top$.

\subsubsection{Consistency}
\begin{definition}\label{def:relprojection}
  If $P$ is a nontrivial, proper, essential parabolic subgroup of
  $(W,S)$, recall that $\Upsilon_{N(P)}$ denotes the unique parabolic
  subcomplex of $\Davis$ with stabilizer equal to $\normalizer_W(P)$.

For $P\in\HHindex\setminus\{\top\}$ and $Q\in\HHindex$, define $\relprojection{P}{Q}=\proj{Q}(\Upsilon_{N(P)})$.
\end{definition}

When $Q$ is reducible the sets $\relprojection{P}{Q}$ have diameter 0.
When $\top$ is irreducible, the fact that $\relprojection{P}{\top}$ has diameter 1 in $\coordinate{\top}$ for all
  $P\nestedneq\top$ is the easy direction of
  \fullref{lem:essential_with_wide_normalizer_equals_unbounded_in_coneoff}.
The following handles the remaining case:
\begin{lemma}\label{lem:diameter_of_relative_projection}
For $P\transverse Q$ or $P\nestedneq Q\nestedneq\top$ with irreducible $Q$, the set
$\relprojection{P}{Q}$ has diameter at most 1 in $\coordinate{Q}$.
\end{lemma}
\begin{proof}
  First observe:
\[  (P\transverse Q)\lor(P\nestedneq Q)\implies (P\not\orthogonal
    Q)\land(Q\not\nested P)
  \]
  Pick $\Xi\in\fibers(Q)$, and let $\Upsilon_{N(P)}$ be the unique
element of $\fibers(\normalizer_W(P))$.
The canonical isomorphism between $\coordinate{Q}$ and $C(\Xi)$ takes
$\relprojection{P}{Q}$ to $\gate_\Xi(\Upsilon_{N(P)})$.
If $\gate_\Xi(\Upsilon_{N(P)})$ is all of $\Xi^{(0)}$ then, using \fullref{lem:properties_of_gate_maps}\eqref{item:gate_bridge}:
\[Q=\stab(\Xi)=\stab(\gate_{\Xi}(\Upsilon_{N(P)}))=Q\cap\normalizer_W(P)\leq
  \normalizer_W(P)=P\times P^\perp\]
Since $Q$ is irreducible, $Q\leq P\times P^\perp$ implies that either
 $Q\leq P$, which contradicts $Q\not\nested P$, or 
$Q\leq P^\perp$, which contradicts $P\not\orthogonal Q$.
Since both outcomes are contradictions, $\gate_\Xi(\Upsilon_{N(P)})$
is not all of $\Xi^{(0)}$. Rather, by
\fullref{lem:properties_of_gate_maps}\eqref{item:gate_residue}, it spans a proper parabolic subcomplex
of $\Xi$, so it is coned off in $C(\Xi)$.
\end{proof}

The consistency condition
\ref{ax:consistency}\ref{ax:consistency:small} is immediate:
$\left(P\nested Q\nestedneq R\right)\lor \left(R\not\orthogonal
  P\nested Q\transverse R\right)$ implies $\relprojection{P}{R}$ and
$\relprojection{Q}{R}$ are both defined, and that $P\leq Q$ are
parabolic subgroups. Choose $\Xi\in\fibers(Q)$. There exists
$\Upsilon\in\fibers(P)$ such that $\Upsilon\subset\Xi$.
Thus, the image of $\Upsilon$ is contained in $\relprojection{P}{R}\cap\relprojection{Q}{R}$, so $\domaindist{R}{\relprojection{P}{R}}{\relprojection{Q}{R}}=0$.

The consistency condition
\ref{ax:consistency}\ref{ax:consistency:Behrstock} is given by the
following, whose proof is a ``grid of grids'' argument inspired by the
proof of \cite[Theorem~2.16]{CapMar13}.
Recall that the wall structure of the Niblo-Reeves cube complex $\NR$
of $(W,S)$ is $W$--equivariantly isomorphic to that of $\Davis$.
Since $\NR$ is a CAT(0) cube complex, its walls have the Helly property.
It follows that the maximum size of a set of pairwise transverse walls
of $\Davis$ is $\dim(\NR)$. 
\begin{proposition}
  For all $x\in\Davis^{(0)}$ and $P,Q\in\HHindex$ with $P\transverse
  Q$:
  \[\min\{\domaindist{P}{x}{\relprojection{Q}{P}},\domaindist{Q}{x}{\relprojection{P}{Q}}\}\leq 1+2(2\caprace-1)\dim(\NR)\]
\end{proposition}
\begin{proof}
  If either of $P$ or $Q$ is reducible then the corresponding
  coordinate space is a single point, so the minimum in the statement
  is 0. Thus, we may assume that both $P$ and $Q$ are irreducible.
  
  Suppose the proposition is not satisfied, so both domain distances
  are strictly greater than the given bound.
  Choose $\Upsilon\in\fibers(P)$ and $\Xi\in\fibers(Q)$.
  Define $a:=\gate_\Upsilon(x)$, $b:=\gate_{\gate_\Upsilon(\Xi)}(x)$,
  $c:=\gate_\Xi(x)$, and $d:=\gate_{\gate_\Xi(\Upsilon)}(x)$.
  Then
  $d_{C(\Upsilon)}(a,b)\geq d_{C(\Upsilon)}(a,\gate_\Upsilon(\Xi))\geq\domaindist{P}{x}{\Xi}\geq\domaindist{P}{x}{\relprojection{Q}{P}}$,
  and, similarly,  $d_{C(\Xi)}(c,d)\geq\domaindist{Q}{x}{\relprojection{P}{Q}}$.

  Every wall $V\in\walls(a,b)$ separating $a$ from $b$ cuts
  $\Upsilon$, by convexity of $\Upsilon$.
  Furthermore, we claim that $V$ separates $x$ from $\Xi$. To see
  this, suppose $V$ does not separate $x$ and some
  $z\in\Xi$.
  Observe, by
  \fullref{lem:properties_of_gate_maps}\eqref{item:gate_nested}, that
  $b=\gate_{\gate_\Upsilon(\Xi)}(a)$, so by
  \fullref{lem:properties_of_gate_maps}\eqref{item:wall_gate}, $V$
  separates $a=\gate_\Upsilon(x)$ from $\gate_\Upsilon(\Xi)$. Since
  $\gate_\Upsilon(\Xi)$ contains $\gate_\Upsilon(z)$, 
   either $V$ separates $z$ from
  $\gate_\Upsilon(z)$ or $V$ separates $x$ from $\gate_\Upsilon(x)$.
  By \fullref{lem:properties_of_gate_maps}\eqref{item:wall_gate}, this
  implies that $V$ separates $z$ from $\Upsilon$ or $V$ separates $x$
  from $\Upsilon$.
  Both are impossible, since $V$ cuts $\Upsilon$.

  Similarly, every wall $H\in \walls(c,d)$ cuts $\Xi$ and separates
  $x$ from $\Upsilon$.

  If we label the complementary open halfspaces
  so that $x\in V^+\cap H^+$, then we have $a\in V^+\cap H^-$, $b,d\in
  V^-\cap H^-$, and $c\in V^-\cap H^+$. Since all four halfspace
  intersections are nonempty, $V\transverse H$.
  Since $V$ and $H$ were arbitrary, every wall of $\walls(a,b)$ is
  transverse to every wall of $\walls(c,d)$. 

Each wall of $\walls(a,b)$ cuts $\Upsilon$,
so $\pc(\walls(a,b))\leq P$.
If it were a proper parabolic subgroup of $P$ then $a$ and $b$ would be
contained in a proper parabolic subcomplex of $\Upsilon$.
That is impossible, because then $d_{C(\Upsilon)}(a,b)\leq 1$,
contrary to hypothesis.
Thus $\pc(\walls(a,b))=P$.
In fact, we can throw away some of the walls without the parabolic
closure getting smaller: suppose $\mathcal{M}'$ is a subset of
$\walls(a,b)$ for which $\pc(\walls(a,b)\setminus\mathcal{M}')\lneq
P$.
Take a combinatorial geodesic $\gamma$ from $a$ to $b$, which is
contained in $\Upsilon$, by convexity.
The edges of $\gamma$ are in bijection with $\walls(a,b)$, so $\mathcal{M}'$ gives a subdivision of $\gamma$ into a
concatenation of (possibly trivial) subsegments with all edges dual to
a wall of $\walls(a,b)\setminus\mathcal{M}'$, separated by single
edges dual to a wall of $\mathcal{M}'$.
By hypothesis, each of the former type of subpath has every reflection
of its dual walls in $\pc(\walls(a,b)\setminus\mathcal{M}')\lneq
P$. Thus, its endpoints are contained in a proper parabolic subcomplex
of $\Upsilon$, so they are at distance at most 1 in $C(\Upsilon)$.
This implies that:
\[ 1+2(2\caprace-1)\dim(\NR)<d_{C(\Upsilon)}(a,b)\leq 2|\mathcal{M}'|+1\]
Hence, $|\mathcal{M}'|>(2\caprace-1)\dim(\NR)$.
Conversely:
\begin{equation}
  \label{eq:1}
|\mathcal{M}'|\leq(2\caprace-1)\dim(\NR)\implies\pc(\walls(a,b)\setminus\mathcal{M}')=P  
\end{equation}

Consider the set of complementary open halfspaces of walls of
$\walls(a,b)$ containing $a$, partially ordered by inclusion. 
By Dilworth's Theorem, that partially ordered set can be partitioned
into a collection of disjoint chains, the number of which is at most
the maximal size of an antichain in the collection.
In this setup, chains of nested halfspaces correspond to chains of
walls of $\walls(a,b)$, and maximal antichains come from maximal
subsets of $\walls(a,b)$ of pairwise crossing walls.
The cardinality of any set of pairwise crossing walls of $\Davis$ is
at most $\dim(\NR)$.
Conclude that $\walls(a,b)$ can be partitioned into
at most $\dim(\NR)$ disjoint chains.

Take such a partition of $\walls(a,b)$ into at most $\dim(\NR)$ disjoint chains, and throw away any whose
length is at most $2\caprace-1$. Label the remaining chains
$\mathcal{V}_1$,\dots,$\mathcal{V}_m$.
Consider the set of walls $\mathcal{M}''$ in the interior of one of
the $\mathcal{V}_i$, and take
$\mathcal{M}':=\walls(a,b)\setminus\mathcal{M}''$.
Then the total number $|\mathcal{M}'|$ of walls discarded, either for
belonging to a short chain or for being near an end of a long chain, is:
\[|\mathcal{M}'|\leq (\dim(\NR)-m)\cdot(2\caprace-1)+m(2\caprace-2)\leq(2\caprace-1)\dim(\NR)\]
By \eqref{eq:1}, $\pc(\mathcal{M}'')=\pc(\interior{\mathcal{V}}_1,\dots,\interior{\mathcal{V}}_m)=P$.
This gives  $P=\pc(\interior{\mathcal{V}}_1,\dots,\interior{\mathcal{V}}_m)\leq
\pc({\mathcal{V}}_1,\dots,{\mathcal{V}}_m)\leq P$.

By a symmetric argument applied to $\walls(c,d)$, we get disjoint
chains $\mathcal{H}_1$,\dots,$\mathcal{H}_n$ in $\walls(c,d)$ that
each have length at least $2\caprace$ and such that
$\pc({\mathcal{H}}_1,\dots,{\mathcal{H}}_n)=\pc(\interior{\mathcal{H}}_1,\dots,\interior{\mathcal{H}}_n)=Q$.

For each $i,j$, $(\mathcal{V}_i,\mathcal{H}_j)$ is a grid with
$|\mathcal{V}_i|,|\mathcal{H}_j|\geq 2\caprace$, so by
\cite[Proposition~3.20]{CasChuTao}, either
$\pc(\mathcal{V}_i)=\pc(\mathcal{H}_j)$ is irreducible higher rank
affine, or $\pc(\interior{\mathcal{H}}_j)\leq\pc(\mathcal{V}_i)^\perp$.

First suppose the latter outcome holds for every pair of indices
$i,j$.
Then for all $i$, for all $j$,
$\pc(\interior{\mathcal{H}}_j)\leq\pc(\mathcal{V}_i)^\perp$, and since
$\pc(\mathcal{V}_i)^\perp$ is a parabolic subgroup,
$Q=\pc(\interior{\mathcal{H}}_1,\dots, \interior{\mathcal{H}}_n)\leq
\pc(\mathcal{V}_i)^\perp$.
Thus, for all $i$, $\pc(\mathcal{V}_i)\leq Q^\perp$.
But $Q^\perp$ is a parabolic subgroup, so
$P=\pc(\mathcal{V}_1,\dots,\mathcal{V}_m)\leq Q^\perp$.
This says $P$ and $Q$ commute, which is a contradiction, because that
means $P\orthogonal Q$, but we started with $P\transverse Q$.

The other option is that there is some pair of indices, which, up to
relabelling, we may assume are both 1, such that
$\pc(\mathcal{V}_1)=\pc(\mathcal{H}_1)$ is irreducible higher rank
affine.
For all $1<i\leq m$, consider the grid
$(\mathcal{H}_1,\mathcal{V}_i)$. By \cite[Proposition~3.20]{CasChuTao},
either $\pc(\interior{\mathcal{V}}_i)\leq \pc(\mathcal{H}_1)^\perp$,
or
$\pc(\interior{\mathcal{V}}_i)\leq\pc(\mathcal{V}_i)=\pc(\mathcal{H}_1)$.
Then $\pc(\mathcal{H}_1)\times\pc(\mathcal{H}_1)^\perp$ is a parabolic
subgroup that contains each $\pc(\interior{\mathcal{V}}_i)$, so contains the irreducible parabolic  subgroup
$P=\pc(\interior{\mathcal{V}}_1,\dots,\interior{\mathcal{V}}_m)$.
Since $P$ is irreducible it is contained in either $\pc(\mathcal{H}_1)$ or
$\pc(\mathcal{H}_1)^\perp$, but since we already know that $P$ has
nontrivial intersection with $\pc(\mathcal{H}_1)$, as 
$\pc(\mathcal{H}_1)=\pc(\mathcal{V}_1)\leq P$, it follows that $P=\pc(\mathcal{H}_1)\leq Q$, which implies $P\nested Q$,
which contradicts $P\transverse Q$.
\end{proof}

\subsubsection{Bounded geodesic image}
If $P\nestedneq Q$ and $P$ is reducible then $\coordinate{P}$ is a point, so
$\domaindist{P}{x}{y}>E$ never occurs. If $Q$ is reducible then
$\coordinate{Q}$ is a point, so the conclusion is immediate. If $P$ is
irreducible and $Q=\top$ is irreducible then, up to taking large enough
$E$, this is \cite[Corollary~5.25]{CasChuTao}.
Thus, the remaining case is that both $P$ and $Q$ are irreducible with
$Q\neq \top$:
\begin{proposition}\label{prop:pre-wide_BGI}
  There exist $B$ and $K$ such that if  $P\lneq Q$ are pre-wide
  parabolic subgroups, $x,y\in\Davis^{(0)}$, and $\gamma$ is a
  geodesic from $\proj{Q}(x)$ to $\proj{Q}(y)$ in $\coordinate{Q}$, then:
  \[\domaindist{Q}{\gamma}{\relprojection{P}{Q}}>B\implies\domaindist{P}{x}{y}\leq K\]
\end{proposition}
\begin{proof}
  Note that $Q$ is not higher rank affine, since it contains the proper
  nonspherical subparabolic $P$, so $\coordinate{Q}$ is
  $\epsilon$--hyperbolic, by
  \fullref{lem:nonaffine_domains_have_hyperbolic_coordinates}.
  Let $D$ be the Hausdorff distance bound of
  \fullref{lem:irreducible_coordinates_hyperbolic}.
We will show it suffices to take $K:=2D+2\epsilon$ and $B:=C+2\epsilon+1$
for $C:=2D+2\epsilon+2$.
  
  Choose $\Xi\in\fibers(Q)$.
  The coordinate space $\coordinate{Q}$ is canonically isomorphic to
  $C(\Xi)$, by a map that identifies $\proj{Q}(x)$ to $\gate_\Xi(x)$
  and $\relprojection{P}{Q}$ with
  $\relprojection{P}{\Xi}:=\bigcup_{\Upsilon\in\fibers(P)}\gate_\Xi(\Upsilon)$.
  Suppose that $x,y\in\Davis^{(0)}$ are vertices and $\gamma$ is a
  geodesic from $x':=\gate_\Xi(x)$ to $y':=\gate_\Xi(y)$ in $C(\Xi)$
  such that $d_{C(\Xi)}(\gamma,\relprojection{P}{\Xi})>B$.
Choose $\Upsilon\in\fibers(P)$ with
  $\Upsilon\subset\Xi$, and let
  $x'':=\gate_\Upsilon(x)=\gate_\Upsilon(x')$ and
  $y'':=\gate_\Upsilon(y)=\gate_\Upsilon(y')$, where the equalities
  are from \fullref{lem:properties_of_gate_maps}\eqref{item:gate_nested}.

  Let $\beta_x$ be a combinatorial geodesic in $\Davis^{(1)}$ from
  $x'$ to $x''$, and let $\beta_y$ be combinatorial geodesic from $y'$
  to $y''$.
  Let $\bar\beta_x$ be a $C(\Xi)$--geodesic from $x'$ to $x''$, and
  let $\bar\beta_y$ be a $C(\Xi)$--geodesic from $y'$ to $y''$.
  Let $\alpha$ be a $C(\Xi)$--geodesic from $x''$ to $y''$, so that
  $\alpha$, $\bar\beta_x$, $\bar\beta_y$, and $\gamma$ are the sides of
  a geodesic quadrilateral in the $\epsilon$--hyperbolic space $C(\Xi)$.
  That quadrilateral is $2\epsilon$--thin.

  Observe that $\Upsilon=\gate_\Xi(\Upsilon)\subset \relprojection{P}{\Xi}$,
  so $d_{C(\Xi)}(x',x'')\geq d_{C(\Xi)}(\gamma, \relprojection{P}{\Xi})>B$
  so there is a point $\bar b_x\in\bar\beta_x$ at distance $C$ from
  $x''$.
We claim that $\bar b_x$ is more than $2\epsilon$ from $\alpha$ and
$\gamma$, so it is $2\epsilon$--close to some point $\bar b_y\in\bar\beta_y$.
First, since $\diam_{d_{C(\Xi)}}(\relprojection{P}{\Xi})\leq 1$ and since
  $\alpha$ connects points $x''$ and $y''$ of $\relprojection{P}{\Xi}$, the
  length of $\alpha$ is at most $1$, so $d_{C(\Xi)}(\bar
  b_x,\alpha)\geq d_{C(\Xi)}(\bar b_x,x'')-1=C-1>2\epsilon$.
 Second, $d_{C(\Xi)}(\bar
 b_x,\gamma)\geq d_{C(\Xi)}(x'',\gamma)-d_{C(\Xi)}(x'',\bar b_x)>B-C>2\epsilon$.
  
By \fullref{lem:irreducible_coordinates_hyperbolic}, there
exist $b_x\in\beta_x$ and $b_y\in\beta_y$ such that $d_{C(\Xi)}(\bar
b_x,b_x)\leq D$ and $d_{C(\Xi)}(\bar b_y,b_y)\leq D$.
Observe, $d_{C(\Xi)}(b_x,b_y)\leq 2\epsilon+2D=K$.

Also, $d_{C(\Xi)}(b_x,\relprojection{P}{\Xi})\geq d_{C(\Xi)}(\bar
b_x,\relprojection{P}{\Xi})-D\geq d_{C(\Xi)}(\bar b_x,x'')-1-D\geq
C-1-D$, and, similarly, $d_{C(\Xi)}(b_y,\relprojection{P}{\Xi})\geq
C-1-D-2\epsilon$.
Let $\eta$ be a  $C(\Xi)$--geodesic from $b_x$ to $b_y$. 
It has length at most $K$, and for all $c\in\eta$:
\begin{align*}
  d_{C(\Xi)}(c,\relprojection{P}{\Xi})&\geq\frac{d_{C(\Xi)}(b_x,\relprojection{P}{\Xi})+d_{C(\Xi)}(b_y,\relprojection{P}{\Xi})-d_{C(\Xi)}(b_x,b_y)}{2}\\
  &\geq\frac{(C-1-D)+(C-1-D-2\epsilon)-K}{2}=1
\end{align*}

Suppose there is a cone vertex $v_\Psi$ on $\eta$, for some proper
parabolic subcomplex $\Psi$ of $\Xi$. 
If $P\leq\stab(\Psi)$ then $\Psi$ contains a $P$--fiber $\Upsilon'$. 
That is impossible, because then
$\Upsilon'\subset\relprojection{P}{\Xi}$, which would give a
contradiction:
\[1\leq d_{C(\Xi)}(v_\Psi,\relprojection{P}{\Xi})\leq d_{C(\Xi)}(v_\Psi,\Upsilon')=\sfrac{1}{2}\]
Thus, $P\not\leq\stab(\Psi)$.
By 
\fullref{lem:properties_of_gate_maps}\eqref{item:gate_residue} and
\eqref{item:gate_bridge},
$\stab(\gate_\Upsilon(\Psi))=\stab(\Upsilon)\cap\stab(\Psi)\leq
\stab(\Psi)$, so $\stab(\gate_\Upsilon(\Psi))\neq P$, so
$\gate_\Upsilon(\Psi)$ spans a proper parabolic subcomplex of $\Upsilon$.
Thus, $P\not\leq\stab(\Psi)$ implies
$\diam_{\delta_\Upsilon}(\Psi)\leq 1$, since $\gate_\Upsilon(\Psi)$ is
coned off in $C(\Upsilon)$. 

Consider the sequence $b_x=z_0,z_1,\dots,z_n=b_y$ of $\Davis^{(0)}$
vertices along $\eta$.
Each consecutive pair $z_i$, $z_{i+1}$ either are mutually adjacent to
a cone vertex $v_\Psi$ of $\eta$, or are an adjacent pair in
$\Davis^{(1)}$.
In the former case, the previous paragraph says
$\delta_{\Upsilon}(z_i,z_{i+1})\leq 1$. In the latter case, 
the fact that $\gate_\Upsilon$ is $1$--Lipschitz, gives the same
conclusion. 
Thus, $\gate_\Upsilon(\eta)$ has length in $C(\Upsilon)$ bounded above
by the length of $\eta$ in $C(\Xi)$, which was at most $K$.
Now, one endpoint $b_x$ of $\eta$ sits on $\beta_x$, which is a
combinatorial geodesic from $x'$ to $x''=\gate_\Upsilon(x')$, so
$\gate_\Upsilon(b_x)=x''$. Similarly, the other endpoint $b_y$ of
$\eta$ projects to $y''$.
Thus, $\delta_\Upsilon(x'',y'')\leq K$.
\end{proof}

\subsubsection{Large links}
Let $B$ and $K$ be common BGI constants satisfying both
\cite[Corollary~5.25]{CasChuTao} and  \fullref{prop:pre-wide_BGI}.
Let $C:=2B+3$. Assume $E>\max\{2K+2C+2,\,\sfrac{|S|}{2}\}$.
If $Q\in\HHindex$ is an irreducible higher rank affine parabolic
subgroup then it is $\nested$--minimal, so there is nothing to prove.

Suppose $Q\in\HHindex$ is reducible, and write $Q=Q_1\times\cdots
\times Q_n$, where $2\leq n\leq \lfloor \sfrac{|S|}{2}\rfloor\leq E$ and the $Q_i$ are
nonspherical irreducible parabolic subgroups. 
Since $Q$ is wide, each $Q_i$ is pre-wide, by \fullref{lem:pre-wide}, so each $Q_i$ is a domain
that is properly nested in $Q$. 
For each $i$, $\domaindist{Q}{x}{\relprojection{Q_i}{Q}}=0$, since $\coordinate{Q}$ is a point.
If $P\nestedneq Q$ is reducible then $\coordinate{P}$ is a point and
$\domaindist{P}{x}{y}=0$, so if  $\domaindist{P}{x}{y}\geq E$ then $P$ is
irreducible. An irreducible parabolic subgroup of $Q$ is contained in
one of the factors $Q_i$.
A similar proof applies for reducible $\top$ by replacing $Q$ by $W^{ess}$.
This completes the verification of the Large Links Axiom for reducible
domains. 

Now suppose $Q\in\HHindex$ is irreducible and not higher rank affine.
If $Q=\top$ then $\coordinate{Q}=\coneoff$. Otherwise, for any choice
of $\Xi\in\fibers(Q)$ there is a canonical identification of
$\coordinate{Q}$ with $C(\Xi)$. 
Let $\gamma$ be a geodesic in $\coordinate{Q}$ from $\proj{Q}(x)$ to
$\proj{Q}(y)$.
For each cone vertex $v_\Psi$ of $\gamma$, we claim that every
nonspherical irreducible factor in the canonical
product decomposition of $\stab(\Psi)$ is a domain. 
If there is more than one nonspherical irreducible factor, then those
factors are pre-wide, by \fullref{lem:pre-wide}, hence are domains. 
If $\stab(\Psi)$ has exactly one nonspherical irreducible factor $R$ and $Q=\top$
then for $v_\Psi$ to have occurred as a cone vertex, $R$ must be higher rank
affine, so $R$ is a domain. If $Q\neq\top$ then $\normalizer_W(Q)$ is
wide, so $\normalizer_W(R)$ is wide, hence $R$ is pre-wide, by
\fullref{lem:pre-wide}, so $R$ is a domain. 
Take $\mathcal{R}$ to be all nonspherical irreducible factors of
$\stab(\Psi)$, as $v_\Psi$ ranges over the cone vertices on $\gamma$. 
Then $|\mathcal{R}|\leq
\domaindist{Q}{x}{y}\cdot\lfloor \sfrac{|S|}{2}\rfloor$. This is less than $N$,
since $E\geq\lfloor \sfrac{|S|}{2}\rfloor$.

Each $R\in\mathcal{R}$ is properly nested in $Q$. This is clear if
$Q=\top$. If $Q\neq\top$ then $R\leq\stab(\Psi)$ for some cone vertex
$v_\Psi\in\coordinate{Q}=C(\Xi)$, so $\Psi$ is proper parabolic
subcomplex of $\Xi$, so $R\leq\stab(\Psi)\lneq Q$.

We check the location of relative projections: if $R\in\mathcal{R}$ arose
as a nonspherical irreducible factor of $\stab(\Psi)$, then choose
$\Upsilon\in\fibers(R)$ with $\Upsilon\subset\Psi$. Then
$\Upsilon\subset\relprojection{R}{Q}$ and every vertex of $\Upsilon$
is at distance $\sfrac{1}{2}$ from $v_\Psi\in\coordinate{Q}$.
Thus,
$\domaindist{Q}{x}{\relprojection{R}{Q}}\leq\domaindist{Q}{x}{v_\Psi}+\sfrac{1}{2}\leq\domaindist{Q}{x}{y}+\sfrac{1}{2}\leq
N$.

Finally, suppose $P\in\nestedin{Q}\setminus\{Q\}$ and
$\domaindist{P}{x}{y}\geq E$.
This implies $P$ is irreducible, since otherwise
$\domaindist{P}{x}{y}=0$.
Since $P\in\HHindex$ is irreducible it is pre-wide.
Since $E>K$, BGI implies $\gamma$ passes within distance $B$ of
$\relprojection{P}{Q}$ in $\coordinate{Q}$.
Let $p\in\relprojection{P}{Q}$ and $q\in\gamma$ be points such that
$\domaindist{Q}{p}{q}\leq B$.
If $\domaindist{Q}{x}{q}<C$, set $u:=\proj{Q}(x)$.
Otherwise, take $u$ to be a vertex of $\Davis^{(0)}$ on the subsegment
of $\gamma$ between $\proj{Q}(x)$ and $q$ such that
$C\leq\domaindist{Q}{q}{u}\leq C+1$.
Similarly, take $v$ to be a $\Davis^{(0)}$ vertex on the subsegment of
$\gamma$ between $q$ and $\proj{Q}(y)$ that is either at distance
between $C$ and $C+1$ from $q$, or equal to $\proj{Q}(y)$ if the
subsegment is too short.
If $\gamma_{uv}$ is the segment of $\gamma$ from $u$ to $v$:
\begin{equation}
  \label{eq:4}
  |\gamma_{uv}|_{\coordinate{Q}}\leq 2C+2
\end{equation}

Suppose $u\neq\proj{Q}(x)$.
For every point $u'$ on the subsegment of $\gamma$ between
$\proj{Q}(x)$ and $u$ and every $p'\in\relprojection{P}{Q}$:
\[\domaindist{Q}{u'}{p'}\geq\domaindist{Q}{u'}{q}-\domaindist{Q}{q}{p}-\domaindist{Q}{p}{p'}\geq C-B-1> B\]
Thus, BGI gives $\domaindist{P}{x}{u}\leq K$ when $u\neq\proj{Q}(x)$.
If $u=\proj{Q}(x)$ then $\domaindist{P}{x}{u}=0\leq K$.

Symmetrically, $\domaindist{P}{y}{v}\leq K$. Thus:
\begin{equation}
  \label{eq:3}
  \domaindist{P}{u}{v}\geq\domaindist{P}{x}{y}-2K\geq E-2K>2C+2
\end{equation}
Now we claim that $P\leq\stab(\Psi)$ for some cone vertex $v_\Psi$ of
$\gamma_{uv}$.
Otherwise, as in the proof of \fullref{prop:pre-wide_BGI}, we find
that $\gamma_{uv}$ has $\delta_P$--diameter at most equal to its
 $\coordinate{Q}$--length, which gives a contradiction between
\eqref{eq:4} and \eqref{eq:3}.
Thus, $P\leq\stab(\Psi)$.
Since $P$ is nonspherical and irreducible, it is contained in some
nonspherical irreducible factor $R$ of the canonical product decomposition
of $\stab(\Psi)$. By construction, $R\in\mathcal{R}$, so this
completes the verification of the Large Links Axiom.
   
\subsubsection{Partial realization}
In Coxeter groups the Partial Realization Axiom, with $E=0$, follows
easily from:
\begin{lemma}\label{lem:simultaneous_conjugation}
  If $P_1$,\dots,$P_n$ are pairwise commuting essential parabolic
  subgroups then there exist $w\in W$ and disjoint $T_i\subset S$ such that
  $P_i=wW_{T_i}w^{-1}$ for all $i$. 
\end{lemma}
\begin{proof}
  Since $P_i$ is essential, its canonical product decomposition
  contains only nonspherical irreducible components. These commute
  with each other and with every other $P_j$, so by replacing each
  $P_i$ with its factors, we may assume that the original collection
  contains only nonspherical irreducible parabolic subgroups.

  The proof is by induction on $n$.
  Since $P_1$ is parabolic, there exist $v\in W$ and $T_1\subset S$
  such that $P_1=vW_{T_1}v^{-1}$.
  If $n=1$ we are done, by taking $w:=v$.
  Otherwise, since $P_1$ is nonspherical and irreducible, $\centralizer_W(P_1)$
  is the parabolic subgroup $P_1^\perp=vW_{T_1^\perp}v^{-1}$.
  Since the subgroups $P_2$,\dots,$P_n$ commute with $P_1$, they are
  all in $P_1^\perp$, so for $i>1$ we have that
  $v^{-1}P_2v$,\dots,$v^{-1}P_nv$ is a collection of $n-1$ pairwise
  commuting parabolic subgroups of $W_{T_1^\perp}$.
  By induction, there is a common conjugator $u\in W_{T_1^\perp}$ and
  disjoint subsets
  $T_2,\dots,T_n\subset T_1^\perp\subset S$ such that
  $uv^{-1}P_ivu^{-1}=W_{T_i}$ for all $2\leq i\leq n$.
  Since $u\in W_{T_1^\perp}$, we have
  $uv^{-1}P_1vu^{-1}=uW_{T_1}u^{-1}=W_{T_1}$.
  This proves the lemma, by taking $w:=vu^{-1}$.
\end{proof}
Suppose $P_1$,\dots,$P_n$ is a collection of pairwise orthogonal
domains.
If some $P_i=\top$ then $i=n=1$, since
no domain is orthogonal to $\top$.
If $\top$ is reducible then $\coordinate{\top}$ is a point, so any $x$
works.
If $\top$ is irreducible then $\coordinate{\top}=\coneoff$.
Every $p\in\proj{\top}(\Davis^{(0)})$ is itself a vertex in
$\Davis^{(0)}$, since the projection map in this case is just the
identity on $\Davis^{(0)}$, so take $x=p$.
Then $\domaindist{\top}{p}{x}=0$.
There are no domains transverse to $\top$ or in which $\top$ is
properly nested, so the second realization condition is vacuous for $\top$.

Otherwise, all $P_i$ are essential parabolic subgroups. 
By \fullref{lem:simultaneous_conjugation}, there exist $w\in W$ and
disjoint $T_i\subset S$ such that $P_i=wW_{T_i}w^{-1}$ for all $i$.
For each $i$ set $\Upsilon_i:=w\Davis_{T_i}\in\fibers(P_i)$.
If $P_i$ is irreducible then $\coordinate{P_i}$ is canonically
isomorphic to $C(\Upsilon_i)$ and $\proj{P_i}(\Davis^{(0)})$ is identified
with $\Upsilon_i^{(0)}$, so $p_i$ is a vertex of
$\Upsilon_i\subset\Davis$, hence $p_i=wa_i$ for some $a_i\in W_{T_i}$.
If $P_i$ is reducible set $a_i:=1$.
Define $x=wa_1\cdots a_n\in W$.
The $a_i$ all commute, since
the subgroups $W_{T_i}$ pairwise commute, so $b_i:=a_1\cdots
a_na_i^{-1}\in W_{T_i^\perp}$.
Let $\Xi_i:=wb_i\Davis_{T_i}\in\fibers(P_i)$.
Then we have $x=wb_ia_i\in\Xi_i$, and we use this to verify the two
realization conditions.

If $P_i$ is reducible then Condition~\ref{ax:realization}\ref{ax:realization:location} is immediate from the fact that
$\coordinate{P_i}$ is a single point, so suppose $P_i$ is
irreducible.
We have points in two
$P_i$--fibers, $x\in\Xi_i$ and $p_i=wa_i\in\Upsilon_i$.
But $x=wb_iw^{-1}p_i$, so
\fullref{lem:properties_of_gate_maps}\eqref{item:gate_perp} says $\gate_{\Upsilon_i}(x)=\gate_{\Upsilon_i}(p_i)=p_i$.
Thus, $\domaindist{P_i}{p_i}{x}=0$.

For Condition~\ref{ax:realization}\ref{ax:realization:consistency}, if $Q\in\HHindex$ with $(P_i\nestedneq Q)\lor(P_i\transverse Q)$ then $\proj{Q}(x)\in\proj{Q}(\Xi_i)\subset\relprojection{P_i}{Q}$, so $\domaindist{Q}{x}{\relprojection{P_i}{Q}}=0$.

  \subsubsection{Uniqueness}
Fix $r$. Take $M(r):=M(r,\dots,r)$ as in
\cite[Proposition~5.23]{CasChuTao}.

Suppose that $\gamma$ is a combinatorial geodesic in $\Davis^{(1)}$
between $x,y\in\Davis^{(0)}$, with $\dcomb(x,y)>M(r)$.

If $\gamma$ contains a subsegment $\gamma'$ of length greater than
$M(r)$ that
is contained in a wide parabolic subcomplex $\Upsilon$ then every wall
dual to an edge of $\gamma'$ cuts $\Upsilon$ and separates $x$ and
$y$, by \fullref{lem:properties_of_gate_maps}\eqref{item:wall_gate}, so
$d_{\Upsilon}(x,y)\geq |\gamma'|> M(r)$.
By \cite[Proposition~5.23]{CasChuTao}, there exists a pre-wide
$P\leq\stab(\Upsilon)$ such that $\domaindist{P}{x}{y}>r$.
We have already excluded the case that $W^{ess}$ is higher rank
affine, so every pre-wide parabolic subgroup belongs to $\HHindex$.
This case necessarily holds when $W$ is wide, since then $\gamma$ is
contained in the wide parabolic subcomplex $\Davis$. 

If no subsegment of $\gamma$ of length at least $M(r)+1$ is contained in a
wide parabolic subcomplex, then
\cite[Theorem~3.31]{CasChuTao} says $\gamma$ is Morse,
with Morse gauge depending only on $M(r)$, and then
\cite[Corollary~4.20]{CasChuTao} says that
the inclusion of $\gamma$ into $\coneoff=\coordinate{\top}$ is
quasigeodesic, with quasigeodesic constants determined by $M(r)$.
In particular, $\domaindist{\top}{x}{y}$ is bounded below by a linear
function of $\dcomb(x,y)$, with constants determined by $M(r)$.
Thus, $\domaindist{\top}{x}{y}\geq r$ once $\dcomb(x,y)$ is
sufficiently large.

\subsection{Verifying coarse equivariance}\label{sec:verify_equivariance}
So far we have verified that $\Davis^{(1)}$, which is the Cayley graph of
$W$ with respect to $S$, is a relatively hierarchically hyperbolic
\emph{space}.
To say that $W$ is a relatively hierarchically hyperbolic
\emph{group}, we must verify the additional requirements of \fullref{def:relHHG}.
Define the action of $W$ on $\HHindex$ as follows.
For all $w\in W$ and $P\in\HHindex\setminus\{\top\}$,
define $w.\top:=\top$ and $w.P:=wPw^{-1}$.
Cofiniteness follows because each $P\in\HHindex\setminus\{\top\}$ is
conjugate to a special subgroup $W_T$, and there are only
finitely many subsets $T\subset S$.
Conjugation of parabolic subgroups preserves inclusion and
commutation, so the action on $\HHindex$ preserves nesting and
orthogonality.

The coordinate isometries are induced by the isometric action of
$W$ on $\Davis$. For reducible $P\in\HHindex$ there is nothing to
prove, since the coordinate space $\coordinate{P}$ is a single point.
For irreducible $P$ and $\Upsilon\in\fibers(P)$, let
$\phi_\Upsilon\from C(\Upsilon)\to \coordinate{P}$ be the canonical
isomorphism.
Since the $g$--action on $\Davis$ takes parabolic subcomplexes to
parabolic subcomplexes, it induces an isomorphism $g_\Upsilon^*\from
C(\Upsilon)\to C(g\Upsilon)$.
Define $g_P\from \coordinate{P}\to\coordinate{g.P} : \phi_{g\Upsilon}\circ
g_\Upsilon^*\circ \phi_\Upsilon^{-1}$.
This is independent of the choice of $\Upsilon$, because any other 
element of $\fibers(P)$ is of the form $z\Upsilon$ for $z\in P^\perp$.
But $\phi_{z\Upsilon}=\phi_{\Upsilon}\circ (z^{-1})^*_{z\Upsilon}$ and
$\phi_{gz\Upsilon}=\phi_{g\Upsilon}\circ (gz^{-1}g^{-1})^*_{gz\Upsilon}$. Thus:
\begin{align*}
  \phi_{gz\Upsilon}\circ g_{z\Upsilon}^*\circ
  \phi_{z\Upsilon}^{-1}&=\phi_{g\Upsilon}\circ (gz^{-1}g^{-1})^*_{gz\Upsilon}\circ    g^*_{z\Upsilon}     \circ
                         z^*_{\Upsilon}\circ \phi_\Upsilon^{-1}\\
&=\phi_{g\Upsilon}\circ g^*_{\Upsilon}\circ \phi_\Upsilon^{-1}
\end{align*}

Composability of these isometries follows because for
$\Upsilon\in\fibers(P)$ we have $h\Upsilon\in\fibers(hPh^{-1})=\fibers(h.P)$ and
$gh\Upsilon\in\fibers(ghPh^{-1}g^{-1})=\fibers((gh).P)$, so:
\begin{align*}
  (gh)_P&:=\phi_{gh\Upsilon}\circ
  (gh)^*_\Upsilon\circ\phi^{-1}_\Upsilon\\
  &=\phi_{gh\Upsilon}\circ
    g^*_{h\Upsilon}\circ\phi^{-1}_{h\Upsilon}\circ
    \phi_{h\Upsilon}\circ h^*_\Upsilon\circ \phi^{-1}_\Upsilon\\
  &=g_{h.P}\circ h_P
\end{align*}

Similar arguments apply to irreducible $\top$, since the $W$--action on
$\Davis$ induces an isometric action on $\coneoff=\coordinate{\top}$.

The coarse equivariance properties of \fullref{def:relHHG} likewise
follow from tracing through the definitions and using the $W$--action
on $\Davis$. In fact, they are equivariant on the nose: we can take $E=0$.



\section{Promoting the Coxeter relative HHG structure}\label{sec:virtual_HHG}
In this section we extend the relative HHG structure provided by
\fullref{thm:Cox_rHHG}.
As described in the introduction, we have two different approaches to
arranging the necessary equivariance conditions, one by passing to a
finite index subgroup of $W$ and one by restricting the affine
parabolics. 

We start by choosing a finite index subgroup of $W$ and showing
it admits an HHG structure. The structure has an additional property
that is needed to deduce \fullref{thm:Cox_quasicubical} from \fullref{thm:Cox_vHHG}:
\begin{definition}[{\cite{Pet26}}]\label{def:colorable}
  An HHG $(G,\HHindex)$ is \emph{colorable} if there exists an
  assignment of finitely many colors to the domains of $\HHindex$ such
  that every pair of distinct domains with the same color is transverse, and
  such that the $G$--action on $\HHindex$ induces a permutation action
  of $G$ on the set of colors. 
\end{definition}

In the proof of \fullref{thm:Cox_vHHG} we start with a relative HHS structure
$\HHindex$ whose peripheral coordinate spaces are HHSs,
build a ``proto-HHS structure'' $\HHindex^+$ by replacing the
peripheral domains by their corresponding coordinate HHS structures, and then complete the proto-HHS
to a full HHS structure $\HHindex^{++}$ by adding ``dummy
container domains''. This is a general relative HHS promotion recipe,
whose details appear in \fullref{sec:promotion}.
The main work is to arrange that the
peripheral stabilizers respect their peripheral HHS structures, so that
the promotion can be carried out equivariantly, as in
\fullref{prop:stabilizer_consistency}. We achieve this by passing to a
congruence subgroup of $W$.
Separate arguments establish cofiniteness of the action on
$\HHindex^{++}$ and colorability.

\medskip

Let $(W,S)$ be a finite rank Coxeter system with $S=\{s_1,\dots,s_n\}$
and Coxeter matrix $M:=(m_{ij})$. Let $C:=(-\cos(\sfrac{\pi}{m_{ij}}))=(c_{ij})$
be its cosine matrix. The \emph{Tits form} $\Titsform$ on
$V:=\mathbb{R}^{|S|}$ is the $W$--invariant symmetric bilinear form
$\Titsform(\vec{e}_i,\vec{e}_j)=c_{ij}$, where $\vec{e}_i$ is the
basis vector of $V$ corresponding to $s_i\in S$.
The \emph{Tits representation}\footnote{The Tits representation is
  called the ``canonical representation'' in \cite{davisbook} and the
  ``geometric representation'' in \cite{Bou02}.} is the faithful representation $\canonicalrep\from W\to
\GL(V)$ defined by $\canonicalrep(s_i):=\vec{v}\mapsto \vec{v}-2\Titsform(\vec{v},\vec{e}_i)\vec{e}_i$.

Let $R:=\mathbb{Z}[2c_{ij}\mid 1\leq i,j\leq n]$ be the ring generated
by the finitely many numbers $2c_{ij}$.
The numbers $2c_{ij}$ are algebraic integers, so $R$ is a finitely
generated $\mathbb{Z}$--module, and an order in the number field
$\mathbb{Q}(2c_{ij}\mid 1\leq i,j\leq n)$ contained in its ring of integers.
Thus, for any nontrivial ideal $\mathcal{I}$ in $R$, for example a principal ideal $\mathcal{I}=(a)=\{ar\mid r\in R\}$ generated by some $a\in R$, the quotient $R/\mathcal{I}$ is a finite ring.
The matrix of $\canonicalrep(s_i)$ has entries in $R$, 
so we get a finite index normal subgroup of $W$:
\begin{equation}
  \label{eq:congruence_subgroup}
  K_\mathcal{I}:=\ker(W\xrightarrow{\canonicalrep}\GL_{|S|}(R)\to\GL_{|S|}(R/\mathcal{I}))
\end{equation}

The following argument is extracted from the proof of the well-known
result that arithmetic groups have a finite index ``neat'' subgroup. We give a
proof for completeness. 
\begin{lemma}\label{lem:no_torsion}
  Let $p$ be an odd prime. Elements of $\ker\left(\GL_{n}(R)\to\GL_{n}(R/(p))\right)$ have no nontrivial root of
  unity eigenvalues.
\end{lemma}
\begin{proof}
  Suppose $A\in \ker\left(\GL_{n}(R)\to\GL_{n}(R/(p))\right)$.
    There exists $B\in M_n(R)$ such that $A=I_n+pB$.
   If $\zeta$ is an eigenvalue of $A$ then $\beta:=\frac{1-\zeta}{p}\in\mathbb{Q}(\zeta)$ is an eigenvalue of $-B$.
 Since $\beta$ is integral over $R$, which is integral over
 $\mathbb{Z}$, $\beta$ is an algebraic
  integer, so its field norm $N_{\mathbb{Q}(\zeta)/\mathbb{Q}}(\beta)$ is an
  integer.
    Multiplicativity of the norm gives:
   \[1-\zeta=p\beta\implies
    N_{\mathbb{Q}(\zeta)/\mathbb{Q}}(1-\zeta)=N_{\mathbb{Q}(\zeta)/\mathbb{Q}}(p)
    N_{\mathbb{Q}(\zeta)/\mathbb{Q}}(\beta)\]
  Thus, $N_{\mathbb{Q}(\zeta)/\mathbb{Q}}(p)$ divides
  $N_{\mathbb{Q}(\zeta)/\mathbb{Q}}(1-\zeta)$.
  If $\zeta$ is a primitive $m$--th root of unity, for some $m>1$,
  then  $N_{\mathbb{Q}(\zeta)/\mathbb{Q}}(p)=p^{\varphi(m)}$ divides $N_{\mathbb{Q}(\zeta)/\mathbb{Q}}(1-\zeta)=\Phi_m(1)$.
  Since $\Phi_m(1)$ is 1 or prime, we must have $\varphi(m)=1$, which
  implies $m=2$.
  But $\Phi_2(1)=2$, so $p=2$,
  contradicting that $p$ is odd. 
\end{proof}
Let $P=W_T$ be an irreducible affine special subgroup of $(W,S)$. It
splits as $W_T=\Lambda_{W_T}\rtimes F_{W_T}$, where
$\Lambda_{W_T}\cong\mathbb{Z}^{|T|-1}$ is the translation subgroup and
$F_{W_T}$ is the finite linear part.
In crystallographic group theory $F_{W_T}$ is called the \emph{point group}.
If $V_T\leq V$ is the subspace spanned by the basis vectors
corresponding to elements of $T$, the restriction of $\Titsform$ to
$V_T$ is positive semidefinite with one-dimensional radical,
\cite[Theorem~6.8.12(ii)]{davisbook}.
The $W_T$--action preserves the radical. The $W_{T^\perp}$--action
fixes $V_T$ pointwise.
Thus, the restriction of the Tits
representation to $\normalizer_W(W_T)=W_T\times W_{T^\perp}$ descends to a representation:
\[\bar\canonicalrep\from \normalizer_W(W_T)\to \GL(V_T/\operatorname{rad}(V_T))\]
Moreover, $\ker(\bar\canonicalrep)=\Lambda_{W_T}\times W_{T^\perp}$ and $\bar\canonicalrep(\normalizer_W(W_T))\cong F_P$.

\begin{lemma}\label{lem:cong_tran}
  Let $(W,S)$ be a finite rank Coxeter system, and let $\mathfrak{P}$ be its set of irreducible higher
  rank affine parabolic subgroups.
  For any odd prime $p$, let $K_{(p)}\leq W$ as in \eqref{eq:congruence_subgroup}.
  For every $P\in\mathfrak{P}$:
  \[K_{(p)}\cap\normalizer_W(P)\leq \ker(\normalizer_W(P)\to
    \normalizer_W(P)/P^\perp\cong P\to P/\Lambda_P=F_P)\]
In particular, $K_{(p)}\cap \normalizer_W(P)\leq \Lambda_P\times P^\perp$. 
\end{lemma}
\begin{proof}
  Since $K_{(p)}$ is normal, we may assume $P=W_T$ is special. 
  Suppose $g=\tau f z\in K_{(p)}\cap\normalizer_W(W_T)$ with
  $\tau\in\Lambda_{W_T}$, $f\in F_{W_T}$, and $z\in W_{T^\perp}$.
Then $\bar\canonicalrep(g)=\bar\canonicalrep(f)$ is a finite order matrix.
If $f$ has order $m$ then every eigenvalue $\lambda$ of
$\bar\canonicalrep(g)$ satisfies $\lambda^m=1$.
Suppose $m>1$.
  Then  $\bar\canonicalrep(g)=\bar\canonicalrep(f)$ is a nontrivial
  matrix of finite order $m$ dividing $|F_{W_T}|$.
  Since $\bar\canonicalrep(g)$ is nontrivial, not all of its
  eigenvalues are 1, so it has some eigenvalue of order
  $d>1$ with $d$ dividing $m$.
  The flag $0<\operatorname{rad}(V_T)<V_T<V$ is $\normalizer_W(W_T)$--invariant, 
so every eigenvalue of $\bar\canonicalrep(g)$ is an eigenvalue of
$\canonicalrep(g)$.
On the other hand, for $g\in K_{(p)}$, \fullref{lem:no_torsion} says $\canonicalrep(g)$ has no nontrivial root of unity eigenvalues. 
Thus, $f=1$, so $g\in\Lambda_{W_T}\times W_{T^\perp}$.
\end{proof}

\begin{proof}[Proof of \fullref{thm:Cox_vHHG}]
  If $W$ is hyperbolic it is already an HHG. If $W^{ess}$ is
  irreducible higher rank affine then $W$ is virtually Abelian, so it
  is virtually an HHG. Both of these cases are trivially colorable.
  Assume we are not in one of these cases.
  Then the HHG structure $\HHindex$ for $W$ relative to the set
  $\mathfrak{P}$ of the irreducible higher rank affine parabolic subgroups of
  $(W,S)$ given by \fullref{thm:Cox_rHHG} has $\HHindex\neq\{\top\}$.

  For $P\in\mathfrak{P}$, let $\Lambda_P\cong\mathbb{Z}^{\rank(P)-1}$ be its translation subgroup and
    let $F_P$ be its point group.
    Let $\mathfrak{Q}$ be the
    special subgroups in $\mathfrak{P}$.
    Pick any odd prime $p$ and take the congruence subgroup
    $K_{(p)}\leq W$
    as defined above. 

    Since $K_{(p)}$ is a finite index subgroup of $W$, the relative HHG
    structure on $W$ transfers to a relative HHG structure on $K_{(p)}$
    via the $K_{(p)}$--equivariant quasiisometry $K_{(p)}\into W$.
   Our goal is to choose HHS structures on representatives of
   $K_{(p)}$--orbits of peripherals, and show that the $K_{(p)}$--stabilizers
  of those representatives satisfy the hypotheses of \fullref{prop:stabilizer_consistency}.

  Since the action of $W$ on $\mathfrak{P}$ is by conjugation,
  $\stab_W(P)=\normalizer_W(P)$.
  Thus, \fullref{lem:cong_tran} gives
  $\stab_{K_{(p)}}(P)=K_{(p)}\cap\normalizer_W(P)\leq \Lambda_P\times P^\perp$.
  If $g\in\stab_{K_{(p)}}(P)$ then we can write $g=yz$ for $y\in \Lambda_P$ and
  $z\in P^\perp$.
  Since $P^\perp$ acts trivially on $\coordinate{P}=C(P)$, the action
  of $g$ on $\coordinate{P}$ is identical to the action of $y$ on
  $\coordinate{P}$.
  Thus, the $\stab_{K_{(p)}}(P)$--action on $\coordinate{P}$ factors
  through the $\Lambda_P$--action.

For each $K_{(p)}$--orbit in $\mathfrak{P}$, choose a representative $P$.
Define the standard product HHS structure on $\Lambda_P\cong
\mathbb{Z}^{\rank(P)-1}$, obtained as the product of $\rank(P)-1$
$\mathbb{Z}$ factors. This structure is finite and $\Lambda_P$--invariant.
There are only finitely many such $P$, because there are finitely many
$W$--orbits and $[W:K_{(p)}]<\infty$.

Since $\Lambda_P$ has finite index in $P$ and $P$ is equivariantly
quasiisometric to $C(P)=\coordinate{P}$, an orbit map $q$ of $\Lambda_P$ into $\coordinate{P}$ is
a $\Lambda_P$--equivariant quasiisometry.
Let $\mathfrak{T}_P$ be the HHS structure on $\coordinate{P}$ obtained
by transferring the product HHS structure on $\Lambda_P$ via the
quasiisometry; that is, keep the same domains and coordinate spaces,
and pre-compose projections with a chosen inverse quasiisometry $\bar q$.
Since 
$q$ is $\Lambda_P$--equivariant, $\bar q$ is coarsely
$\Lambda_P$--equivariant.
Moreover, we have only used finitely many $P$, so the HHS constants
are uniform.

Now we can verify the hypotheses of
\fullref{prop:stabilizer_consistency}.
Hypothesis \eqref{pp:action} is true because the $\stab_{K_{(p)}}(P)$--action on
$\mathfrak{T}_P$ factors through the $\Lambda_P$--action, which fixes
every domain of $\mathfrak{T}_P$.
If $g=yz\in\stab_{K_{(p)}}(P)$ for $y\in\Lambda_P$ and $z\in
P^\perp$ then we define the isometry $g_U\from
\coordinate{U}\to\coordinate{g.U}=\coordinate{U}$ required by
\eqref{pp:isometries} as follows.
If $\coordinate{U}$ is a single point then $g_U$ is the identity map.
Otherwise, $\coordinate{U}=\mathbb{R}$ and $g_U$ is the translation
given by the $U$--coordinate of $y$. 
Hypothesis \eqref{pp:composition} follows immediately.

To see that projections are coarsely $\stab_{K_{(p)}}(P)$--equivariant,
hypothesis \eqref{pp:projection}, recall that for domain $U$ with
nonsingleton coordinate space,  the projection $\proj{U}$ was
defined by precomposing a coordinate projection
$\proj{U}^{\Lambda_P}$ in $\Lambda_P$ with  the quasiisometry $\bar
q\from \coordinate{P}\to\Lambda_P$.
Thus, $g_U\circ \proj{U}=y_U\circ \proj{U}^{\Lambda_P}\circ \bar
q=\proj{U}^{\Lambda_P}\circ y\circ\bar q$,
while $\proj{g.U}\circ g_P=\proj{U}^{\Lambda_P}\circ\bar q\circ y_P$.
Since $\bar q$ is coarsely $\Lambda_P$--equivariant, these two
uniformly coarsely agree. 

Hypothesis~\eqref{pp:relprojection} is trivial, because in
$\mathfrak{T}_P$ only the $\nested$--minimal domains have nonsingleton
coordinate spaces, and those domains are pairwise orthogonal, so
$(U\nestedneq V)\lor(U\transverse V)$ implies $\coordinate{V}$ is a
point.

\fullref{prop:stabilizer_consistency} gives a $K_{(p)}$--HHS
structure $\HHindex^{++}$ on $K_{(p)}$ built from the restriction of
$\HHindex$ to $K_{(p)}$, promoted to a proto-HHS structure $\HHindex^+$ by
adding the standard product HHG structures on the
translation subgroups of the irreducible higher rank affine parabolic
subgroups, and then completed to $\HHindex^{++}$ by adding dummy
container domains. To conclude that $\HHindex^{++}$ is actually an HHG
structure, not just a $K_{(p)}$--HHS, we show that the $K_{(p)}$--action
satisfies the hypothesis of \fullref{cor:dummy_cofiniteness}.

The $K_{(p)}$--action on $\HHindex^+$ is cofinite, so we show that for
each $R\in\HHindex^+$ the $\stab_{K_{(p)}}(R)$--action on the set of
nonmaximal sets $\mathfrak{O}$ of pairwise orthogonal domains in
$\HHindex^+_R\setminus\{R\}$ is cofinite.

The peripheral HHG structures $\mathfrak{T}_P$ are finite, so if
$R\in\mathfrak{T}_P$ we are done, so we may assume
$R\in\HHindex\setminus\mathfrak{P}$.
Let
$\{Q_1,\dots,Q_n\}:=\mathfrak{O}\cap(\HHindex\setminus\mathfrak{P})$.
Note that $Q_j\nestedneq R$, so $Q_j\neq\top$, so $Q_j$ is some
essential parabolic subgroup of $W$. 
Let $P_1,\dots,P_m$ be the peripheral domains such that
$\mathfrak{O}$ contains some domain from $\mathfrak{T}_{P_i}$.
Then $\mathfrak{O}':=\{P_1,\dots,P_m,Q_1,\dots,Q_n\}$ is a family of pairwise
commuting essential parabolic subgroups of $W$.
Since the $\mathfrak{T}_P$ are finite, the map
$\mathfrak{O}\to\mathfrak{O}'$ is finite-to-one, so it suffices to
show that $\stab_{K_{(p)}}(R)$ acts cofinitely on sets of the form
$\mathfrak{O}'$. 

If $R=\top$ then $\stab_{K_{(p)}}(R)=K_{(p)}$, and
\fullref{lem:simultaneous_conjugation}, says
$\mathfrak{O}'$ can be simultaneously conjugated in
$W$ to a collection of pairwise commuting special subgroups.
There are only finitely many special subgroups, so only finitely many
subsets of pairwise commuting ones, so there are finitely many
$W$--orbits of sets of the form of $\mathfrak{O}'$.
Since $K_{(p)}$ has finite index in $W$, there are also finitely many
$K_{(p)}$--orbits. 

If $R\neq\top$ then $R\leq \stab_W(R)$.
The argument is the same as in the $\top$--case, applying
\fullref{lem:simultaneous_conjugation} with $R$ as the ambient Coxeter
group. Since $K_{(p)}\cap R$ has finite index in $R$, the
$\stab_{K_{(p)}}(R)$--action is also cofinite.

Finally, we show that $\HHindex^{++}$ is colorable.
We will actually produce a $K_{(p)}$--invariant coloring, which is
stronger than what is required by \fullref{def:colorable}.
Fix the power set of $S$ as one color palette.
If $W^{ess}=W_T$ give $\top$ color $T$. 
For $P\in\HHindex\setminus(\{\top\}\cup\mathfrak{P})$, color $P$ the unique subset $T\subset S$ such that
$P$ is conjugate to $W_T$.
If $P$ and $P'$ are both conjugate to $W_T$ then they are not properly
nested, since then they would have different ranks, and they do not
commute, because if they did then $P'$ would be conjugate into
$W_{T^\perp}$, but it cannot be since the unique special subgroup to
which it is conjugate is $W_T$ and $T\cap T^\perp=\emptyset$.
Thus, if $P$ and $P'\in \HHindex\setminus(\{\top\}\cup\mathfrak{P})$
have the same color then they are transverse.

Now consider the domains of $\HHindex^+\setminus\HHindex$.
They come from $\mathfrak{T}_P$ for some $P\in\mathfrak{P}$.
The $\mathfrak{T}_P$ are finite.
Take disjoint color palettes for the distinct $K_{(p)}$--orbits in
$\mathfrak{P}$, such that for a chosen representative $P$ of each
$K_{(p)}$ orbit of $\mathfrak{P}$, we can assign every domain in $\mathfrak{T}_P$ its own color.
Since the $\stab_{K_{(p)}}(P)$--action on $\mathfrak{T}_P$ factors through
the $\Lambda_P$--action, which is trivial on $\mathfrak{T}_P$, the
$\stab_{K_{(p)}}(P)$--action is color preserving.
Thus, there is a well-defined $K_{(p)}$--invariant coloring on the entire
union $\bigcup_{w\in K_{(p)}}w\mathfrak{T}_P$.
Moreover, distinct domains have the same color if and only if they
are $w.U$ and $w'.U$ for some $U\in\mathfrak{T}_P$ and $w,w'\in K_{(p)}$
with $w.P\neq w'.P$. By the argument of the previous paragraph,
$w.P\transverse w'.P$, so by definition of $\HHindex^+$, we also have
$w.U\transverse w'.U$.

The remaining domains are the dummy container domains, indexed by pairs
$(R,\mathfrak{O})$ for $R\in\HHindex^+$ and $\mathfrak{O}$ a
nonempty, nonmaximal set of pairwise orthogonal domains of
$\HHindex^+_R\setminus\{R\}$. Color these by a pair whose first entry
is the color of $R$ and whose second entry is the set of colors
occurring in $\mathfrak{O}$.
These are different from all of the colors assigned to
$\HHindex^{+}$, so it remains only to check that if
$(R,\mathfrak{O})$ and $(R',\mathfrak{O}')$ have the same color then
they are transverse.
Note that since elements of $\mathfrak{O}$ are pairwise orthogonal
they have distinct colors. 

Consider the orthogonality and nesting relations for dummy container
domains as constructed in \cite[Theorem~A.1]{AbbBehDur21}.
$(R,\mathfrak{O})$ and $(R',\mathfrak{O}')$ are orthogonal if and only
if $R\orthogonal R'$, but we already know that if $R$ and $R'$ have
the same color they are either equal or transverse, so if $(R,\mathfrak{O})$ and
$(R',\mathfrak{O}')$ have the same color they are not orthogonal.
Suppose $(R,\mathfrak{O})\nested(R',\mathfrak{O}')$ and they have the
same color. The nesting is true if and only if
$R\nested R'$ and for every $U'\in\mathfrak{O}'$ either $U'\orthogonal
R$ or $U'\nested U$ for some $U\in \mathfrak{O}$.
Since $(R,\mathfrak{O})$ and $(R',\mathfrak{O}')$ have the same color,
so do $R$ and $R'$, which they cannot if $R\nestedneq R'$, so
$R=R'$. Then for every $U'\in\mathfrak{O}'$ there exists
$U,V\in\mathfrak{O}$ with $U'\nested U$ and $V$ of the same color as
$U'$.
If $U\neq V$ then $U\orthogonal V$, by construction.
But then $U'\nested U$ and $U\orthogonal V$ implies $U'\orthogonal V$,
which contradicts that domains with the same color are transverse.
Thus, $U=V$, which further implies $U'=U$, since $U'$ cannot have the
same color as $U$ if it is properly nested.
Since this was true for every $U'\in\mathfrak{O}'$, we have $\mathfrak{O}=\mathfrak{O}'$.
Thus, distinct nested $(R,\mathfrak{O})$ and $(R',\mathfrak{O}')$
cannot have the same color.
Conclude that two dummy container domains with the same colors are
transverse.
Moreover, the coloring is $K_{(p)}$--invariant since its component colorings
were. 
\end{proof}

The second idea we pursue is to improve  the HHG
structure on $W$ relative to all irreducible higher rank affine
parabolic subgroups to an HHG structure relative to
the nonhyperoctahedral irreducible affine parabolic subgroups.
This is \fullref{thm:refined_rHHG}. \fullref{thm:diagram_HHG} shows
that this is the best we can do; the presence of a nonhyperoctahedral
irreducible affine parabolic subgroup obstructs the
existence of an HHG structure on $W$. 
\fullref{thm:refined_rHHG} follows by a similar
argument to \fullref{thm:Cox_vHHG}:

\begin{proof}[Proof of \fullref{thm:refined_rHHG}]
  Let $\mathfrak{P}_0$ and $\mathfrak{P}_1$ be the sets of irreducible
  higher rank affine parabolic subgroups of $(W,S)$ that are
  nonhyperoctahedral and hyperoctahedral, respectively.
  These sets are invariant under $W$--conjugation, so 
    $\mathfrak{P}=\mathfrak{P}_0\sqcup\mathfrak{P}_1$ is a
    $W$--invariant decomposition of the peripheral domains of the relative HHG structure $\HHindex$ of \fullref{thm:Cox_rHHG}.

    If $W$ is hyperbolic then it has no higher rank affine parabolic
    subgroups, so $\mathfrak{P}=\emptyset$ and $\HHindex$ is already
    an HHG structure for $W$, so we are done. 

    Suppose $P\in\mathfrak{P}_1$.
    Since $P$ is affine, it acts as a Euclidean reflection group on a
    Euclidean space $\mathbb{E}^n$, with $n=\rank(P)-1>1$.
    Since $P$ is hyperoctahedral, there is a choice of coordinates
    identifying $\mathbb{R}^n$ with $\mathbb{E}^n$ such that the
    $P$--action is by translation followed by a signed permutation of
    the coordinates of $\mathbb{R}^n$.
    Let $\mathfrak{T}^{prod}_P$ be the standard product HHS structure on
    $E_P:=(\mathbb{R}^n,\ell^1)$: its domains are indexed by nonempty
    subsets $I\subset \{1,\dots,n\}$ with nesting given by inclusion
    and orthogonality given by disjointness.
    The nesting minimal domains, indexed by singletons $\{i\}$, have
    coordinate space $\mathbb{R}$ with projection the $i$--th
    coordinate projection in $\mathbb{R}^n$.
    All other coordinate spaces are singletons.
    The HHS structure $\mathfrak{T}^{prod}_P$ makes $P$ an HHG. 

    If $W^{ess}$ is irreducible higher rank affine
    then we are done: if $W^{ess}$ is hyperoctahedral then its product
    HHS gives an HHG structure on $W$, and
    if not then we have a trivial relative HHG structure for $W$ whose
    unique domain is peripheral.

    We proceed, assuming $W$ is not hyperbolic and has essential
    part that is not irreducible affine, so $\HHindex\neq\{\top\}$.  

    Let $\mathfrak{Q}$ be the set of special subgroups in
    $\mathfrak{P}_1$, which accounts for exactly one representative of each
    $W$--orbit in $\mathfrak{P}_1$.
    For each $P\in\mathfrak{Q}$, let $\mathfrak{T}^{prod}_P$ be the product
    HHG structure on $E_P$ described above.
The orbit maps of $P$ on $E_P$ and $\coordinate{P}=C(P)$ induce a
coarsely $P$--equivariant quasiisometry $\phi_P\from \coordinate{P}\to
E_P$.
Let $\mathfrak{T}_P$ denote the HHS structure on
$\coordinate{P}$ pulled back from the structure
$\mathfrak{T}_P^{prod}$ on $E_P$ via $\phi_P$.
Since there are only finitely many conjugacy classes of parabolic
subgroup in $W$, the set $\mathfrak{Q}$ is finite, so the
$\mathfrak{T}_P$ have uniform HHS parameters.

$\stab_W(P)=\normalizer_W(P)=P\times P^\perp$, and the $P^\perp$--action on
$\coordinate{P}$ is trivial. Thus, $\stab_W(P)$ acts on
$\mathfrak{T}_P$ by factoring through the $P$--action.
Unlike in the proof of \fullref{thm:Cox_vHHG}, the $\stab_W(P)$--action
is not purely translational, but we used hyperoctahedrality to choose
$E_P$ such that the linear component of the $P$--action on $\mathfrak{T}_P^{prod}$ is by signed
permutations of the coordinate axes, so we still get the desired
permutations of the domains of $\mathfrak{T}_P^{prod}$ and corresponding
isometries between its coordinate spaces. 
Since $\phi_P$ is coarsely $P$--equivariant, the desired properties
also hold for $\mathfrak{T}_P$. 
The rest of the proof proceeds exactly  as in the proof of
\fullref{thm:Cox_vHHG}, using \fullref{prop:stabilizer_consistency} to
propagate the chosen $\mathfrak{T}_P$ for $P\in\mathfrak{Q}$ to all of
$\mathfrak{P}_1$ in a compatible way, and then using partial promotion
to get a new structure $\HHindex^{++}$ relative to $\mathfrak{P}_0$.
Cofiniteness of the $W$--action on $\HHindex^{++}$ follows
by the same argument as in \fullref{thm:Cox_vHHG}.
Thus, $\HHindex^{++}$ is an HHG structure on $W$ relative to $\mathfrak{P}_0$.
\end{proof}

Azuelos and Hagen \cite[Definition~7.1]{AzuHag} define a virtually Abelian subgroup
$A$ of $G$ to be \emph{highest} if $A$ is virtually $\mathbb{Z}^n$ and
no finite index subgroup of $A$ is contained in
a subgroup of $G$ isomorphic to $\mathbb{Z}^{n+1}$.

The proof of \fullref{thm:diagram_HHG} hinges on the following Coxeter
nonhyperoctahedral escalation lemma:
\begin{lemma}\label{lem:highest_nonhyperoctahedral}
   Every nonhyperoctahedral irreducible higher rank affine
   special subgroup $W_T$ of a finite rank Coxeter system $(W,S)$ is contained in a nonhyperoctahedral highest virtually
   Abelian subgroup of $W$. 
 \end{lemma}
\begin{proof}
  Let $T_0:=T$.
  If $T_0^\perp$ is nonspherical, choose a minimal nonspherical subset
  $T_1\subset T_0^\perp$.
  Continue inductively, choosing a minimal nonspherical subset
  $T_{j+1}$ of $\bigcap_{i=0}^j T_i^\perp$ until the intersection is
  spherical. 
  This process stops in finitely many steps since $S$ is finite and
  the nested intersections are strictly decreasing in cardinality.
Let $U:=T_0\sqcup\cdots\sqcup T_m$, so that the $W_{T_i}$ are pairwise
commuting and $U^\perp$ is spherical. 

 For each $i$, if  $W_{T_i}$ is irreducible affine, set $B_i:=W_{T_i}$.
  Otherwise, minimality requires that $W_{T_i}$ has no proper nonspherical
  parabolic subgroups, so if $W_{T_i}$ is not itself irreducible
  affine then it has no wide parabolic subgroups, so it is
  hyperbolic, by Moussong's Theorem.
Every nonspherical Coxeter group contains infinite
dihedral reflection subgroups. Choose $B_i$ to be such a subgroup of $W_{T_i}$.
  Define $A:=\prod_{i=0}^m B_i$.

The group $A$ is virtually Abelian: each factor $B_i$ has a finite
  index translation subgroup $\Lambda_i$ isomorphic to some
  $\mathbb{Z}^{n_i}$, where $n_i=1$ if $B_i$ is infinite dihedral and
  $n_i>1$ otherwise.
  Let $\Lambda:=\prod_{i=0}^m\Lambda_i$, so that
  $\Lambda$ is a finite index subgroup of $A$ that is isomorphic to
  $\mathbb{Z}^n$ for $n:=\sum_{i=0}^mn_i$.

We claim that $A$ is highest.
  Let $\Lambda'\leq \Lambda$ be a finite index subgroup, and define
  $\Lambda_i':=\Lambda'\cap \Lambda_i$.
  Each $\Lambda_i'$ has finite index in $\Lambda_i$, and
  $\prod_{i=0}^m\Lambda_i'$ has finite index in $\Lambda'$. 
For each $i$, $\pc(\Lambda_i')=W_{T_i}$, since $\Lambda_i'$ is
infinite and $W_{T_i}$ contains no proper nonspherical parabolic
subgroups.  
Thus, 
$\pc(\Lambda')=W_U$.
By construction $W_U$ is essential, so
$\centralizer_W(\Lambda')\leq\normalizer_W(W_U)=W_U\times W_U^\perp$.

The projection of $\centralizer_W(\Lambda')$ to each factor $W_{T_i}$
centralizes $\Lambda_i'$. 
If $W_{T_i}$ is affine, the centralizer of $\Lambda_i'$ in $B_i$ is
the full translation subgroup $\Lambda_i$.
Otherwise, $\Lambda_i'$ is an infinite cyclic subgroup of the
hyperbolic group $W_{T_i}$, so its centralizer in $W_{T_i}$ is
virtually cyclic.
Thus, in either case $\Lambda_i'$ has finite index in
$\centralizer_{W_{T_i}}(\Lambda_i')$. Since $W_U^\perp$ is spherical,
we conclude that $\Lambda'$ has finite index in $\centralizer_W(\Lambda')$.

This shows that $A$ is highest, because a finite index subgroup $A'$
of $A$ would give a finite index subgroup $\Lambda':=A'\cap \Lambda$ of
$\Lambda$.
Any free Abelian subgroup of $W$ containing $A'$ centralizes
$\Lambda'$, but we just showed that $\Lambda'$ has finite index in its
centralizer, so the centralizer has the same free Abelian rank as
$\Lambda'$ and $\Lambda$.

For each $i$, set $V_i:=\Lambda_i\otimes_{\mathbb{Z}}\mathbb{R}$, and
let $F_i\leq\GL(V_i)$ be the point group of $B_i$.
The point representation of $A$ acts on $V:=\bigoplus_{i=0}^m V_i$.
Its image is $F:=\prod_{i=0}^m F_i$, where each $F_i$ acts trivially
on $V_j$ for all $j\neq i$.
For each $i$, $F_i$ acts on $V_i$ with no nonzero fixed vectors: for
an irreducible affine factor this is its finite reflection
representation, and for an infinite dihedral factor it is the sign
action on $V_i\cong\mathbb{R}$.
Set $K:=\prod_{i=1}^mF_i$, with $K$ trivial if $m=0$.
Since $K$ acts trivially on $V_0$ and has no nonzero fixed vectors in
 $V_i$ for $i>0$, we have:
\begin{equation}
  \label{eq:8}
V^K=V_0\oplus\bigoplus_{i=1}^mV_i^{F_i}=V_0  
\end{equation}

Suppose that $A$ is hyperoctahedral.
Choose a basis $e_1,\cdots,e_n$ of $V$ such that $F$ acts by signed
permutation matrices on that basis.
Consider a $K$--orbit $\mathcal{O}$ of coordinate lines, and let
$E_{\mathcal{O}}$ be their span in $V$.
Choose $e_i$ such that $\mathbb{R}e_i\in\mathcal{O}$.
If $-e_i\in Ke_i$ then $K$ fixes
no nonzero vectors in $E_{\mathcal{O}}$.
In this case set $v_{\mathcal{O}}=0$. 
Otherwise, for every $\mathbb{R}e_j\in\mathcal{O}$ exactly one of $e_j$ or
$-e_j$ is in $Ke_i$, and the sum $v_{\mathcal{O}}$ of the included
vectors spans the $K$--fixed subspace of $E_{\mathcal{O}}$.

Elements of $F_0$ act by signed permutation on
$\{e_1,\dots,e_n\}$, so they send coordinate lines to coordinate
lines. Since $F_0$ commutes with $K$, it maps $K$--orbits of
coordinate lines to $K$--orbits of coordinate lines.
Thus, each $f\in F_0$ sends each $E_{\mathcal{O}}$ to some
$E_{\mathcal{O}'}$, and therefore sends $E_{\mathcal{O}}^K$ to some
$E_{\mathcal{O}'}^K$.
In particular, $F_0$ permutes the subspaces $\mathbb{R}v_{\mathcal{O}}$, respecting their dimensions.
Now, the nonzero $v_{\mathcal{O}}$ are a basis for the fixed
subspace of $K$ in $V$, but that is exactly $V_0$, by \eqref{eq:8}, so $F_0$ acts by
signed permutation on a basis of $V_0$.
This contradicts that $W_T=B_0$ is nonhyperoctahedral.
Therefore, $A$ is  nonhyperoctahedral.
\end{proof}
\begin{remark}\label{rem:standardvabelian}
  In the proof of \fullref{lem:highest_nonhyperoctahedral}, 
    the group $\Lambda$ is a standard free Abelian subgroup in the sense
    of Krammer \cite{Kra09}.
    The group $A$ is constructed in an analogous way, but choosing
    virtually Abelian factors instead of free Abelian. Azuelos and
    Hagen made a similar definition of ``standard virtually Abelian
    subgroup'' immediately before \cite[Corollary~7.11]{AzuHag},
    except that they chose the pieces contained in nonaffine factors to be
    $\mathbb{Z}$ instead of infinite dihedral. This difference turns
    out to be critical, even though both choices would yield highest
    virtually Abelian subgroups in our construction.
    We need the 1--dimensional $B_i$ to be infinite dihedral
    so that the corresponding $F_i$ has no nonzero fixed vectors on
    $V_i$. Replacing any of these infinite dihedral factors by its
    translation subgroup makes $V^K$ strictly larger than $V_0$, and this
    extra space can change hyperoctahedrality.
Consider, for example, $W$ that is the product of the $(3,3,3)$--triangle group
and the $(4,4,4)$--triangle group. The latter is hyperbolic with no
proper nonspherical parabolic subgroups. Our $A$ is of the form
$\widetilde{A}_2\times\widetilde{A}_1$, which is nonhyperoctahedral, by \fullref{ex:hyperoctahedral_times_Dinfty}.
Azuelos and Hagen choose $A'$ of the form $\widetilde{A}_2\times\mathbb{Z}$, which
\emph{is} hyperoctahedral, by \fullref{ex:hyperoctahedral_times_Z}.
Both $A$ and $A'$ are highest virtually Abelian subgroups of $W$, but
of the two, only $A$ witnesses that $W$ is not an HHG. 
\end{remark}

\begin{proof}[Proof of \fullref{thm:diagram_HHG}]
  The equivalence between the hyperoctahedral parabolic condition and
  the diagrammatic condition is established by
  \fullref{thm:hyperoct_classification}.
  
  By \fullref{thm:refined_rHHG}, if the Coxeter system
  contains no nonhyperoctahedral irreducible higher rank affine
  parabolic subgroups then $W$ is an
  HHG.
  Every HHG admits a coarsely equivariant coarse
  median: the standard HHS coarse median constructed in
  \cite[Theorem~7.3]{BehHagSis19} is coarsely equivariant
  by equivariance of the HHS projections and coarse
  uniqueness of realization.
 
Suppose $W$ admits a coarsely equivariant coarse
  median $\mu$.
  Suppose that $(W,S)$ has a nonhyperoctahedral
  irreducible higher rank affine parabolic subgroup.
  Up to conjugation, we may assume it is a special subgroup.
  By \fullref{lem:highest_nonhyperoctahedral}, it is
  contained in a nonhyperoctahedral highest virtually
  Abelian subgroup $A$ of $W$.
  Moreover, in  the proof of \fullref{lem:highest_nonhyperoctahedral},
  we constructed such an $A$ with a finite index free Abelian subgroup
  $\Lambda\leq A$ that also has finite
  index in $\centralizer_W(\Lambda)$.
  Fioravanti, Levcovitz, and Sageev \cite[Definition~2.13 and
  Lemma~2.15]{FioLevSag24} show $\centralizer_W(\Lambda)$ is
  ``quasi-submedian'', in the sense that
$\mu(\centralizer_W(\Lambda), \centralizer_W(\Lambda),
\centralizer_W(\Lambda))$ is at finite Hausdorff distance from
$\centralizer_W(\Lambda)$. 
  Since $A$ and $\centralizer_W(\Lambda)$
  contain $\Lambda$ with finite index, they are at finite
  Hausdorff distance from each other in $W$.
  Thus $A$ is also quasi-submedian.
  Define $\mu'\from A^3\to A$ by first restricting $\mu$ to $A^3$ and
  then choosing a closest point in
  $A$. Since $A$ is quasi-submedian, $\mu|_{A^3}$ and $\mu'$ are 
  at uniformly bounded distance from one another as maps $A^3\to W$.
  Since $A$ is a virtually Abelian subgroup of the CAT(0) group $W$,
  it is undistorted, so a word metric on $A$ is comparable to the
  restriction to $A$ of a word metric on $W$.
  It follows that $\mu'$ is a coarsely equivariant coarse median
  on $A$; the metric compatibility conditions of
  \fullref{def:NWZ} and \fullref{def:cecm} are satisfied by $\mu$ with
  respect to the $W$--metric, so by bounded distance between $\mu$ and
  $\mu'$ they are satisfied by $\mu'$ with respect to the $W$--metric,
  thus by undistortion they are satisfied by $\mu'$ with respect to the $A$--metric. 
The existence of a coarsely equivariant coarse median on $A$ contradicts  \fullref{thm:vAbelian_coarse_median}, since $A$ was
constructed to be nonhyperoctahedral.
Thus, if $W$ admits a coarsely equivariant coarse median then $(W,S)$ does not contain a nonhyperoctahedral irreducible higher
rank affine parabolic subgroup. 
\end{proof}

\section*{AI Use}
The authors conceived the arguments, developed the proofs, and wrote
the entirety of this article.
Improvements were made after considering comments from Anthropic's
Claude Opus 5.5 on earlier drafts.

\appendix
\section{Promotion of relative HHS structures}\label{sec:promotion}
In the main text we produce a relative HHG structure on Coxeter
groups, observe that the peripheral domains have coordinate spaces
that are themselves HHSs, and then need a promotion theorem that
combines these structures to make a full HHS/HHG structure.
Special cases of such promotion theorems appear in the literature for
relatively hyperbolic groups \cite[Theorem~9.1]{BehHagSis19} and graph
products \cite[Theorem~5.1]{BerRus22}.
Rather than make a Coxeter-specific special case, in this section we
distill a reusable, general HHS statement.
The argument follows that of Berlyne and Russell \cite{BerRus22}, who
already note that the bulk of their proof only uses general HHS
properties, not the specifics of their graph product setting.

\begin{definition}[{\cite[Definition~3.2]{AbbBehDur21}}]
  An HHS has \emph{bounded domain dichotomy} if there exists $B$ such
  that every bounded coordinate space has diameter at most $B$. 
\end{definition}
\begin{definition}[{\cite[Definition~3.3]{AbbBehRus25}}]
  Let $G$ be a group. A \emph{$G$--HHS} is an HHS structure on a
  Cayley graph of $G$ that has bounded domain dichotomy and satisfies all of the requirements of an
  HHG, \fullref{def:relHHG}, except that the $G$ action on the domains
  with bounded coordinate space is not required to be cofinite.
\end{definition}

\begin{theorem}[Partial Berlyne-Russell Promotion]\label{thm:relHHS_promotion}
  Let $(X,\HHindex)$ be an HHS structure relative to a set of
  peripheral domains $\mathfrak{P}\subset\HHindex$.
  Suppose that $\mathfrak{P}=\mathfrak{P}_0\sqcup\mathfrak{P}_1$, and for each $P\in\mathfrak{P}_1$ the coordinate space
  $\coordinate{P}$ admits an HHS structure
  $(\coordinate{P},\mathfrak{T}_P)$, and that
  the HHS parameters are
  uniformly bounded for all $(\coordinate{P},\mathfrak{T}_P)$.
  Then $X$ admits an HHS structure relative to $\mathfrak{P}_0$.

  Suppose, moreover, that $(G,\HHindex)$ is an HHG structure relative
  to $\mathfrak{P}$, that the decomposition
  $\mathfrak{P}=\mathfrak{P}_0\sqcup\mathfrak{P}_1$ is $G$--invariant,
  and that the HHS structures on the peripheral
  coordinate spaces of $\mathfrak{P}_1$ satisfy the following $G$--compatibility
  conditions, where for all $g\in G$ and $P\in\mathfrak{P}_1$ the map $g_P\from
    \coordinate{P}\to\coordinate{g.P}$ is the isometry from the
    definition of relative HHG, \fullref{def:relHHG}.
  \begin{enumerate}
  \item There is an $\nested$ and $\orthogonal$--preserving, cofinite action of $G$ on
    $\coprod_{P\in\mathfrak{P}_1}\mathfrak{T}_P$ with $g$ sending
    $\mathfrak{T}_P$ bijectively to $\mathfrak{T}_{g.P}$ by $U\mapsto
    g.U$.
    \label{item:promotion:bijection}
    \item For all $g,h\in G$, $P\in\mathfrak{P}_1$, 
      $U,V\in\mathfrak{T}_P$, and $x\in\coordinate{P}$:
      \begin{itemize}
        \item There is an isometry $g_U\from\coordinate{U}\to\coordinate{g.U}$.
      \item $(gh)_U=g_{h.U}\circ h_U$
      \item $\domaindist{g.U}{g_U(\proj{U}(x))}{\proj{g.U}(g_P(x))}$ is
        uniformly bounded.
      \item $(U\nestedneq V)\lor(U\transverse
        V)\implies\domaindist{g.V}{g_V(\relprojection{U}{V})}{\relprojection{g.U}{g.V}}$
        is uniformly bounded. 
      \end{itemize}\label{item:promotion:compatibility}
      \end{enumerate}
      Then the promoted structure making $G$ an HHS relative to
      $\mathfrak{P}_0$ satisfies all of the requirements of being an
      HHG relative to $\mathfrak{P}_0$ except possibly that the
      $G$--action on domains with bounded coordinate spaces may not be
      cofinite. In particular, if $\mathfrak{P}_0=\emptyset$ then
      the promoted relative HHS structure on $G$ is a $G$--HHS.
\end{theorem}
\begin{proof}
Define $\HHindex^+:=(\HHindex\setminus\mathfrak{P}_1)\cup
 \coprod_{P\in\mathfrak{P}_1}\mathfrak{T}_P$.
 The new nesting relation is unchanged within
 $(\HHindex\setminus\mathfrak{P}_1)$ and within each $\mathfrak{T}_P$.
 When $Q\in (\HHindex\setminus\mathfrak{P}_1)$ and $R\in
 \mathfrak{T}_P$, define $R\nested Q$ if and only if $P\nested Q$.
Similarly, the new orthogonality relation is unchanged within
$(\HHindex\setminus\mathfrak{P}_1)$ and within each $\mathfrak{T}_P$.
If $Q\in (\HHindex\setminus\mathfrak{P}_1)$ and $R\in
 \mathfrak{T}_P$, define $Q\orthogonal R$ if and only if
 $Q\orthogonal P$. If $R_0\in\mathfrak{T}_{P_0}$ and
 $R_1\in\mathfrak{T}_{P_1}$ with $P_0\neq P_1$, then define
 $R_0\orthogonal R_1$ if and only if $P_0\orthogonal P_1$.
 Note that domains of $\mathfrak{P}_0$ remain nesting minimal in $\HHindex^+$.
 Berlyne and Russell then define projections and relative projections
 and check, \cite[Claims 5.2--5.4]{BerRus22}, that
 all of the relative HHS axioms are satisfied except possibly for the
 Containers Axiom. They call $\HHindex^+$ a ``proto-hierarchy structure''.
 They also say, \cite[Lemma~5.5]{BerRus22}, citing
 \cite[Lemma~1.5]{BehHagSis17}, that from the uniform HHS constants one can
 derive a bound on the maximal cardinality of a set of pairwise
 orthogonal domains in $\HHindex^+$.
 With this condition, their \cite[Theorem~A.1]{AbbBehDur21} adds ``dummy''
 container domains whose coordinate spaces are single points to make a
 larger structure $\HHindex^{++}$.
 By their construction, nesting minimal domains of $\HHindex^+$ remain
 nesting minimal in $\HHindex^{++}$, so the domains of $\mathfrak{P}_0$ are nesting minimal in $\HHindex^{++}$.
 Thus, $\HHindex^{++}$ is an HHS structure relative
 to $\mathfrak{P}_0$.
 Moreover, $\HHindex^{++}$ keeps the original nesting and
 orthogonality relations between domains from $\HHindex^+$.

 The additional $G$--compatibility conditions are to guarantee that the
 coarse equivariance conditions of \fullref{def:relHHG}, which are
 satisfied by $\HHindex$, by hypothesis, also are satisfied by
 $\HHindex^{++}$.
 \cite[Claim~5.6]{BerRus22} shows that the compatibility conditions
 imply that $\HHindex^+$ satisfies the definition of relative HHG,
 except for the Containers Axiom. 
 \cite[Theorem~5.7]{BerRus22} says that after the addition of the
 dummy containers, $\HHindex^{++}$ is an HHG structure on $G$,
 relative to $\mathfrak{P}_0$.
 However, as noted in \cite[Remark~3.4]{AbbBehRus25}, that theorem has a
 gap: it does not establish cofiniteness
 of the group action on the dummy container domains. 
 The definition of ``$G$--HHS'' circumvents this problem
 by not requiring cofiniteness on those domains.
 Bounded domain dichotomy is clear
 because the $G$--action on $\HHindex^+$ is cofinite and the domains
 of $\HHindex^{++}\setminus\HHindex^+$ have single point coordinate
 spaces. 
\end{proof}
\begin{corollary}\label{cor:dummy_cofiniteness}
  The structure $\HHindex^{++}$ produced by the second branch
  of \fullref{thm:relHHS_promotion} by adding dummy container domains
  to the proto-hierarchy structure $\HHindex^+$ is a relative HHG structure if
  and only if $G$ acts cofinitely on the set of pairs
  $(R,\mathfrak{O})$, where $R\in\HHindex^+$ and $\mathfrak{O}$ is a
  nonempty, nonmaximal collection of pairwise orthogonal domains in $\HHindex^+_R\setminus\{R\}$.
\end{corollary}
\begin{proof}
  Pairs of the form $(R,\mathfrak{O})$ index the dummy container
  domains added to $\HHindex^+$ in \cite[Theorem~A.1]{AbbBehDur21}. 
\end{proof}

The next result gives a sufficient
condition that lets us generate $G$--compatible HHS structures on the
peripheral coordinate spaces.
Given a $G$--invariant decomposition
$\mathfrak{P}=\mathfrak{P}_0\sqcup\mathfrak{P}_1$, 
the idea is to pick a representative $P$ of each $G$--orbit in
$\mathfrak{P}_1$, take an HHS structure $\mathfrak{T}_P$ on
$\coordinate{P}$, and use the $G$--action to transport that HHS
structure to the coordinate spaces of the other peripherals in the $G$--orbit of $P$.
The desired $G$--compatibility will follow, by construction, but for
this $G$--transport to be well-defined we need conditions
that say that $\stab_G(P)$ respects the $\mathfrak{T}_P$ structure.

\begin{proposition}\label{prop:stabilizer_consistency}
  Suppose $(G,\HHindex)$ is a relative HHG structure with peripherals
  $\mathfrak{P}$, and suppose there is 
a $G$--invariant decomposition
$\mathfrak{P}=\mathfrak{P}_0\sqcup\mathfrak{P}_1$.
  Let $\mathfrak{Q}\subset\mathfrak{P}_1$ be a set of
  representatives for the $G$--orbits in $\mathfrak{P}_1$.
  Let $(\coordinate{P},\mathfrak{T}_P)$ be HHS structures for each
  $P\in\mathfrak{Q}$.
  Suppose, for all $P\in\mathfrak{Q}$, $g,h\in\stab_G(P)$, and $U,V\in\mathfrak{T}_P$:
  \begin{enumerate}
  \item There is a cofinite action of $\stab_G(P)$ on $\mathfrak{T}_P$ by
    bijections $\mathfrak{T}_P\to\mathfrak{T}_P:U\mapsto
    g.U$ that preserve
    the $\nested$ and $\orthogonal$ relations.\label{pp:action}
    \item There are isometries $g_U\from \coordinate{U}\to\coordinate{g.U}$.\label{pp:isometries}
    \item $(gh)_U=g_{h.U}\circ h_U$\label{pp:composition}
    \item $g_U\circ \proj{U}\ceq \proj{g.U}\circ g_P$\label{pp:projection}
      \item $(U\nestedneq V)\lor(U\transverse V)\implies g_V(\relprojection{U}{V})\ceq \relprojection{g.U}{g.V}$\label{pp:relprojection}
      \end{enumerate}
Suppose, moreover, that the coarse equivalences of the last two items
are uniform.
  Then HHS structures for the coordinate spaces of all $P\in\mathfrak{P}_1$ can be defined by
  transporting the chosen structures from $\mathfrak{Q}$ by the $G$--action, and this
  collection of HHS structures on the peripheral coordinate spaces
  satisfies the hypotheses of \fullref{thm:relHHS_promotion}, so
  $\HHindex^{++}$ satisfies all of the requirements of being an HHG
  relative to $\mathfrak{P}_0$, except possibly for cofiniteness of
  the $G$--action on domains with bounded coordinate spaces. In
  particular, if $\mathfrak{P}_0=\emptyset$ then $\HHindex^{++}$ is a $G$--HHS.
\end{proposition}
\begin{proof}
  Suppose $P\in\mathfrak{Q}$.
Define the transported HHS structures on the $G$--orbit of $P$ by:
\[G\times_{\stab_G(P)}\mathfrak{T}_P:=(G\times\mathfrak{T}_P)/(g,U)\sim
  (gh,h^{-1}\!.U)\text{ for }h\in\stab_G(P)\]
This collection contains the original $\mathfrak{T}_P$ via the
identification $U=[1,U]$.
The $G$--action on this collection of domains is $a.[g,U]=[ag,U]$.
Define the set $\mathfrak{T}_{g.P}$ of domains of $\coordinate{(g.P)}$
to be  $\{[g,U]\mid U\in
\mathfrak{T}_P\}$.
This only depends on $g.P$, since if $g.P=k.P$ then
$h^{-1}:=k^{-1}g\in\stab_G(P)$,
so: \[\mathfrak{T}_{k.P}=\mathfrak{T}_{(gh).P}=\{[gh,U]\mid
  U\in\mathfrak{T}_P\}=\{[g,U']\mid U'\in h.(\mathfrak{T}_P)\}=\{[g,U']\mid U'\in \mathfrak{T}_P\}=\mathfrak{T}_{g.P}\]
For $[g,U],[g,V]\in\mathfrak{T}_{g.P}$ define the nesting and
orthogonality relations to be the same as the corresponding relations
for $U,V\in\mathfrak{T}_P$.
These relations are well-defined and independent of choice of
representatives because the action of $\stab_G(P)$ on $\mathfrak{T}_P$
preserves nesting and orthogonality.
Thus,  Condition~\eqref{item:promotion:bijection} of
\fullref{thm:relHHS_promotion} is satisfied.
Cofiniteness is also built in, since $\mathfrak{Q}$ is finite and
the action of $\stab_G(P)$ on $\mathfrak{T}_P$ was assumed to be
cofinite for each $P\in\mathfrak{Q}$.

For $[g,U]\in\mathfrak{T}_{g.P}$, pick a preferred representative
$(g,U)\in [g,U]$, with the convention that for a domain of the form
$[1,U]$ we choose $(1,U)$ to be the
preferred representative.
For a domain of the form $[1,U]$ define its coordinate space to be
$\coordinate{[1,U]}:=\coordinate{U}$.
For $h\in\stab_G(P)$ define $h_{[1,U]}:=h_U$.
For a domain $[g,U]$ with preferred representative $(g,U)$ with $g\neq
1$, which implies $g\notin\stab_G(P)$, define $\coordinate{[g,U]}$ to
be a copy of $\coordinate{[1,U]}$ with the identification
$g_{[1,U]}\from \coordinate{[1,U]}\to\coordinate{[g,U]}$.

For any $k\in G$ and domain $[g,U]$ with preferred representative
$(g,U)$, we have $k.[g,U]=[kg,U]$, which has some preferred
representative $(kgh,h^{-1}\!.U)$, for a unique $h\in\stab_G(P)$.
Define the isometry $k_{[g,U]}$ as a composition of the isometries of
the previous paragraph (in particular, if $g=1$ this definition
agrees with that of the previous paragraph):
\[k_{[g,U]}\from\coordinate{[g,U]}\to\coordinate{(k.[g,U])}=\coordinate{[kgh,h^{-1}\!.U]}:=(kgh)_{[1,h^{-1}\!.U]}\circ
  (h^{-1})_{[1,U]}\circ (g_{[1,U]})^{-1}\]
This definition satisfies the first two bullets of
 Condition~\eqref{item:promotion:compatibility} of \fullref{thm:relHHS_promotion}.

 Define projection maps by:
  \[\proj{[g,U]}\from \coordinate{(g.P)}\to\coordinate{[g,U]} :
    x\mapsto\bigcup_{h\in\stab_G(P)}
    (gh)_{[1,h^{-1}\!.U]}\circ \proj{h^{-1}\!.U}\circ ((gh)_P)^{-1}(x)\]

To see that this projection is coarsely well-defined, we show, using the coarse equivalence of hypothesis \eqref{pp:projection}, that the
uniformly bounded sets of the union are coarsely equivalent, so their
union is uniformly bounded:
\begin{align*}
(gh)_{[1,h^{-1}\!.U]}\circ\proj{h^{-1}\!.U}\circ
                       ((gh)_P)^{-1}    &=g_{[1,U]}\circ h_{[1,h^{-1}\!.U]}\circ\proj{h^{-1}\!.U}\circ (h_P)^{-1}\circ
                                        (g_{P})^{-1}  \\
  &=g_{[1,U]}\circ h_{h^{-1}\!.U}\circ\proj{h^{-1}\!.U}\circ (h_P)^{-1}\circ
    (g_{P})^{-1}  \\
  &\ceq g_{[1,U]}\circ \proj{U}\circ h_{P} \circ (h_P)^{-1}\circ
    (g_{P})^{-1}  \\
  &=g_{[1,U]}\circ \proj{U}\circ
    (g_{P})^{-1}  
\end{align*}

We further claim these definitions imply the third bullet of
\fullref{thm:relHHS_promotion}\eqref{item:promotion:compatibility}.
Specifically, for 
$[g,U]\in\mathfrak{T}_{g.P}$, $x\in\coordinate{(g.P)}$, and $k\in G$,
we claim:
\[\domaindist{k.[g,U]}{k_{[g,U]}(\proj{[g,U]}(x))}{\proj{k.[g,U]}(k_{g.P}(x))}=0\]
This is because the two sets intersect, as follows.
If $(g,U)\in
[g,U]$ and $(kgh_0,h_0^{-1}\!.U)\in [kgh_0,h_0^{-1}\!.U]=k.[g,U]$, where
$h_0\in\stab_G(P)$,  are the preferred representatives, we observe
from the definitions of the isometries:
\[(gh_0)_{[1,h_0^{-1}\!.U]}=g_{[1,U]}\circ((h_0^{-1})_{[1,U]})^{-1}\]
Using this,  we produce a nontrivial intersection by specializing each union
to the $h_0$ term:
\begin{align*}
  k_{[g,U]}(\proj{[g,U]}(x))&=(kgh_0)_{[1,h_0^{-1}\!.U]}\circ(h_0^{-1})_{[1,U]}\circ
                             (g_{[1,U]})^{-1}\left(\bigcup_{h\in\stab_G(P)}(gh)_{[1,h^{-1}\!.U]}\circ\proj{h^{-1}\!.U}\circ((gh)_P)^{-1}(x)\right)\\
  &\supset (kgh_0)_{[1,h_0^{-1}\!.U]}\circ(h_0^{-1})_{[1,U]}\circ
                             (g_{[1,U]})^{-1}\circ
    (gh_0)_{[1,h_0^{-1}\!.U]}\circ\proj{h_0^{-1}\!.U}\circ((gh_0)_P)^{-1}(x)\\
  &=(kgh_0)_{[1,h_0^{-1}\!.U]}\circ\proj{h_0^{-1}\!.U}\circ((gh_0)_P)^{-1}(x)\\
  &=(kgh_0)_{[1,h_0^{-1}\!.U]}\circ\proj{h_0^{-1}\!.U}\circ
    ((kgh_0)_{P})^{-1}\circ k_{g.P}(x)\\
                           &\subset\bigcup_{h\in\stab_G(P)}(kgh)_{[1,h^{-1}\!.U]}\circ\proj{h^{-1}\!.U}\circ
    ((kgh)_{P})^{-1}(k_{g.P}(x))\\
  &=\proj{k.[g,U]}(k_{g.P}(x))
\end{align*}

Similarly, if $U,V\in\mathfrak{T}_P$ with $U\nestedneq V$ or
$U\transverse V$, define the relative projection by:
\[\relprojection{[g,U]}{[g,V]}:=\bigcup_{h\in\stab_G(P)}(gh)_{[1,h^{-1}\!.V]}(\relprojection{h^{-1}\!.U}{h^{-1}\!.V})\]
Again, this is a union of uniformly bounded diameter, uniformly
coarsely equivalent sets, so it has uniformly bounded diameter: for $h\in \stab_G(P)$, the pairs $(gh,h^{-1}\!.U)$ and $(gh,h^{-1}\!.V)$ are different representatives
of $[g,U]$ and $[g,V]$,
and, using hypothesis \eqref{pp:relprojection}:
\begin{align*}
  (gh)_{[1,h^{-1}\!.V]}(\relprojection{h^{-1}\!.U}{h^{-1}\!.V})&=g_{[1,V]}\circ h_{[1,h^{-1}\!.V]} (\relprojection{h^{-1}\!.U}{h^{-1}\!.V})\\
  &=g_{[1,V]}\circ h_{h^{-1}\!.V}  (\relprojection{h^{-1}\!.U}{h^{-1}\!.V})\\
  &\ceq g_{[1,V]}(\relprojection{U}{V})
\end{align*}

Finally, the fourth bullet of
\fullref{thm:relHHS_promotion}\eqref{item:promotion:compatibility} follows, because again we claim
that for any $P\in\mathfrak{Q}$,  $g,k\in G$ and $[g,U],[g,V]\in\mathfrak{T}_{g.P}$ with
$[g,U]\nestedneq[g,V]$ or $[g,U]\transverse[g,V]$ we have:
\[\domaindist{k.[g,V]}{k_{[g,V]}(\relprojection{[g,U]}{[g,V]})}{\relprojection{k.[g,U]}{k.[g,V]}}=0\]
Again, this is because the two sets intersect.
Let $(g,V)$ and $(kgh_0,h_0^{-1}\!.V)$ be the preferred
representatives of $[g,V]$ and $k.[g,V]=[kgh_0,h_0^{-1}\!.V]$.
Specialize the unions defining $\relprojection{k.[g,U]}{k.[g,V]}=\relprojection{[kg,U]}{[kg,V]}$
and $\relprojection{[g,U]}{[g,V]}$ to their $h_0$ terms. Then the identity
$(gh_0)_{[1,h_0^{-1}\!.V]}=g_{[1,V]}\circ((h_0^{-1})_{[1,V]})^{-1}$ and the
definition of $k_{[g,V]}$, as in the previous calculation, show that both
$\relprojection{k.[g,U]}{k.[g,V]}$ and $k_{[g,V]}(\relprojection{[g,U]}{[g,V]})$ contain
$(kgh_0)_{[1,h_0^{-1}\!.V]}(\relprojection{h_0^{-1}\!.U}{h_0^{-1}\!.V})$.
\end{proof}

\section{Coarsely equivariant coarse medians on virtually Abelian groups}\label{sec:equivariantcoarsemedians}
Recall that the proof of \fullref{thm:diagram_HHG} uses
\fullref{lem:highest_nonhyperoctahedral} to escalate a
nonhyperoctahedral irreducible higher rank affine parabolic subgroup to
a nonhyperoctahedral highest virtually Abelian subgroup. 
In this section we show that for finitely generated virtually Abelian
groups hyperoctahedrality is equivalent to the existence of a coarsely
equivariant coarse median. 

\begin{definition}\label{def:NWZ}
  A $C$--\emph{coarse median space} $(X,d,\mu)$ is
  a metric space $(X,d)$ equipped with a ternary operation
  $\mu\from X^3\to X$ such that, writing $x\sim_C y$ if
  $d(x,y)\leq C$, the following hold for all $v,w,x,y,z\in X$:
  \begin{enumerate}[label={(C\arabic*)}]
  \item $\mu(x,x,y)=x$\label{NWZ1}
  \item $\mu(x,y,z)=\mu(y,x,z)=\mu(x,z,y)$\label{NWZ2}
  \item $d(\mu(x,y,z),\mu(x,y,w))\leq Cd(z,w)+C$\label{NWZ3}
  \item $\mu(\mu(x,w,y),w,z)\sim_C
    \mu(x,w,\mu(y,w,z))$\label{NWZ4}
  \item $\mu(v,w,\mu(x,y,z))\sim_C
    \mu(\mu(v,w,x),\mu(v,w,y),z)$\label{NWZ5}
  \end{enumerate}
  It is a \emph{coarse median space} if it is a $C$--coarse median
  space for some $C$.
\end{definition}

This is the normalized formulation of Niblo, Wright, and Zhang
\cite[Theorem~3.1]{NibWriZha19}, which is equivalent to Bowditch's formulation
\cite{Bow13b} up to a uniformly bounded change of the operation.
Such a change allows \ref{NWZ1} and \ref{NWZ2} to hold exactly
\cite[Remark~2.10]{NibWriZha19}.
By \cite[Proposition~3.11]{NibWriZha19}, given \ref{NWZ1}, \ref{NWZ2},
and \ref{NWZ3}, each of \ref{NWZ4} and \ref{NWZ5} implies the other,
after possibly increasing $C$, so it suffices to take either one
in the definition. 

A coarse median can be pulled back by a quasiisometry, with a
uniformly bounded adjustment to enforce \ref{NWZ1} and \ref{NWZ2}
if necessary. Thus, the existence of a coarse median is a
quasiisometry invariant. In particular, it makes sense to say that a finitely generated group
is a \emph{coarse median group} if one, equivalently every, Cayley
graph with respect to a finite generating set is a coarse median
space.

\begin{definition}\label{def:cecm}
  Let $G$ be a finitely generated group equipped with a word metric
  $d$ and a $C$--coarse median $\mu$. We call $\mu$ \emph{coarsely equivariant}
  if for all $g,x,y,z\in G$:
  \begin{enumerate}[label={(C6)}]
  \item $g\mu(x,y,z)\sim_C\mu(gx,gy,gz)$\label{NWZ6}
  \end{enumerate}
\end{definition}

\begin{definition}
  A \emph{cubulation} of a group $G$ is a proper action
  $G\curvearrowright Y$ by cubical automorphisms on a CAT(0) cube
  complex. 
\end{definition}

\begin{definition}[{\cite[Section~2.4]{FioLevSag24}}]  
  A coarse median $\mu$ on a finitely generated group $G$ is
  \emph{cubical} if there exist a cocompact cubulation
  $G\curvearrowright Y$, a vertex $y\in Y^{(0)}$, and a constant
  $D\geq 0$ such that, for all $f,g,h\in G$, if $m_Y$ is the median on $Y^{(0)}$:
  \[d_Y\big(\mu(f,g,h)y,m_Y(fy,gy,hy)\big)\leq D\]
  Equivalently, $\mu$ is at uniformly bounded distance from the
pullback of $m_Y$ by the orbit quasiisometry $g\mapsto gy$.
\end{definition}

\begin{theorem}\label{thm:vAbelian_coarse_median}
  Let $G$ be virtually $\mathbb{Z}^n$, for
  some $n\geq 1$. The following are equivalent:
  \begin{enumerate}
  \item $G$ is hyperoctahedral.\label{item:hyperoctahedral}
    \item $G$ is cocompactly cubulable. \label{item:cubulable}
\item $G$ admits a cubical coarse median.\label{item:cubicalcoarsemedian}
\item $G$ admits a coarsely equivariant coarse median. \label{item:equivariantcoarsemedian}
\end{enumerate}
\end{theorem}
The main work for us is the implication
\eqref{item:equivariantcoarsemedian}$\implies$\eqref{item:hyperoctahedral}.
The idea is to use two limiting arguments: first an argument of
Bowditch turns a coarsely equivariant coarse
median on $G$ into an equivariant topological median on an asymptotic
cone of $G$, which is homeomorphic to $\mathbb{R}^n$ (see
\cite[Chapter~10]{DruKap18} for background on asymptotic cones). Then
we use arguments modeled on those of Haettel, \cite[Propositions~G and
4.1]{Hae16} to say that
there is a biLipschitz median embedding of a standard $n$--cube into
the asymptotic cone.
A differentiation and rescaling argument produces a median on the cone that is
  $F$--equivariant and linearly conjugate to the standard $n$--cube
  median.
  The second part is broken off separately as an auxiliary lemma:
  \begin{lemma}\label{lem:linear_median_obstruction}
  For $n\geq 1$, let $F\leq \GL_n(\mathbb{R})$ be finite.
  Suppose $\mu$ is a Lipschitz median on $\mathbb{R}^n$ that is
  equivariant with respect to the actions of $F$ and the translations of $\mathbb{R}^n$. 
  Then $F$ is conjugate into $\orthog_n(\mathbb{Z})$.
\end{lemma}
\begin{proof}
  Since $\mu$ is Lipschitz it is, in particular, continuous, so
  \cite[Corollary~3.10]{BesBroSag} implies that it is locally convex
  with rank at most $n$.
  By \cite[Corollary~3.7]{Hae16}, there does not exist a Lipschitz,
  locally convex median of rank strictly less than $n$ on
  $\mathbb{R}^n$, so the rank of $\mu$ is exactly $n$.
  Thus, \cite[Proposition~G]{Hae16} says there exists a biLipschitz
  median embedding $\phi\from([0,1]^n,\nu)\to (\mathbb{R}^n,\mu)$,
  where $\nu$ is the standard coordinate median on $\mathbb{R}^n$.
The goal now is to transfer $F$--equivariance from $\mu$ to $\nu$.

  By Rademacher's theorem, $\phi$ is differentiable almost everywhere.
  Choose an interior point $p$ of the $n$--cube at which $\phi$ is
  differentiable, and let $D_p \phi$ be the derivative of $\phi$ at $p$,
  which is invertible, since $\phi$ is biLipschitz.  
For all $u\in \mathbb{R}^n$ and all sufficiently small $\epsilon>0$
the point $p+u\epsilon$ is in the cube, so we can define
$\phi_\epsilon(u):=\left(\phi (p+\epsilon u)-\phi
  (p)\right)/\epsilon$, and we have
$\phi_\epsilon(u)\xrightarrow{\epsilon\to 0} (D_p \phi)(u)$.
  For $\epsilon>0$, define $\mu_\epsilon(x,y,z):=\mu(\epsilon
  x,\epsilon y,\epsilon z)/\epsilon$.

  Translation equivariance of $\mu$ and median preservation of $\phi$
  give, for all $u,v,w\in\mathbb{R}^n$ and small enough $\epsilon$:
  \begin{align*}
\epsilon    \mu_\epsilon(\phi_\epsilon(u), \phi_\epsilon(v), \phi_\epsilon(w))
    &=\mu\left(\phi (p+\epsilon u)-\phi (p), \phi (p+\epsilon
v)-\phi (p), \phi (p+\epsilon w)-\phi (p)\right)\\
    &=\mu(\phi (p+\epsilon
      u), \phi (p+\epsilon v), \phi (p+\epsilon w))-\phi (p)\\
     &=\phi (\nu(p+\epsilon
       u, p+\epsilon v, p+\epsilon w))-\phi (p)\\
     &=\phi (p+\epsilon\nu(u,v,w))-\phi (p)\\
    &=\epsilon\phi_\epsilon(\nu(u,v,w))
  \end{align*}
 Since the $\mu_\epsilon$ share the same Lipschitz constant as $\mu$,
  the family is equicontinuous, so:
  \begin{align*}
    (D_p\phi)(\nu(u,v,w))&=\lim_{\epsilon\to
                            0}\phi_\epsilon(\nu(u,v,w))\\
    &=\lim_{\epsilon\to
      0}\mu_\epsilon(\phi_\epsilon(u),\phi_\epsilon(v),\phi_\epsilon(w))\\
    &=\lim_{\epsilon\to
      0}\mu_\epsilon\left((D_p\phi)(u),(D_p\phi)(v),(D_p\phi)(w)\right)
  \end{align*}
 Thus, we can define a limit median:
\[    \mu_0(a,b,c) :=\lim_{\epsilon\to 0}\mu_\epsilon\left(a,b,c\right)=(D_p\phi)\left(\nu\left((D_p\phi)^{-1}(a), (D_p\phi)^{-1}(b), (D_p\phi)^{-1}(c)\right)\right)\]

  Since $\mu$ is equivariant with respect to each linear $f\in F$, so
  is each $\mu_\epsilon$, and hence so is $\mu_0$.
  Thus, $\nu$ is $(D_p\phi)^{-1}F (D_p\phi)$--equivariant.

  Since each $f'\in (D_p\phi)^{-1}F (D_p\phi)$ is linear it fixes
  0.
  By equivariance, it takes the median interval $[0,e_i]_{\nu}$ to the
  median interval $[0,f'e_i]_{\nu}$.
  Since $f'$ is a linear automorphism and $[0,e_i]_{\nu}$ is
  a nondegenerate line segment, so is $[0,f'e_i]_{\nu}$, which
implies $f'e_i$ has exactly 1 nonzero
  coordinate.
  Since $f'$ is invertible, this implies that $f'$ permutes the
  coordinate lines. 
 By averaging an inner product on $\mathbb{R}^n$ over the action of
 the finite group $(D_p\phi)^{-1}F (D_p\phi)$ and normalizing the
 coordinate vectors to have unit length with respect to this new inner
 product, conclude that $(D_p\phi)^{-1}F (D_p\phi)$ permutes these
 coordinate vectors up to sign. 
\end{proof}

\begin{proof}[Proof of \fullref{thm:vAbelian_coarse_median}]
Hagen \cite[Theorem~B]{Hag14crystal} proved
\eqref{item:hyperoctahedral}$\iff$\eqref{item:cubulable}
when $G$ is crystallographic. We argue that the virtually Abelian case
of \eqref{item:hyperoctahedral}$\implies$\eqref{item:cubulable}
reduces to the crystallographic case.
  As above, let $A$ be a finite index normal $\mathbb{Z}^n$ subgroup of
  $G$, let $V:=A\otimes_{\mathbb{Z}}\mathbb{R}\cong \mathbb{R}^n$, and
  let $\rho\from G\to\GL(V)$ be the representation induced by
  conjugation on $A$, with $F:=\rho(G)$.
  Let $Q:=G/A$, which is a finite group.

A classical construction of Zassenhaus relates virtually
$\mathbb{Z}^n$ and crystallographic groups.
This is well-known, compare   \cite[Section~3.2]{Hod23} or
\cite[Section~2.13]{AzuHag}, but we will re-use the following notation in a
later part of the proof. 
The construction in the proof of \cite[Satz~1]{Zas48} gives a map $\tau\from G\to V$ and a homomorphism $\alpha\from G\to V\rtimes F$ such that:
\[\alpha:g\mapsto (\tau(g),\rho(g))\]
Furthermore, the restriction of $\alpha$ to $A$ is $\alpha|_A\from
a\mapsto (a,\operatorname{I}_n)$.
Since $\alpha$ is a homomorphism:
\begin{equation}
  \label{eq:2}
  \tau(gh)=\tau(g)+\rho(g)(\tau(h))
\end{equation}
Together with the $A$--normalization, this gives, for $a\in A$ and $h\in G$:
\begin{equation}
  \label{eq:5}
  \tau(ah)=a+\tau(h)
\end{equation}

Define an affine action of $G$ on $V$ by:
\[g.x:=\tau(g)+\rho(g)(x)\]
By averaging over the finite $F$--action, we can 
construct an $F$--invariant inner product on $V$ that makes this
action isometric.
Observe that $A$ is a uniform lattice of translations.
Since $A$ has finite index in $G$, the $G$--action is geometric.
The kernel $N$ of the $G$--action is finite, since it intersects each
of the finitely many cosets of $A$ in at most one element.
Zassenhaus assumes that the conjugation action of $Q$ on $A$ is faithful
and concludes that the affine action is faithful, so $G$ is
crystallographic, but the same argument applies to say that $G/N$ is
crystallographic.
Moreover, the point group of $G/N$ is again $F$, acting on $V$ by the
original representation $\rho$. 
Hence,  if $G$ is hyperoctahedral then so is $G/N$.
In that case,  \cite[Theorem~B]{Hag14crystal} establishes that $G/N$ acts
geometrically by cubical automorphisms on
  a CAT(0) cube complex.
  Since $N$ is finite, there is a geometric, cubical $G$--action on the same
  cube complex, by factoring through $G/N$.
This proves \eqref{item:hyperoctahedral}$\implies$\eqref{item:cubulable}.
  
  Pulling back the cube complex median gives a cubical coarse median
  on $G$. 
  This proves \eqref{item:cubulable}$\implies$\eqref{item:cubicalcoarsemedian}.

  An orbit map is equivariant, and cubical automorphisms preserve the
  standard cubical median.
  Pulling that median back by a quasiinverse to the orbit map is
  therefore coarsely equivariant, proving
  \eqref{item:cubicalcoarsemedian}$\implies$\eqref{item:equivariantcoarsemedian}.

  Now suppose $\mu$ is a coarsely equivariant coarse median on $G$.
  Fix a nonprincipal ultrafilter $\omega$ and scaling sequence
  $r_i\to\infty$.
Take the constant sequence at the identity element $1_G$ of $G$ to be the
observation sequence.
There is an extended pseudometric on sequences $(g_i)$ and $(h_i)$ in
$G$ defined by $\lim_\omega d_G(g_i,h_i)/r_i$. The asymptotic cone of
$G$ with respect to $(\omega,(r_i),(1_G))$ is
the metric quotient of the sequences at finite distance from
$(1_G)$.
Since $G$ is virtually $\mathbb{Z}^n$, this
cone is biLipschitz equivalent to $\mathbb{R}^n$.
Moreover, we can make an explicit biLipschitz equivalence to $V$ using
the maps $\alpha$, $\tau$, and $\rho$ from the first part of the proof.
First observe that $g.0=\tau(g)$, so $\tau\from G\to V$ is the orbit map
based at $0\in V$. Since the action is geometric, $\tau$ is a
quasiisometry.
Define a map from the asymptotic cone to $V$ by:
\[ [(g_i)]\mapsto\lim_\omega\frac{\tau(g_i)}{r_i}\]
The scaling kills off the additive quasiisometry constant of $\tau$,
so this defines a biLipschitz map from the cone to its image.
It is surjective because for every $v\in V$ and $r_i$ there is an
$a_i\in A$ at uniformly bounded distance from $r_iv$, so $a_i/r_i\to v$, which
implies that the point $[(a_i)]$ of the cone maps to $v$.
By \cite[Proposition~9.1]{Bow13b}, $\mu$ induces a Lipschitz median
  $\mu^\omega$ on the cone, and thus on $V$.

To apply \fullref{lem:linear_median_obstruction},  we must verify that
$\mu^\omega$  is equivariant with respect to translations and $F$.
  Given $v\in V$, choose $a_i\in A$ with $a_i/r_i\to v$.
  Since $\tau(a_i g_i)=a_i+\tau(g_i)$, by \eqref{eq:5}, the sequence
  of left translations by $a_i$ induces translation by $v$ on the cone.
  The median equivariance error for $a_i$ is the constant $C$ of
  \ref{NWZ6}, so the median equivariance error for $v$ is $\lim_\omega C/r_i=0$.
  Thus, $\mu^\omega$ is translation equivariant.
  For $g\in G$, by \eqref{eq:2}:
  \[\frac{\tau(gg_i)}{r_i}=\frac{\tau(g)}{r_i}+\rho(g)\left(\frac{\tau(g_i)}{r_i}\right)\]
  Thus, left multiplication by $g$ induces the linear map $\rho(g)$ on
  the cone.
The same error estimate as in the translation case implies that
$\mu^\omega$ is $\rho(g)$--equivariant.
Since $F=\rho(G)$, conclude that $\mu^\omega$ is $F$--equivariant.
  Now \fullref{lem:linear_median_obstruction} says that $F$ is
  conjugate into $\orthog_n(\mathbb{Z})$, 
so $G$ is hyperoctahedral.
  This proves
  \eqref{item:equivariantcoarsemedian}$\implies$\eqref{item:hyperoctahedral}.
\end{proof}

\begin{corollary}\label{cor:triangle_no_equivariant_median}
  The irreducible affine Coxeter groups of types $\widetilde{A}_2$,
  $\widetilde{A}_n$ for $n\geq 4$, $\widetilde{E}_6$, $\widetilde{E}_7$,
  $\widetilde{E}_8$, $\widetilde{F}_4$, and $\widetilde{G}_2$ admit no
  coarsely equivariant coarse medians. 
\end{corollary}
\begin{proof}
  Apply \fullref{thm:vAbelian_coarse_median} and \fullref{thm:hyperoct_classification}.
\end{proof}

There are other results in the literature building off of Hagen's
equivalence between hyperoctahedrality and cocompact cubulation
that give additional equivalent
conditions, either for crystallographic groups or more generally for
virtually Abelian groups. These include a result of Petyt and Spriano
\cite[Theorem~4.4]{PetSpr23}, who add ``being an HHG'' for
crystallographic groups, and Hoda \cite[Theorem~7.1]{Hod23}, who adds
``being Helly'' and ``being coarsely injective'' for virtually
Abelian groups.
\fullref{thm:vAbelian_coarse_median} adds to these results the
condition of ``admitting a coarsely equivariant
coarse median''.
Hoda's proof also runs through the asymptotic cone, so
in that sense our arguments are similar. 

We note an imprecision in \cite[Theorem~4.4]{PetSpr23}: their goal is
to show that for crystallographic groups being an HHG, being
cocompactly cubulated, and being hyperoctahedral are equivalent, but
in place of hyperoctahedrality they use that ``the point group  embeds into the
hyperoctahedral group $\orthog_n(\mathbb{Z})$''.
This is strictly weaker than hyperoctahedrality, recall
\fullref{rem:conjugation}, so their invocation of
\cite[Theorem~B]{Hag14crystal} is not fully justified.
\cite[Corollary~7.9]{AzuHag} closes this gap.
Alternatively, one can complete their proof by observing that an HHG
admits a coarsely equivariant coarse median and applying \fullref{thm:vAbelian_coarse_median}.
Establishing the required linear conjugacy, rather than just an
abstract embedding of the point group, is precisely the purpose of the
differentiation and rescaling argument in the proof of \fullref{lem:linear_median_obstruction}.


\bibliographystyle{hypersshort}
\bibliography{CoxHHS}

\providecommand{\bysame}{\leavevmode\hbox to3em{\hrulefill}\thinspace}
\providecommand{\MR}{\relax\ifhmode\unskip\space\fi MR }
\providecommand{\MRhref}[2]{%
  \href{http://www.ams.org/mathscinet-getitem?mr=#1}{#2}
}
\providecommand{\href}[2]{#2}
\providecommand{\doi}[1]{doi: #1}
\begin{thebibliography}{10}

\bibitem{AbbBehDur21}
C.~Abbott, J.~Behrstock, and M.~G. Durham,
  \emph{\href{http://dx.doi.org/10.1090/btran/50}{Largest acylindrical actions
  and stability in hierarchically hyperbolic groups}}, Trans. Amer. Math. Soc.
  Ser. B \textbf{8} (2021), 66--104, With an appendix by Daniel Berlyne and
  Jacob Russell.

\bibitem{HHSProblemList}
C.~Abbott, J.~Behrstock, and J.~Russell,
  \emph{\href{https://www.wescac.net/HHG_Question_Session-dec-2024.pdf}{Advances
  in hierarchical hyperbolicity problem list}}, 2024, Compiled by the
  participants of the BIRS workshop \emph{Advances in Hierarchical
  Hyperbolicity} (24w5254), Banff International Research Station; version dated
  November 26, 2024,
  \href{https://www.wescac.net/HHG_Question_Session-dec-2024.pdf}{{\texttt{\detokenize{https://www.wescac.net/HHG_Question_Session-dec-2024.pdf}}}}.

\bibitem{AbbBehRus25}
C.~Abbott, J.~Behrstock, and J.~Russell,
  \emph{\href{http://dx.doi.org/10.1142/S1793525323500516}{Structure invariant
  properties of the hierarchically hyperbolic boundary}}, J. Topol. Anal.
  \textbf{17} (2025), no.~3, 859--903.

\bibitem{AzuHag}
P.~Azuelos and M.~Hagen, \emph{Periodic quasiflats in hierarchically hyperbolic
  spaces}, preprint,
  \href{http://arXiv.org/abs/2608.01513v1}{{\texttt{arXiv:2608.01513v1}}}.

\bibitem{BehHagSis17}
J.~Behrstock, M.~F. Hagen, and A.~Sisto,
  \emph{\href{http://dx.doi.org/10.2140/gt.2017.21.1731}{Hierarchically
  hyperbolic spaces, {I}: {C}urve complexes for cubical groups}}, Geom. Topol.
  \textbf{21} (2017), no.~3, 1731--1804.

\bibitem{BehHagSis19}
J.~Behrstock, M.~F. Hagen, and A.~Sisto,
  \emph{\href{http://dx.doi.org/10.2140/pjm.2019.299.257}{Hierarchically
  hyperbolic spaces {II}: {C}ombination theorems and the distance formula}},
  Pacific J. Math. \textbf{299} (2019), no.~2, 257--338.

\bibitem{BehHagSis21}
J.~Behrstock, M.~F. Hagen, and A.~Sisto,
  \emph{\href{http://dx.doi.org/10.1215/00127094-2020-0056}{Quasiflats in
  hierarchically hyperbolic spaces}}, Duke Math. J. \textbf{170} (2021), no.~5,
  909--996.

\bibitem{BerRus22}
D.~Berlyne and J.~Russell,
  \emph{\href{http://dx.doi.org/10.4171/ggd/652}{Hierarchical hyperbolicity of
  graph products}}, Groups Geom. Dyn. \textbf{16} (2022), no.~2, 523--580.

\bibitem{BesBroSag}
M.~Bestvina, K.~Bromberg, and M.~Sageev, \emph{Topological median algebra
  structures on {ER} homology manifolds {I}: local cubulation},
  \href{http://arXiv.org/abs/2512.09875}{{\texttt{arXiv:2512.09875}}}.

\bibitem{Bou02}
N.~Bourbaki, \emph{\href{http://dx.doi.org/10.1007/978-3-540-89394-3}{Lie
  groups and {L}ie algebras. {C}hapters 4--6}}, Elements of Mathematics
  (Berlin), Springer-Verlag, Berlin, 2002, Translated from the 1968 French
  original by Andrew Pressley.

\bibitem{Bow13b}
B.~H. Bowditch, \emph{\href{http://dx.doi.org/10.2140/pjm.2013.261.53}{Coarse
  median spaces and groups}}, Pacific J. Math. \textbf{261} (2013), no.~1,
  53--93.

\bibitem{Bow20}
B.~H. Bowditch, \emph{\href{http://dx.doi.org/10.1017/s0305004118000555}{Median
  and injective metric spaces}}, Math. Proc. Cambridge Philos. Soc.
  \textbf{168} (2020), no.~1, 43--55.

\bibitem{CapFuj10}
P.-E. Caprace and K.~Fujiwara,
  \emph{\href{http://dx.doi.org/10.1007/s00039-009-0042-2}{Rank-one isometries
  of buildings and quasi-morphisms of {K}ac-{M}oody groups}}, Geom. Funct.
  Anal. \textbf{19} (2010), no.~5, 1296--1319.

\bibitem{CapMar13}
P.-E. Caprace and T.~Marquis,
  \emph{\href{http://dx.doi.org/10.1007/s00209-012-1070-4}{Open subgroups of
  locally compact {K}ac-{M}oody groups}}, Math. Z. \textbf{274} (2013),
  no.~1-2, 291--313.

\bibitem{CapMarRei22}
P.-E. Caprace, T.~Marquis, and C.~D. Reid,
  \emph{\href{http://dx.doi.org/10.17879/21089688074}{Locally normal subgroups
  and ends of locally compact {K}ac-{M}oody groups}}, M\"{u}nster J. Math.
  \textbf{15} (2022), no.~2, 473--498.

\bibitem{CapMuh05}
P.-E. Caprace and B.~M\"{u}hlherr,
  \emph{\href{http://dx.doi.org/10.1515/jgth.2005.8.4.467}{Reflection triangles
  in {C}oxeter groups and biautomaticity}}, J. Group Theory \textbf{8} (2005),
  no.~4, 467--489.

\bibitem{CasChuTao}
C.~H. Cashen, M.~Chu, and J.~Tao, \emph{Hyperbolic models and stability
  recognition for {C}oxeter groups}, preprint,
  \href{http://arXiv.org/abs/2608.17036v2}{{\texttt{arXiv:2608.17036v2}}}.

\bibitem{davisbook}
M.~W. Davis, \emph{The geometry and topology of {C}oxeter groups}, London
  Mathematical Society Monographs Series, vol.~32, Princeton University Press,
  Princeton, NJ, 2008.

\bibitem{DruKap18}
C.~Dru\c{t}u and M.~Kapovich, \emph{Geometric group theory}, American
  Mathematical Society Colloquium Publications, vol.~63, American Mathematical
  Society, Providence, RI, 2018, With an appendix by Bogdan Nica.

\bibitem{FioLevSag24}
E.~Fioravanti, I.~Levcovitz, and M.~Sageev,
  \emph{\href{http://dx.doi.org/10.1112/topo.12353}{Coarse cubical rigidity}},
  J. Topol. \textbf{17} (2024), no.~3, Paper No. e12353, 50.

\bibitem{Hae16}
T.~Haettel, \emph{\href{http://dx.doi.org/10.2140/agt.2016.16.2895}{Higher rank
  lattices are not coarse median}}, Algebr. Geom. Topol. \textbf{16} (2016),
  no.~5, 2895--2910.

\bibitem{HaeHodPet23}
T.~Haettel, N.~Hoda, and H.~Petyt,
  \emph{\href{http://dx.doi.org/10.2140/gt.2023.27.1587}{Coarse injectivity,
  hierarchical hyperbolicity and semihyperbolicity}}, Geom. Topol. \textbf{27}
  (2023), no.~4, 1587--1633.

\bibitem{Hag14crystal}
M.~F. Hagen, \emph{\href{http://dx.doi.org/10.1112/jlms/jdu017}{Cocompactly
  cubulated crystallographic groups}}, J. Lond. Math. Soc. (2) \textbf{90}
  (2014), no.~1, 140--166.

\bibitem{HagWis10}
F.~Haglund and D.~T. Wise,
  \emph{\href{http://dx.doi.org/10.1016/j.aim.2010.01.011}{Coxeter groups are
  virtually special}}, Adv. Math. \textbf{224} (2010), no.~5, 1890--1903.

\bibitem{Hod23}
N.~Hoda, \emph{\href{http://dx.doi.org/10.1112/blms.12906}{Crystallographic
  {H}elly groups}}, Bull. Lond. Math. Soc. \textbf{55} (2023), no.~6,
  2991--3011.

\bibitem{Hum90}
J.~E. Humphreys,
  \emph{\href{http://dx.doi.org/10.1017/CBO9780511623646}{Reflection groups and
  {C}oxeter groups}}, Cambridge Studies in Advanced Mathematics, vol.~29,
  Cambridge University Press, Cambridge, 1990.

\bibitem{KapRaf14}
I.~Kapovich and K.~Rafi, \emph{\href{http://dx.doi.org/10.4171/GGD/231}{On
  hyperbolicity of free splitting and free factor complexes}}, Groups Geom.
  Dyn. \textbf{8} (2014), no.~2, 391--414.

\bibitem{Kra09}
D.~Krammer, \emph{\href{http://dx.doi.org/10.4171/GGD/52}{The conjugacy problem
  for {C}oxeter groups}}, Groups Geom. Dyn. \textbf{3} (2009), no.~1, 71--171.

\bibitem{MasMin99}
H.~A. Masur and Y.~N. Minsky,
  \emph{\href{http://dx.doi.org/10.1007/s002220050343}{Geometry of the complex
  of curves. {I}. {H}yperbolicity}}, Invent. Math. \textbf{138} (1999), no.~1,
  103--149.

\bibitem{MasMin00}
H.~A. Masur and Y.~N. Minsky,
  \emph{\href{http://dx.doi.org/10.1007/PL00001643}{Geometry of the complex of
  curves. {II}. {H}ierarchical structure}}, Geom. Funct. Anal. \textbf{10}
  (2000), no.~4, 902--974.

\bibitem{Mou88}
G.~Moussong,
  \emph{\href{http://gateway.proquest.com/openurl?url_ver=Z39.88-2004&rft_val_fmt=info:ofi/fmt:kev:mtx:dissertation&res_dat=xri:pqdiss&rft_dat=xri:pqdiss:8824577}{Hyperbolic
  {C}oxeter groups}}, Ph.D. thesis, Ohio State University, 1988, p.~55.

\bibitem{MunPet}
Z.~Munro and H.~Petyt, \emph{Coarse obstructions to cocompact cubulation},
  preprint,
  \href{http://arXiv.org/abs/2407.09275}{{\texttt{arXiv:2407.09275}}}.

\bibitem{NibRee03}
G.~A. Niblo and L.~D. Reeves,
  \emph{\href{http://dx.doi.org/10.1515/jgth.2003.028}{Coxeter groups act on
  {${\rm CAT}(0)$} cube complexes}}, J. Group Theory \textbf{6} (2003), no.~3,
  399--413.

\bibitem{NibWriZha19}
G.~A. Niblo, N.~Wright, and J.~Zhang,
  \emph{\href{http://dx.doi.org/10.4171/GGD/510}{A four point characterisation
  for coarse median spaces}}, Groups Geom. Dyn. \textbf{13} (2019), no.~3,
  939--980.

\bibitem{Par07}
L.~Paris, \emph{\href{http://dx.doi.org/10.1142/S0218196707003779}{Irreducible
  {C}oxeter groups}}, Internat. J. Algebra Comput. \textbf{17} (2007), no.~3,
  427--447.

\bibitem{Per26}
J.~Perlmutter, \emph{The {M}orse local-to-global property for graph products},
  preprint,
  \href{http://arXiv.org/abs/2601.09901}{{\texttt{arXiv:2601.09901}}}.

\bibitem{Pet26}
H.~Petyt, \emph{Mapping class groups are quasicubical}, Amer. J. Math.
  \textbf{148} (2026), no.~3, 771--806.

\bibitem{PetSpr23}
H.~Petyt and D.~Spriano,
  \emph{\href{http://dx.doi.org/10.4171/ggd/706}{Unbounded domains in
  hierarchically hyperbolic groups}}, Groups Geom. Dyn. \textbf{17} (2023),
  no.~2, 479--500.

\bibitem{Spe09}
D.~E. Speyer,
  \emph{\href{http://dx.doi.org/10.1090/S0002-9939-08-09638-X}{Powers of
  {C}oxeter elements in infinite groups are reduced}}, Proc. Amer. Math. Soc.
  \textbf{137} (2009), no.~4, 1295--1302.

\bibitem{Zas48}
H.~Zassenhaus, \emph{\href{http://dx.doi.org/10.1007/BF02568029}{{\"U}ber einen
  {A}lgorithmus zur {B}estimmung der {R}aumgruppen}}, Comment. Math. Helv.
  \textbf{21} (1948), 117--141.

\end{thebibliography}
\end{document}